
\ifx\shlhetal\undefinedcontrolsequence\let\shlhetal\relax\fi
\def\fmtname{AmS-TeX}

\def\fmtversion{2.1}
\catcode`\@=11
\ifx\amstexloaded@\relax\catcode`\@=\active
  \endinput\else\let\amstexloaded@\relax\fi
\newlinechar=`\^^J
\def\W@{\immediate\write\sixt@@n}
\def\CR@{\W@{^^J\fmtname - Version \fmtversion^^J%
COPYRIGHT 1985, 1990, 1991 - AMERICAN MATHEMATICAL SOCIETY^^J%
Use of this macro package is not restricted provided^^J%
each use is acknowledged upon publication.^^J}}
\CR@ \everyjob{\CR@}
\message{Loading definitions for}
\message{misc utility macros,}
\toksdef\toks@@=2
\long\def\rightappend@#1\to#2{\toks@{\\{#1}}\toks@@
 =\expandafter{#2}\xdef#2{\the\toks@@\the\toks@}\toks@{}\toks@@{}}
\def\alloclist@{}
\newif\ifalloc@
\def\showallocations{{\def\\{\immediate\write\m@ne}\alloclist@}\alloc@true}
\def\alloc@#1#2#3#4#5{\global\advance\count1#1by\@ne
 \ch@ck#1#4#2\allocationnumber=\count1#1
 \global#3#5=\allocationnumber
 \edef\next@{\string#5=\string#2\the\allocationnumber}%
 \expandafter\rightappend@\next@\to\alloclist@}
\newcount\count@@
\newcount\count@@@
\def\FN@{\futurelet\next}
\def\DN@{\def\next@}
\def\DNii@{\def\nextii@}
\def\RIfM@{\relax\ifmmode}
\def\RIfMIfI@{\relax\ifmmode\ifinner}
\def\setboxz@h{\setbox\z@\hbox}
\def\wdz@{\wd\z@}
\def\boxz@{\box\z@}
\def\setbox@ne{\setbox\@ne}
\def\wd@ne{\wd\@ne}
\def\iterate{\body\expandafter\iterate\else\fi}
\def\err@#1{\errmessage{AmS-TeX error: #1}}
\newhelp\defaulthelp@{Sorry, I already gave what help I could...^^J
Maybe you should try asking a human?^^J
An error might have occurred before I noticed any problems.^^J
``If all else fails, read the instructions.''}
\def\Err@{\errhelp\defaulthelp@\err@}
\def\eat@#1{}
\def\in@#1#2{\def\in@@##1#1##2##3\in@@{\ifx\in@##2\in@false\else\in@true\fi}%
 \in@@#2#1\in@\in@@}
\newif\ifin@
\def\space@.{\futurelet\space@\relax}
\space@. %
\newhelp\athelp@
{Only certain combinations beginning with @ make sense to me.^^J
Perhaps you wanted \string\@\space for a printed @?^^J
I've ignored the character or group after @.}
{\catcode`\~=\active 
 \lccode`\~=`\@ \lowercase{\gdef~{\FN@\at@}}}
\def\at@{\let\next@\at@@
 \ifcat\noexpand\next a\else\ifcat\noexpand\next0\else
 \ifcat\noexpand\next\relax\else
   \let\next\at@@@\fi\fi\fi
 \next@}
\def\at@@#1{\expandafter
 \ifx\csname\space @\string#1\endcsname\relax
  \expandafter\at@@@ \else
  \csname\space @\string#1\expandafter\endcsname\fi}
\def\at@@@#1{\errhelp\athelp@ \err@{\Invalid@@ @}}
\def\atdef@#1{\expandafter\def\csname\space @\string#1\endcsname}
\newhelp\defahelp@{If you typed \string\define\space cs instead of
\string\define\string\cs\space^^J
I've substituted an inaccessible control sequence so that your^^J
definition will be completed without mixing me up too badly.^^J
If you typed \string\define{\string\cs} the inaccessible control sequence^^J
was defined to be \string\cs, and the rest of your^^J
definition appears as input.}
\newhelp\defbhelp@{I've ignored your definition, because it might^^J
conflict with other uses that are important to me.}
\def\define{\FN@\define@}
\def\define@{\ifcat\noexpand\next\relax
 \expandafter\define@@\else\errhelp\defahelp@                               
 \err@{\string\define\space must be followed by a control
 sequence}\expandafter\def\expandafter\nextii@\fi}                          
\def\undefined@@@@@@@@@@{}
\def\preloaded@@@@@@@@@@{}
\def\next@@@@@@@@@@{}
\def\define@@#1{\ifx#1\relax\errhelp\defbhelp@                              
 \err@{\string#1\space is already defined}\DN@{\DNii@}\else
 \expandafter\ifx\csname\expandafter\eat@\string                            
 #1@@@@@@@@@@\endcsname\undefined@@@@@@@@@@\errhelp\defbhelp@
 \err@{\string#1\space can't be defined}\DN@{\DNii@}\else
 \expandafter\ifx\csname\expandafter\eat@\string#1\endcsname\relax          
 \global\let#1\undefined\DN@{\def#1}\else\errhelp\defbhelp@
 \err@{\string#1\space is already defined}\DN@{\DNii@}\fi
 \fi\fi\next@}

\def\predefine#1#2{\let#1#2}
\def\undefine#1{\let#1\undefined}
\message{page layout,}
\newdimen\captionwidth@
\captionwidth@\hsize
\advance\captionwidth@-1.5in
\def\pagewidth#1{\hsize#1\relax
 \captionwidth@\hsize\advance\captionwidth@-1.5in}
\def\pageheight#1{\vsize#1\relax}
\def\hcorrection#1{\advance\hoffset#1\relax}
\def\vcorrection#1{\advance\voffset#1\relax}
\message{accents/punctuation,}

\let\graveaccent\`
\let\acuteaccent\'
\let\tildeaccent\~
\let\hataccent\^
\let\underscore\_
\let\B\=
\let\D\.
\let\ic@\/
\def\/{\unskip\ic@}
\def\textfonti{\the\textfont\@ne}
\def\t#1#2{{\edef\next@{\the\font}\textfonti\accent"7F \next@#1#2}}
\def~{\unskip\nobreak\ \ignorespaces}
\def\.{.\spacefactor\@m}
\atdef@;{\leavevmode\null;}
\atdef@:{\leavevmode\null:}
\atdef@?{\leavevmode\null?}
\edef\@{\string @}
\def\relaxnext@{\let\next\relax}
\atdef@-{\relaxnext@\leavevmode
 \DN@{\ifx\next-\DN@-{\FN@\nextii@}\else
  \DN@{\leavevmode\hbox{-}}\fi\next@}%
 \DNii@{\ifx\next-\DN@-{\leavevmode\hbox{---}}\else
  \DN@{\leavevmode\hbox{--}}\fi\next@}%
 \FN@\next@}
\def\srdr@{\kern.16667em}
\def\drsr@{\kern.02778em}
\def\sldl@{\drsr@}
\def\dlsl@{\srdr@}
\atdef@"{\unskip\relaxnext@
 \DN@{\ifx\next\space@\DN@. {\FN@\nextii@}\else
  \DN@.{\FN@\nextii@}\fi\next@.}%
 \DNii@{\ifx\next`\DN@`{\FN@\nextiii@}\else
  \ifx\next\lq\DN@\lq{\FN@\nextiii@}\else
  \DN@####1{\FN@\nextiv@}\fi\fi\next@}%
 \def\nextiii@{\ifx\next`\DN@`{\sldl@``}\else\ifx\next\lq
  \DN@\lq{\sldl@``}\else\DN@{\dlsl@`}\fi\fi\next@}%
 \def\nextiv@{\ifx\next'\DN@'{\srdr@''}\else
  \ifx\next\rq\DN@\rq{\srdr@''}\else\DN@{\drsr@'}\fi\fi\next@}%
 \FN@\next@}

\def\textfontii{\the\textfont\tw@}
\def\lbrace@{\delimiter"4266308 }
\def\rbrace@{\delimiter"5267309 }
\def\{{\RIfM@\lbrace@\else{\textfontii f}\spacefactor\@m\fi}
\def\}{\RIfM@\rbrace@\else
 \let\@sf\empty\ifhmode\edef\@sf{\spacefactor\the\spacefactor}\fi
 {\textfontii g}\@sf\relax\fi}
\let\lbrace\{
\let\rbrace\}
\def\AmSTeX{{\textfontii A\kern-.1667em%
  \lower.5ex\hbox{M}\kern-.125emS}-\TeX}
\message{line and page breaks,}
\def\vmodeerr@#1{\Err@{\string#1\space not allowed between paragraphs}}
\def\mathmodeerr@#1{\Err@{\string#1\space not allowed in math mode}}
\def\linebreak{\RIfM@\mathmodeerr@\linebreak\else
 \ifhmode\unskip\unkern\break\else\vmodeerr@\linebreak\fi\fi}

\newskip\saveskip@
\def\allowlinebreak{\RIfM@\mathmodeerr@\allowlinebreak\else
 \ifhmode\saveskip@\lastskip\unskip
 \allowbreak\ifdim\saveskip@>\z@\hskip\saveskip@\fi
 \else\vmodeerr@\allowlinebreak\fi\fi}
\def\nolinebreak{\RIfM@\mathmodeerr@\nolinebreak\else
 \ifhmode\saveskip@\lastskip\unskip
 \nobreak\ifdim\saveskip@>\z@\hskip\saveskip@\fi
 \else\vmodeerr@\nolinebreak\fi\fi}
\def\newline{\relaxnext@
 \DN@{\RIfM@\expandafter\mathmodeerr@\expandafter\newline\else
  \ifhmode\ifx\next\par\else
  \expandafter\unskip\expandafter\null\expandafter\hfill\expandafter\break\fi
  \else
  \expandafter\vmodeerr@\expandafter\newline\fi\fi}%
 \FN@\next@}
\def\dmatherr@#1{\Err@{\string#1\space not allowed in display math mode}}
\def\nondmatherr@#1{\Err@{\string#1\space not allowed in non-display math
 mode}}
\def\onlydmatherr@#1{\Err@{\string#1\space allowed only in display math mode}}
\def\nonmatherr@#1{\Err@{\string#1\space allowed only in math mode}}
\def\mathbreak{\RIfMIfI@\break\else
 \dmatherr@\mathbreak\fi\else\nonmatherr@\mathbreak\fi}
\def\nomathbreak{\RIfMIfI@\nobreak\else
 \dmatherr@\nomathbreak\fi\else\nonmatherr@\nomathbreak\fi}
\def\allowmathbreak{\RIfMIfI@\allowbreak\else
 \dmatherr@\allowmathbreak\fi\else\nonmatherr@\allowmathbreak\fi}
\def\pagebreak{\RIfM@
 \ifinner\nondmatherr@\pagebreak\else\postdisplaypenalty-\@M\fi
 \else\ifvmode\removelastskip\break\else\vadjust{\break}\fi\fi}
\def\nopagebreak{\RIfM@
 \ifinner\nondmatherr@\nopagebreak\else\postdisplaypenalty\@M\fi
 \else\ifvmode\nobreak\else\vadjust{\nobreak}\fi\fi}
\def\nonvmodeerr@#1{\Err@{\string#1\space not allowed within a paragraph
 or in math}}
\def\vnonvmode@#1#2{\relaxnext@\DNii@{\ifx\next\par\DN@{#1}\else
 \DN@{#2}\fi\next@}%
 \ifvmode\DN@{#1}\else
 \DN@{\FN@\nextii@}\fi\next@}
\def\newpage{\vnonvmode@{\vfill\break}{\nonvmodeerr@\newpage}}
\def\smallpagebreak{\vnonvmode@\smallbreak{\nonvmodeerr@\smallpagebreak}}
\def\medpagebreak{\vnonvmode@\medbreak{\nonvmodeerr@\medpagebreak}}
\def\bigpagebreak{\vnonvmode@\bigbreak{\nonvmodeerr@\bigpagebreak}}
\def\NoBlackBoxes{\global\overfullrule\z@}
\def\BlackBoxes{\global\overfullrule5\p@}
\def\Invalid@#1{\def#1{\Err@{\Invalid@@\string#1}}}
\def\Invalid@@{Invalid use of }
\message{figures,}
\Invalid@\caption
\Invalid@\captionwidth
\newdimen\smallcaptionwidth@
\def\topspace{\mid@false\ins@}
\def\midspace{\mid@true\ins@}
\newif\ifmid@
\def\captionfont@{}
\def\ins@#1{\relaxnext@\allowbreak
 \smallcaptionwidth@\captionwidth@\gdef\thespace@{#1}%
 \DN@{\ifx\next\space@\DN@. {\FN@\nextii@}\else
  \DN@.{\FN@\nextii@}\fi\next@.}%
 \DNii@{\ifx\next\caption\DN@\caption{\FN@\nextiii@}%
  \else\let\next@\nextiv@\fi\next@}%
 \def\nextiv@{\vnonvmode@
  {\ifmid@\expandafter\midinsert\else\expandafter\topinsert\fi
   \vbox to\thespace@{}\endinsert}
  {\ifmid@\nonvmodeerr@\midspace\else\nonvmodeerr@\topspace\fi}}%
 \def\nextiii@{\ifx\next\captionwidth\expandafter\nextv@
  \else\expandafter\nextvi@\fi}%
 \def\nextv@\captionwidth##1##2{\smallcaptionwidth@##1\relax\nextvi@{##2}}%
 \def\nextvi@##1{\def\thecaption@{\captionfont@##1}%
  \DN@{\ifx\next\space@\DN@. {\FN@\nextvii@}\else
   \DN@.{\FN@\nextvii@}\fi\next@.}%
  \FN@\next@}%
 \def\nextvii@{\vnonvmode@
  {\ifmid@\expandafter\midinsert\else
  \expandafter\topinsert\fi\vbox to\thespace@{}\nobreak\smallskip
  \setboxz@h{\noindent\ignorespaces\thecaption@\unskip}%
  \ifdim\wdz@>\smallcaptionwidth@\centerline{\vbox{\hsize\smallcaptionwidth@
   \noindent\ignorespaces\thecaption@\unskip}}%
  \else\centerline{\boxz@}\fi\endinsert}
  {\ifmid@\nonvmodeerr@\midspace
  \else\nonvmodeerr@\topspace\fi}}%
 \FN@\next@}
\message{comments,}
\def\newcodes@{\catcode`\\12\catcode`\{12\catcode`\}12\catcode`\#12%
 \catcode`\%12\relax}
\def\oldcodes@{\catcode`\\0\catcode`\{1\catcode`\}2\catcode`\#6%
 \catcode`\%14\relax}
\def\comment{\newcodes@\endlinechar=10 \comment@}
{\lccode`\0=`\\
\lowercase{\gdef\comment@#1^^J{\comment@@#10endcomment\comment@@@}%
\gdef\comment@@#10endcomment{\FN@\comment@@@}%
\gdef\comment@@@#1\comment@@@{\ifx\next\comment@@@\let\next\comment@
 \else\def\next{\oldcodes@\endlinechar=`\^^M\relax}%
 \fi\next}}}
\def\pr@m@s{\ifx'\next\DN@##1{\prim@s}\else\let\next@\egroup\fi\next@}
\def\prime{{\null\prime@\null}}
\mathchardef\prime@="0230
\let\dsize\displaystyle

\let\ssize\scriptstyle

\message{math spacing,}
\def\,{\RIfM@\mskip\thinmuskip\relax\else\kern.16667em\fi}
\def\!{\RIfM@\mskip-\thinmuskip\relax\else\kern-.16667em\fi}
\let\thinspace\,
\let\negthinspace\!
\def\medspace{\RIfM@\mskip\medmuskip\relax\else\kern.222222em\fi}
\def\negmedspace{\RIfM@\mskip-\medmuskip\relax\else\kern-.222222em\fi}
\def\thickspace{\RIfM@\mskip\thickmuskip\relax\else\kern.27777em\fi}
\let\;\thickspace
\def\negthickspace{\RIfM@\mskip-\thickmuskip\relax\else
 \kern-.27777em\fi}
\atdef@,{\RIfM@\mskip.1\thinmuskip\else\leavevmode\null,\fi}
\atdef@!{\RIfM@\mskip-.1\thinmuskip\else\leavevmode\null!\fi}
\atdef@.{\RIfM@&&\else\leavevmode.\spacefactor3000 \fi}
\def\and{\DOTSB\;\mathchar"3026 \;}

\message{fractions,}
\def\frac#1#2{{#1\over#2}}

\newdimen\ex@
\ex@.2326ex
\Invalid@\thickness
\def\thickfrac{\relaxnext@
 \DN@{\ifx\next\thickness\let\next@\nextii@\else
 \DN@{\nextii@\thickness1}\fi\next@}%
 \DNii@\thickness##1##2##3{{##2\above##1\ex@##3}}%
 \FN@\next@}

\def\thickfracwithdelims#1#2{\relaxnext@\def\ldelim@{#1}\def\rdelim@{#2}%
 \DN@{\ifx\next\thickness\let\next@\nextii@\else
 \DN@{\nextii@\thickness1}\fi\next@}%
 \DNii@\thickness##1##2##3{{##2\abovewithdelims
 \ldelim@\rdelim@##1\ex@##3}}%
 \FN@\next@}

\def\:{\nobreak\hskip.1111em\mathpunct{}\nonscript\mkern-\thinmuskip{:}\hskip
 .3333emplus.0555em\relax}
\def\snug{\unskip\kern-\mathsurround}
\message{smash commands,}
\def\topsmash{\top@true\bot@false\smash@}
\def\botsmash{\top@false\bot@true\smash@}
\newif\iftop@
\newif\ifbot@
\def\smash{\top@true\bot@true\smash@}
\def\smash@{\RIfM@\expandafter\mathpalette\expandafter\mathsm@sh\else
 \expandafter\makesm@sh\fi}
\def\finsm@sh{\iftop@\ht\z@\z@\fi\ifbot@\dp\z@\z@\fi\leavevmode\boxz@}
\message{large operator symbols,}
\def\LimitsOnSums{\global\let\slimits@\displaylimits}
\def\NoLimitsOnSums{\global\let\slimits@\nolimits}
\LimitsOnSums
\mathchardef\coprod@="1360       \def\coprod{\DOTSB\coprod@\slimits@}
\mathchardef\bigvee@="1357       \def\bigvee{\DOTSB\bigvee@\slimits@}
\mathchardef\bigwedge@="1356     \def\bigwedge{\DOTSB\bigwedge@\slimits@}
\mathchardef\biguplus@="1355     \def\biguplus{\DOTSB\biguplus@\slimits@}
\mathchardef\bigcap@="1354       \def\bigcap{\DOTSB\bigcap@\slimits@}
\mathchardef\bigcup@="1353       \def\bigcup{\DOTSB\bigcup@\slimits@}
\mathchardef\prod@="1351         \def\prod{\DOTSB\prod@\slimits@}
\mathchardef\sum@="1350          \def\sum{\DOTSB\sum@\slimits@}
\mathchardef\bigotimes@="134E    \def\bigotimes{\DOTSB\bigotimes@\slimits@}
\mathchardef\bigoplus@="134C     \def\bigoplus{\DOTSB\bigoplus@\slimits@}
\mathchardef\bigodot@="134A      \def\bigodot{\DOTSB\bigodot@\slimits@}
\mathchardef\bigsqcup@="1346     \def\bigsqcup{\DOTSB\bigsqcup@\slimits@}
\message{integrals,}
\def\LimitsOnInts{\global\let\ilimits@\displaylimits}
\def\NoLimitsOnInts{\global\let\ilimits@\nolimits}
\NoLimitsOnInts
\def\int{\DOTSI\intop\ilimits@}
\def\oint{\DOTSI\ointop\ilimits@}
\def\intic@{\mathchoice{\hskip.5em}{\hskip.4em}{\hskip.4em}{\hskip.4em}}
\def\negintic@{\mathchoice
 {\hskip-.5em}{\hskip-.4em}{\hskip-.4em}{\hskip-.4em}}
\def\intkern@{\mathchoice{\!\!\!}{\!\!}{\!\!}{\!\!}}
\def\intdots@{\mathchoice{\plaincdots@}
 {{\cdotp}\mkern1.5mu{\cdotp}\mkern1.5mu{\cdotp}}
 {{\cdotp}\mkern1mu{\cdotp}\mkern1mu{\cdotp}}
 {{\cdotp}\mkern1mu{\cdotp}\mkern1mu{\cdotp}}}
\newcount\intno@
\def\iint{\DOTSI\intno@\tw@\FN@\ints@}
\def\iiint{\DOTSI\intno@\thr@@\FN@\ints@}
\def\iiiint{\DOTSI\intno@4 \FN@\ints@}
\def\idotsint{\DOTSI\intno@\z@\FN@\ints@}
\def\ints@{\findlimits@\ints@@}
\newif\iflimtoken@
\newif\iflimits@
\def\findlimits@{\limtoken@true\ifx\next\limits\limits@true
 \else\ifx\next\nolimits\limits@false\else
 \limtoken@false\ifx\ilimits@\nolimits\limits@false\else
 \ifinner\limits@false\else\limits@true\fi\fi\fi\fi}
\def\multint@{\int\ifnum\intno@=\z@\intdots@                                
 \else\intkern@\fi                                                          
 \ifnum\intno@>\tw@\int\intkern@\fi                                         
 \ifnum\intno@>\thr@@\int\intkern@\fi                                       
 \int}                                                                      
\def\multintlimits@{\intop\ifnum\intno@=\z@\intdots@\else\intkern@\fi
 \ifnum\intno@>\tw@\intop\intkern@\fi
 \ifnum\intno@>\thr@@\intop\intkern@\fi\intop}
\def\ints@@{\iflimtoken@                                                    
 \def\ints@@@{\iflimits@\negintic@\mathop{\intic@\multintlimits@}\limits    
  \else\multint@\nolimits\fi                                                
  \eat@}                                                                    
 \else                                                                      
 \def\ints@@@{\iflimits@\negintic@
  \mathop{\intic@\multintlimits@}\limits\else
  \multint@\nolimits\fi}\fi\ints@@@}
\def\LimitsOnNames{\global\let\nlimits@\displaylimits}
\def\NoLimitsOnNames{\global\let\nlimits@\nolimits@}
\LimitsOnNames
\def\nolimits@{\relaxnext@
 \DN@{\ifx\next\limits\DN@\limits{\nolimits}\else
  \let\next@\nolimits\fi\next@}%
 \FN@\next@}
\message{operator names,}
\def\newmcodes@{\mathcode`\'"27\mathcode`\*"2A\mathcode`\."613A%
 \mathcode`\-"2D\mathcode`\/"2F\mathcode`\:"603A }
\def\operatorname#1{\mathop{\newmcodes@\kern\z@\fam\z@#1}\nolimits@}
\def\operatornamewithlimits#1{\mathop{\newmcodes@\kern\z@\fam\z@#1}\nlimits@}
\def\qopname@#1{\mathop{\fam\z@#1}\nolimits@}
\def\qopnamewl@#1{\mathop{\fam\z@#1}\nlimits@}
\def\arccos{\qopname@{arccos}}
\def\arcsin{\qopname@{arcsin}}
\def\arctan{\qopname@{arctan}}
\def\arg{\qopname@{arg}}
\def\cos{\qopname@{cos}}
\def\cosh{\qopname@{cosh}}
\def\cot{\qopname@{cot}}
\def\coth{\qopname@{coth}}
\def\csc{\qopname@{csc}}
\def\deg{\qopname@{deg}}
\def\det{\qopnamewl@{det}}
\def\dim{\qopname@{dim}}
\def\exp{\qopname@{exp}}
\def\gcd{\qopnamewl@{gcd}}
\def\hom{\qopname@{hom}}
\def\inf{\qopnamewl@{inf}}
\def\injlim{\qopnamewl@{inj\,lim}}
\def\ker{\qopname@{ker}}
\def\lg{\qopname@{lg}}
\def\lim{\qopnamewl@{lim}}
\def\liminf{\qopnamewl@{lim\,inf}}
\def\limsup{\qopnamewl@{lim\,sup}}
\def\ln{\qopname@{ln}}
\def\log{\qopname@{log}}
\def\max{\qopnamewl@{max}}
\def\min{\qopnamewl@{min}}
\def\Pr{\qopnamewl@{Pr}}
\def\projlim{\qopnamewl@{proj\,lim}}
\def\sec{\qopname@{sec}}
\def\sin{\qopname@{sin}}
\def\sinh{\qopname@{sinh}}
\def\sup{\qopnamewl@{sup}}
\def\tan{\qopname@{tan}}
\def\tanh{\qopname@{tanh}}
\def\varinjlim{\mathop{\vtop{\ialign{##\crcr
 \hfil\rm lim\hfil\crcr\noalign{\nointerlineskip}\rightarrowfill\crcr
 \noalign{\nointerlineskip\kern-\ex@}\crcr}}}}
\def\varprojlim{\mathop{\vtop{\ialign{##\crcr
 \hfil\rm lim\hfil\crcr\noalign{\nointerlineskip}\leftarrowfill\crcr
 \noalign{\nointerlineskip\kern-\ex@}\crcr}}}}
\def\varliminf{\mathop{\underline{\vrule height\z@ depth.2exwidth\z@
 \hbox{\rm lim}}}}

\newdimen\buffer@
\buffer@\fontdimen13 \tenex
\newdimen\buffer
\buffer\buffer@

\def\ResetBuffer{\fontdimen13 \tenex\buffer@\global\buffer\buffer@}
\def\shave#1{\mathop{\hbox{$\m@th\fontdimen13 \tenex\z@                     
 \displaystyle{#1}$}}\fontdimen13 \tenex\buffer}

\message{multilevel sub/superscripts,}
\Invalid@\\
\def\Let@{\relax\iffalse{\fi\let\\=\cr\iffalse}\fi}
\Invalid@\vspace
\def\vspace@{\def\vspace##1{\crcr\noalign{\vskip##1\relax}}}
\def\multilimits@{\bgroup\vspace@\Let@
 \baselineskip\fontdimen10 \scriptfont\tw@
 \advance\baselineskip\fontdimen12 \scriptfont\tw@
 \lineskip\thr@@\fontdimen8 \scriptfont\thr@@
 \lineskiplimit\lineskip
 \vbox\bgroup\ialign\bgroup\hfil$\m@th\scriptstyle{##}$\hfil\crcr}
\def\Sb{_\multilimits@}
\def\endSb{\crcr\egroup\egroup\egroup}
\def\Sp{^\multilimits@}

\def\spreadlines#1{\RIfMIfI@\onlydmatherr@\spreadlines\else
 \openup#1\relax\fi\else\onlydmatherr@\spreadlines\fi}
\def\Mathstrut@{\copy\Mathstrutbox@}
\newbox\Mathstrutbox@
\setbox\Mathstrutbox@\null
\setboxz@h{$\m@th($}
\ht\Mathstrutbox@\ht\z@
\dp\Mathstrutbox@\dp\z@
\message{matrices,}
\newdimen\spreadmlines@
\def\spreadmatrixlines#1{\RIfMIfI@
 \onlydmatherr@\spreadmatrixlines\else
 \spreadmlines@#1\relax\fi\else\onlydmatherr@\spreadmatrixlines\fi}
\def\matrix{\null\,\vcenter\bgroup\Let@\vspace@
 \normalbaselines\openup\spreadmlines@\ialign
 \bgroup\hfil$\m@th##$\hfil&&\quad\hfil$\m@th##$\hfil\crcr
 \Mathstrut@\crcr\noalign{\kern-\baselineskip}}
\def\endmatrix{\crcr\Mathstrut@\crcr\noalign{\kern-\baselineskip}\egroup
 \egroup\,}
\def\format{\crcr\egroup\iffalse{\fi\ifnum`}=0 \fi\format@}
\newtoks\hashtoks@
\hashtoks@{#}
\def\format@#1\\{\def\preamble@{#1}%
 \def\l{$\m@th\the\hashtoks@$\hfil}%
 \def\c{\hfil$\m@th\the\hashtoks@$\hfil}%
 \def\r{\hfil$\m@th\the\hashtoks@$}%
 \edef\preamble@@{\preamble@}\ifnum`{=0 \fi\iffalse}\fi
 \ialign\bgroup\span\preamble@@\crcr}
\def\smallmatrix{\null\,\vcenter\bgroup\vspace@\Let@
 \baselineskip9\ex@\lineskip\ex@
 \ialign\bgroup\hfil$\m@th\scriptstyle{##}$\hfil&&\thickspace\hfil
 $\m@th\scriptstyle{##}$\hfil\crcr}
\def\endsmallmatrix{\crcr\egroup\egroup\,}

\newmuskip\dotsspace@
\dotsspace@1.5mu
\def\strip@#1 {#1}
\def\spacehdots#1\for#2{\multispan{#2}\xleaders
 \hbox{$\m@th\mkern\strip@#1 \dotsspace@.\mkern\strip@#1 \dotsspace@$}\hfill}
\def\hdotsfor#1{\spacehdots\@ne\for{#1}}
\def\multispan@#1{\omit\mscount#1\unskip\loop\ifnum\mscount>\@ne\sp@n\repeat}
\def\spaceinnerhdots#1\for#2\after#3{\multispan@{\strip@#2 }#3\xleaders
 \hbox{$\m@th\mkern\strip@#1 \dotsspace@.\mkern\strip@#1 \dotsspace@$}\hfill}
\def\innerhdotsfor#1\after#2{\spaceinnerhdots\@ne\for#1\after{#2}}
\def\cases{\bgroup\spreadmlines@\jot\left\{\,\matrix\format\l&\quad\l\\}
\def\endcases{\endmatrix\right.\egroup}
\message{multiline displays,}
\newif\ifinany@
\newif\ifinalign@
\newif\ifingather@
\def\strut@{\copy\strutbox@}
\newbox\strutbox@
\setbox\strutbox@\hbox{\vrule height8\p@ depth3\p@ width\z@}
\def\topaligned{\null\,\vtop\aligned@}
\def\botaligned{\null\,\vbox\aligned@}
\def\aligned{\null\,\vcenter\aligned@}
\def\aligned@{\bgroup\vspace@\Let@
 \ifinany@\else\openup\jot\fi\ialign
 \bgroup\hfil\strut@$\m@th\displaystyle{##}$&
 $\m@th\displaystyle{{}##}$\hfil\crcr}
\def\endaligned{\crcr\egroup\egroup}

\def\alignedat#1{\null\,\vcenter\bgroup\doat@{#1}\vspace@\Let@
 \ifinany@\else\openup\jot\fi\ialign\bgroup\span\preamble@@\crcr}
\newcount\atcount@
\def\doat@#1{\toks@{\hfil\strut@$\m@th
 \displaystyle{\the\hashtoks@}$&$\m@th\displaystyle
 {{}\the\hashtoks@}$\hfil}
 \atcount@#1\relax\advance\atcount@\m@ne                                    
 \loop\ifnum\atcount@>\z@\toks@=\expandafter{\the\toks@&\hfil$\m@th
 \displaystyle{\the\hashtoks@}$&$\m@th
 \displaystyle{{}\the\hashtoks@}$\hfil}\advance
  \atcount@\m@ne\repeat                                                     
 \xdef\preamble@{\the\toks@}\xdef\preamble@@{\preamble@}}

\def\gathered{\null\,\vcenter\bgroup\vspace@\Let@
 \ifinany@\else\openup\jot\fi\ialign
 \bgroup\hfil\strut@$\m@th\displaystyle{##}$\hfil\crcr}
\def\endgathered{\crcr\egroup\egroup}
\newif\iftagsleft@
\def\TagsOnLeft{\global\tagsleft@true}
\def\TagsOnRight{\global\tagsleft@false}
\TagsOnLeft
\newif\ifmathtags@
\def\TagsAsMath{\global\mathtags@true}
\def\TagsAsText{\global\mathtags@false}
\TagsAsText
\def\tagform@#1{\hbox{\rm(\ignorespaces#1\unskip)}}
\def\thetag{\leavevmode\tagform@}
\def\tag#1$${\iftagsleft@\leqno\else\eqno\fi                                
 \maketag@#1\maketag@                                                       
 $$}                                                                        
\def\maketag@{\FN@\maketag@@}
\def\maketag@@{\ifx\next"\expandafter\maketag@@@\else\expandafter\maketag@@@@
 \fi}
\def\maketag@@@"#1"#2\maketag@{\hbox{\rm#1}}                                
\def\maketag@@@@#1\maketag@{\ifmathtags@\tagform@{$\m@th#1$}\else
 \tagform@{#1}\fi}
\interdisplaylinepenalty\@M
\def\allowdisplaybreaks{\RIfMIfI@
 \onlydmatherr@\allowdisplaybreaks\else
 \interdisplaylinepenalty\z@\fi\else\onlydmatherr@\allowdisplaybreaks\fi}
\Invalid@\allowdisplaybreak
\Invalid@\displaybreak
\Invalid@\intertext
\def\allowdisplaybreak@{\def\allowdisplaybreak{\crcr\noalign{\allowbreak}}}
\def\displaybreak@{\def\displaybreak{\crcr\noalign{\break}}}
\def\intertext@{\def\intertext##1{\crcr\noalign{%
 \penalty\postdisplaypenalty \vskip\belowdisplayskip
 \vbox{\normalbaselines\noindent##1}%
 \penalty\predisplaypenalty \vskip\abovedisplayskip}}}
\newskip\centering@
\centering@\z@ plus\@m\p@
\def\align{\relax\ifingather@\DN@{\csname align (in
  \string\gather)\endcsname}\else
 \ifmmode\ifinner\DN@{\onlydmatherr@\align}\else
  \let\next@\align@\fi
 \else\DN@{\onlydmatherr@\align}\fi\fi\next@}
\newhelp\andhelp@
{An extra & here is so disastrous that you should probably exit^^J
and fix things up.}
\newif\iftag@
\newcount\and@
\def\align@{\inalign@true\inany@true
 \vspace@\allowdisplaybreak@\displaybreak@\intertext@
 \def\tag{\global\tag@true\ifnum\and@=\z@\DN@{&&}\else
          \DN@{&}\fi\next@}%
 \iftagsleft@\DN@{\csname align \endcsname}\else
  \DN@{\csname align \space\endcsname}\fi\next@}
\def\Tag@{\iftag@\else\errhelp\andhelp@\err@{Extra & on this line}\fi}
\newdimen\lwidth@
\newdimen\rwidth@
\newdimen\maxlwidth@
\newdimen\maxrwidth@
\newdimen\totwidth@
\def\measure@#1\endalign{\lwidth@\z@\rwidth@\z@\maxlwidth@\z@\maxrwidth@\z@
 \global\and@\z@                                                            
 \setbox@ne\vbox                                                            
  {\everycr{\noalign{\global\tag@false\global\and@\z@}}\Let@                
  \halign{\setboxz@h{$\m@th\displaystyle{\@lign##}$}
   \global\lwidth@\wdz@                                                     
   \ifdim\lwidth@>\maxlwidth@\global\maxlwidth@\lwidth@\fi                  
   \global\advance\and@\@ne                                                 
   &\setboxz@h{$\m@th\displaystyle{{}\@lign##}$}\global\rwidth@\wdz@        
   \ifdim\rwidth@>\maxrwidth@\global\maxrwidth@\rwidth@\fi                  
   \global\advance\and@\@ne                                                
   &\Tag@
   \eat@{##}\crcr#1\crcr}}
 \totwidth@\maxlwidth@\advance\totwidth@\maxrwidth@}                       
\def\displ@y@{\global\dt@ptrue\openup\jot
 \everycr{\noalign{\global\tag@false\global\and@\z@\ifdt@p\global\dt@pfalse
 \vskip-\lineskiplimit\vskip\normallineskiplimit\else
 \penalty\interdisplaylinepenalty\fi}}}
\def\black@#1{\noalign{\ifdim#1>\displaywidth
 \dimen@\prevdepth\nointerlineskip                                          
 \vskip-\ht\strutbox@\vskip-\dp\strutbox@                                   
 \vbox{\noindent\hbox to#1{\strut@\hfill}}
 \prevdepth\dimen@                                                          
 \fi}}
\expandafter\def\csname align \space\endcsname#1\endalign
 {\measure@#1\endalign\global\and@\z@                                       
 \ifingather@\everycr{\noalign{\global\and@\z@}}\else\displ@y@\fi           
 \Let@\tabskip\centering@                                                   
 \halign to\displaywidth
  {\hfil\strut@\setboxz@h{$\m@th\displaystyle{\@lign##}$}
  \global\lwidth@\wdz@\boxz@\global\advance\and@\@ne                        
  \tabskip\z@skip                                                           
  &\setboxz@h{$\m@th\displaystyle{{}\@lign##}$}
  \global\rwidth@\wdz@\boxz@\hfill\global\advance\and@\@ne                  
  \tabskip\centering@                                                       
  &\setboxz@h{\@lign\strut@\maketag@##\maketag@}
  \dimen@\displaywidth\advance\dimen@-\totwidth@
  \divide\dimen@\tw@\advance\dimen@\maxrwidth@\advance\dimen@-\rwidth@     
  \ifdim\dimen@<\tw@\wdz@\llap{\vtop{\normalbaselines\null\boxz@}}
  \else\llap{\boxz@}\fi                                                    
  \tabskip\z@skip                                                          
  \crcr#1\crcr                                                             
  \black@\totwidth@}}                                                      
\newdimen\lineht@
\expandafter\def\csname align \endcsname#1\endalign{\measure@#1\endalign
 \global\and@\z@
 \ifdim\totwidth@>\displaywidth\let\displaywidth@\totwidth@\else
  \let\displaywidth@\displaywidth\fi                                        
 \ifingather@\everycr{\noalign{\global\and@\z@}}\else\displ@y@\fi
 \Let@\tabskip\centering@\halign to\displaywidth
  {\hfil\strut@\setboxz@h{$\m@th\displaystyle{\@lign##}$}%
  \global\lwidth@\wdz@\global\lineht@\ht\z@                                 
  \boxz@\global\advance\and@\@ne
  \tabskip\z@skip&\setboxz@h{$\m@th\displaystyle{{}\@lign##}$}%
  \global\rwidth@\wdz@\ifdim\ht\z@>\lineht@\global\lineht@\ht\z@\fi         
  \boxz@\hfil\global\advance\and@\@ne
  \tabskip\centering@&\kern-\displaywidth@                                  
  \setboxz@h{\@lign\strut@\maketag@##\maketag@}%
  \dimen@\displaywidth\advance\dimen@-\totwidth@
  \divide\dimen@\tw@\advance\dimen@\maxlwidth@\advance\dimen@-\lwidth@
  \ifdim\dimen@<\tw@\wdz@
   \rlap{\vbox{\normalbaselines\boxz@\vbox to\lineht@{}}}\else
   \rlap{\boxz@}\fi
  \tabskip\displaywidth@\crcr#1\crcr\black@\totwidth@}}
\expandafter\def\csname align (in \string\gather)\endcsname
  #1\endalign{\vcenter{\align@#1\endalign}}
\Invalid@\endalign
\newif\ifxat@
\def\alignat{\RIfMIfI@\DN@{\onlydmatherr@\alignat}\else
 \DN@{\csname alignat \endcsname}\fi\else
 \DN@{\onlydmatherr@\alignat}\fi\next@}
\newif\ifmeasuring@
\newbox\savealignat@
\expandafter\def\csname alignat \endcsname#1#2\endalignat                   
 {\inany@true\xat@false
 \def\tag{\global\tag@true\count@#1\relax\multiply\count@\tw@
  \xdef\tag@{}\loop\ifnum\count@>\and@\xdef\tag@{&\tag@}\advance\count@\m@ne
  \repeat\tag@}%
 \vspace@\allowdisplaybreak@\displaybreak@\intertext@
 \displ@y@\measuring@true                                                   
 \setbox\savealignat@\hbox{$\m@th\displaystyle\Let@
  \attag@{#1}
  \vbox{\halign{\span\preamble@@\crcr#2\crcr}}$}%
 \measuring@false                                                           
 \Let@\attag@{#1}
 \tabskip\centering@\halign to\displaywidth
  {\span\preamble@@\crcr#2\crcr                                             
  \black@{\wd\savealignat@}}}                                               
\Invalid@\endalignat
\def\xalignat{\RIfMIfI@
 \DN@{\onlydmatherr@\xalignat}\else
 \DN@{\csname xalignat \endcsname}\fi\else
 \DN@{\onlydmatherr@\xalignat}\fi\next@}
\expandafter\def\csname xalignat \endcsname#1#2\endxalignat
 {\inany@true\xat@true
 \def\tag{\global\tag@true\def\tag@{}\count@#1\relax\multiply\count@\tw@
  \loop\ifnum\count@>\and@\xdef\tag@{&\tag@}\advance\count@\m@ne\repeat\tag@}%
 \vspace@\allowdisplaybreak@\displaybreak@\intertext@
 \displ@y@\measuring@true\setbox\savealignat@\hbox{$\m@th\displaystyle\Let@
 \attag@{#1}\vbox{\halign{\span\preamble@@\crcr#2\crcr}}$}%
 \measuring@false\Let@
 \attag@{#1}\tabskip\centering@\halign to\displaywidth
 {\span\preamble@@\crcr#2\crcr\black@{\wd\savealignat@}}}
\def\attag@#1{\let\Maketag@\maketag@\let\TAG@\Tag@                          
 \let\Tag@=0\let\maketag@=0
 \ifmeasuring@\def\llap@##1{\setboxz@h{##1}\hbox to\tw@\wdz@{}}%
  \def\rlap@##1{\setboxz@h{##1}\hbox to\tw@\wdz@{}}\else
  \let\llap@\llap\let\rlap@\rlap\fi                                         
 \toks@{\hfil\strut@$\m@th\displaystyle{\@lign\the\hashtoks@}$\tabskip\z@skip
  \global\advance\and@\@ne&$\m@th\displaystyle{{}\@lign\the\hashtoks@}$\hfil
  \ifxat@\tabskip\centering@\fi\global\advance\and@\@ne}
 \iftagsleft@
  \toks@@{\tabskip\centering@&\Tag@\kern-\displaywidth
   \rlap@{\@lign\maketag@\the\hashtoks@\maketag@}%
   \global\advance\and@\@ne\tabskip\displaywidth}\else
  \toks@@{\tabskip\centering@&\Tag@\llap@{\@lign\maketag@
   \the\hashtoks@\maketag@}\global\advance\and@\@ne\tabskip\z@skip}\fi      
 \atcount@#1\relax\advance\atcount@\m@ne
 \loop\ifnum\atcount@>\z@
 \toks@=\expandafter{\the\toks@&\hfil$\m@th\displaystyle{\@lign
  \the\hashtoks@}$\global\advance\and@\@ne
  \tabskip\z@skip&$\m@th\displaystyle{{}\@lign\the\hashtoks@}$\hfil\ifxat@
  \tabskip\centering@\fi\global\advance\and@\@ne}\advance\atcount@\m@ne
 \repeat                                                                    
 \xdef\preamble@{\the\toks@\the\toks@@}
 \xdef\preamble@@{\preamble@}
 \let\maketag@\Maketag@\let\Tag@\TAG@}                                      
\Invalid@\endxalignat
\def\xxalignat{\RIfMIfI@
 \DN@{\onlydmatherr@\xxalignat}\else\DN@{\csname xxalignat
  \endcsname}\fi\else
 \DN@{\onlydmatherr@\xxalignat}\fi\next@}
\expandafter\def\csname xxalignat \endcsname#1#2\endxxalignat{\inany@true
 \vspace@\allowdisplaybreak@\displaybreak@\intertext@
 \displ@y\setbox\savealignat@\hbox{$\m@th\displaystyle\Let@
 \xxattag@{#1}\vbox{\halign{\span\preamble@@\crcr#2\crcr}}$}%
 \Let@\xxattag@{#1}\tabskip\z@skip\halign to\displaywidth
 {\span\preamble@@\crcr#2\crcr\black@{\wd\savealignat@}}}
\def\xxattag@#1{\toks@{\tabskip\z@skip\hfil\strut@
 $\m@th\displaystyle{\the\hashtoks@}$&%
 $\m@th\displaystyle{{}\the\hashtoks@}$\hfil\tabskip\centering@&}%
 \atcount@#1\relax\advance\atcount@\m@ne\loop\ifnum\atcount@>\z@
 \toks@=\expandafter{\the\toks@&\hfil$\m@th\displaystyle{\the\hashtoks@}$%
  \tabskip\z@skip&$\m@th\displaystyle{{}\the\hashtoks@}$\hfil
  \tabskip\centering@}\advance\atcount@\m@ne\repeat
 \xdef\preamble@{\the\toks@\tabskip\z@skip}\xdef\preamble@@{\preamble@}}
\Invalid@\endxxalignat
\newdimen\gwidth@
\newdimen\gmaxwidth@
\def\gmeasure@#1\endgather{\gwidth@\z@\gmaxwidth@\z@\setbox@ne\vbox{\Let@
 \halign{\setboxz@h{$\m@th\displaystyle{##}$}\global\gwidth@\wdz@
 \ifdim\gwidth@>\gmaxwidth@\global\gmaxwidth@\gwidth@\fi
 &\eat@{##}\crcr#1\crcr}}}
\def\gather{\RIfMIfI@\DN@{\onlydmatherr@\gather}\else
 \ingather@true\inany@true\def\tag{&}%
 \vspace@\allowdisplaybreak@\displaybreak@\intertext@
 \displ@y\Let@
 \iftagsleft@\DN@{\csname gather \endcsname}\else
  \DN@{\csname gather \space\endcsname}\fi\fi
 \else\DN@{\onlydmatherr@\gather}\fi\next@}
\expandafter\def\csname gather \space\endcsname#1\endgather
 {\gmeasure@#1\endgather\tabskip\centering@
 \halign to\displaywidth{\hfil\strut@\setboxz@h{$\m@th\displaystyle{##}$}%
 \global\gwidth@\wdz@\boxz@\hfil&
 \setboxz@h{\strut@{\maketag@##\maketag@}}%
 \dimen@\displaywidth\advance\dimen@-\gwidth@
 \ifdim\dimen@>\tw@\wdz@\llap{\boxz@}\else
 \llap{\vtop{\normalbaselines\null\boxz@}}\fi
 \tabskip\z@skip\crcr#1\crcr\black@\gmaxwidth@}}
\newdimen\glineht@
\expandafter\def\csname gather \endcsname#1\endgather{\gmeasure@#1\endgather
 \ifdim\gmaxwidth@>\displaywidth\let\gdisplaywidth@\gmaxwidth@\else
 \let\gdisplaywidth@\displaywidth\fi\tabskip\centering@\halign to\displaywidth
 {\hfil\strut@\setboxz@h{$\m@th\displaystyle{##}$}%
 \global\gwidth@\wdz@\global\glineht@\ht\z@\boxz@\hfil&\kern-\gdisplaywidth@
 \setboxz@h{\strut@{\maketag@##\maketag@}}%
 \dimen@\displaywidth\advance\dimen@-\gwidth@
 \ifdim\dimen@>\tw@\wdz@\rlap{\boxz@}\else
 \rlap{\vbox{\normalbaselines\boxz@\vbox to\glineht@{}}}\fi
 \tabskip\gdisplaywidth@\crcr#1\crcr\black@\gmaxwidth@}}
\newif\ifctagsplit@
\def\CenteredTagsOnSplits{\global\ctagsplit@true}
\def\TopOrBottomTagsOnSplits{\global\ctagsplit@false}
\TopOrBottomTagsOnSplits
\def\split{\relax\ifinany@\let\next@\insplit@\else
 \ifmmode\ifinner\def\next@{\onlydmatherr@\split}\else
 \let\next@\outsplit@\fi\else
 \def\next@{\onlydmatherr@\split}\fi\fi\next@}
\def\insplit@{\global\setbox\z@\vbox\bgroup\vspace@\Let@\ialign\bgroup
 \hfil\strut@$\m@th\displaystyle{##}$&$\m@th\displaystyle{{}##}$\hfill\crcr}
\def\endsplit{\crcr\egroup\egroup\iftagsleft@\expandafter\lendsplit@\else
 \expandafter\rendsplit@\fi}
\def\rendsplit@{\global\setbox9 \vbox
 {\unvcopy\z@\global\setbox8 \lastbox\unskip}
 \setbox@ne\hbox{\unhcopy8 \unskip\global\setbox\tw@\lastbox
 \unskip\global\setbox\thr@@\lastbox}
 \global\setbox7 \hbox{\unhbox\tw@\unskip}
 \ifinalign@\ifctagsplit@                                                   
  \gdef\split@{\hbox to\wd\thr@@{}&
   \vcenter{\vbox{\moveleft\wd\thr@@\boxz@}}}
 \else\gdef\split@{&\vbox{\moveleft\wd\thr@@\box9}\crcr
  \box\thr@@&\box7}\fi                                                      
 \else                                                                      
  \ifctagsplit@\gdef\split@{\vcenter{\boxz@}}\else
  \gdef\split@{\box9\crcr\hbox{\box\thr@@\box7}}\fi
 \fi
 \split@}                                                                   
\def\lendsplit@{\global\setbox9\vtop{\unvcopy\z@}
 \setbox@ne\vbox{\unvcopy\z@\global\setbox8\lastbox}
 \setbox@ne\hbox{\unhcopy8\unskip\setbox\tw@\lastbox
  \unskip\global\setbox\thr@@\lastbox}
 \ifinalign@\ifctagsplit@                                                   
  \gdef\split@{\hbox to\wd\thr@@{}&
  \vcenter{\vbox{\moveleft\wd\thr@@\box9}}}
  \else                                                                     
  \gdef\split@{\hbox to\wd\thr@@{}&\vbox{\moveleft\wd\thr@@\box9}}\fi
 \else
  \ifctagsplit@\gdef\split@{\vcenter{\box9}}\else
  \gdef\split@{\box9}\fi
 \fi\split@}
\def\outsplit@#1$${\align\insplit@#1\endalign$$}
\newdimen\multlinegap@
\multlinegap@1em
\newdimen\multlinetaggap@
\multlinetaggap@1em
\def\MultlineGap#1{\global\multlinegap@#1\relax}
\def\multlinegap#1{\RIfMIfI@\onlydmatherr@\multlinegap\else
 \multlinegap@#1\relax\fi\else\onlydmatherr@\multlinegap\fi}
\def\nomultlinegap{\multlinegap{\z@}}
\def\multline{\RIfMIfI@
 \DN@{\onlydmatherr@\multline}\else
 \DN@{\multline@}\fi\else
 \DN@{\onlydmatherr@\multline}\fi\next@}
\newif\iftagin@
\def\tagin@#1{\tagin@false\in@\tag{#1}\ifin@\tagin@true\fi}
\def\multline@#1$${\inany@true\vspace@\allowdisplaybreak@\displaybreak@
 \tagin@{#1}\iftagsleft@\DN@{\multline@l#1$$}\else
 \DN@{\multline@r#1$$}\fi\next@}
\newdimen\mwidth@
\def\rmmeasure@#1\endmultline{%
 \def\shoveleft##1{##1}\def\shoveright##1{##1}
 \setbox@ne\vbox{\Let@\halign{\setboxz@h
  {$\m@th\@lign\displaystyle{}##$}\global\mwidth@\wdz@
  \crcr#1\crcr}}}
\newdimen\mlineht@
\newif\ifzerocr@
\newif\ifonecr@
\def\lmmeasure@#1\endmultline{\global\zerocr@true\global\onecr@false
 \everycr{\noalign{\ifonecr@\global\onecr@false\fi
  \ifzerocr@\global\zerocr@false\global\onecr@true\fi}}
  \def\shoveleft##1{##1}\def\shoveright##1{##1}%
 \setbox@ne\vbox{\Let@\halign{\setboxz@h
  {$\m@th\@lign\displaystyle{}##$}\ifonecr@\global\mwidth@\wdz@
  \global\mlineht@\ht\z@\fi\crcr#1\crcr}}}
\newbox\mtagbox@
\newdimen\ltwidth@
\newdimen\rtwidth@
\def\multline@l#1$${\iftagin@\DN@{\lmultline@@#1$$}\else
 \DN@{\setbox\mtagbox@\null\ltwidth@\z@\rtwidth@\z@
  \lmultline@@@#1$$}\fi\next@}
\def\lmultline@@#1\endmultline\tag#2$${%
 \setbox\mtagbox@\hbox{\maketag@#2\maketag@}
 \lmmeasure@#1\endmultline\dimen@\mwidth@\advance\dimen@\wd\mtagbox@
 \advance\dimen@\multlinetaggap@                                            
 \ifdim\dimen@>\displaywidth\ltwidth@\z@\else\ltwidth@\wd\mtagbox@\fi       
 \lmultline@@@#1\endmultline$$}
\def\lmultline@@@{\displ@y
 \def\shoveright##1{##1\hfilneg\hskip\multlinegap@}%
 \def\shoveleft##1{\setboxz@h{$\m@th\displaystyle{}##1$}%
  \setbox@ne\hbox{$\m@th\displaystyle##1$}%
  \hfilneg
  \iftagin@
   \ifdim\ltwidth@>\z@\hskip\ltwidth@\hskip\multlinetaggap@\fi
  \else\hskip\multlinegap@\fi\hskip.5\wd@ne\hskip-.5\wdz@##1}
  \halign\bgroup\Let@\hbox to\displaywidth
   {\strut@$\m@th\displaystyle\hfil{}##\hfil$}\crcr
   \hfilneg                                                                 
   \iftagin@                                                                
    \ifdim\ltwidth@>\z@                                                     
     \box\mtagbox@\hskip\multlinetaggap@                                    
    \else
     \rlap{\vbox{\normalbaselines\hbox{\strut@\box\mtagbox@}%
     \vbox to\mlineht@{}}}\fi                                               
   \else\hskip\multlinegap@\fi}                                             
\def\multline@r#1$${\iftagin@\DN@{\rmultline@@#1$$}\else
 \DN@{\setbox\mtagbox@\null\ltwidth@\z@\rtwidth@\z@
  \rmultline@@@#1$$}\fi\next@}
\def\rmultline@@#1\endmultline\tag#2$${\ltwidth@\z@
 \setbox\mtagbox@\hbox{\maketag@#2\maketag@}%
 \rmmeasure@#1\endmultline\dimen@\mwidth@\advance\dimen@\wd\mtagbox@
 \advance\dimen@\multlinetaggap@
 \ifdim\dimen@>\displaywidth\rtwidth@\z@\else\rtwidth@\wd\mtagbox@\fi
 \rmultline@@@#1\endmultline$$}
\def\rmultline@@@{\displ@y
 \def\shoveright##1{##1\hfilneg\iftagin@\ifdim\rtwidth@>\z@
  \hskip\rtwidth@\hskip\multlinetaggap@\fi\else\hskip\multlinegap@\fi}%
 \def\shoveleft##1{\setboxz@h{$\m@th\displaystyle{}##1$}%
  \setbox@ne\hbox{$\m@th\displaystyle##1$}%
  \hfilneg\hskip\multlinegap@\hskip.5\wd@ne\hskip-.5\wdz@##1}%
 \halign\bgroup\Let@\hbox to\displaywidth
  {\strut@$\m@th\displaystyle\hfil{}##\hfil$}\crcr
 \hfilneg\hskip\multlinegap@}
\def\endmultline{\iftagsleft@\expandafter\lendmultline@\else
 \expandafter\rendmultline@\fi}
\def\lendmultline@{\hfilneg\hskip\multlinegap@\crcr\egroup}
\def\rendmultline@{\iftagin@                                                
 \ifdim\rtwidth@>\z@                                                        
  \hskip\multlinetaggap@\box\mtagbox@                                       
 \else\llap{\vtop{\normalbaselines\null\hbox{\strut@\box\mtagbox@}}}\fi     
 \else\hskip\multlinegap@\fi                                                
 \hfilneg\crcr\egroup}
\def\bmod{\mskip-\medmuskip\mkern5mu\mathbin{\fam\z@ mod}\penalty900
 \mkern5mu\mskip-\medmuskip}
\def\pmod#1{\allowbreak\ifinner\mkern8mu\else\mkern18mu\fi
 ({\fam\z@ mod}\,\,#1)}
\def\pod#1{\allowbreak\ifinner\mkern8mu\else\mkern18mu\fi(#1)}
\def\mod#1{\allowbreak\ifinner\mkern12mu\else\mkern18mu\fi{\fam\z@ mod}\,\,#1}
\message{continued fractions,}
\newcount\cfraccount@
\def\cfrac{\bgroup\bgroup\advance\cfraccount@\@ne\strut
 \iffalse{\fi\def\\{\over\displaystyle}\iffalse}\fi}
\def\lcfrac{\bgroup\bgroup\advance\cfraccount@\@ne\strut
 \iffalse{\fi\def\\{\hfill\over\displaystyle}\iffalse}\fi}
\def\rcfrac{\bgroup\bgroup\advance\cfraccount@\@ne\strut\hfill
 \iffalse{\fi\def\\{\over\displaystyle}\iffalse}\fi}
\def\gloop@#1\repeat{\gdef\body{#1}\iterate}
\def\endcfrac{\gloop@\ifnum\cfraccount@>\z@\global\advance\cfraccount@\m@ne
 \egroup\hskip-\nulldelimiterspace\egroup\repeat}
\message{compound symbols,}
\def\binrel@#1{\setboxz@h{\thinmuskip0mu
  \medmuskip\m@ne mu\thickmuskip\@ne mu$#1\m@th$}%
 \setbox@ne\hbox{\thinmuskip0mu\medmuskip\m@ne mu\thickmuskip
  \@ne mu${}#1{}\m@th$}%
 \setbox\tw@\hbox{\hskip\wd@ne\hskip-\wdz@}}
\def\overset#1\to#2{\binrel@{#2}\ifdim\wd\tw@<\z@
 \mathbin{\mathop{\kern\z@#2}\limits^{#1}}\else\ifdim\wd\tw@>\z@
 \mathrel{\mathop{\kern\z@#2}\limits^{#1}}\else
 {\mathop{\kern\z@#2}\limits^{#1}}{}\fi\fi}
\def\underset#1\to#2{\binrel@{#2}\ifdim\wd\tw@<\z@
 \mathbin{\mathop{\kern\z@#2}\limits_{#1}}\else\ifdim\wd\tw@>\z@
 \mathrel{\mathop{\kern\z@#2}\limits_{#1}}\else
 {\mathop{\kern\z@#2}\limits_{#1}}{}\fi\fi}
\def\oversetbrace#1\to#2{\overbrace{#2}^{#1}}
\def\undersetbrace#1\to#2{\underbrace{#2}_{#1}}
\def\sideset#1\and#2\to#3{%
 \setbox@ne\hbox{$\dsize{\vphantom{#3}}#1{#3}\m@th$}%
 \setbox\tw@\hbox{$\dsize{#3}#2\m@th$}%
 \hskip\wd@ne\hskip-\wd\tw@\mathop{\hskip\wd\tw@\hskip-\wd@ne
  {\vphantom{#3}}#1{#3}#2}}
\def\rightarrowfill@#1{\setboxz@h{$#1-\m@th$}\ht\z@\z@
  $#1\m@th\copy\z@\mkern-6mu\cleaders
  \hbox{$#1\mkern-2mu\box\z@\mkern-2mu$}\hfill
  \mkern-6mu\mathord\rightarrow$}
\def\leftarrowfill@#1{\setboxz@h{$#1-\m@th$}\ht\z@\z@
  $#1\m@th\mathord\leftarrow\mkern-6mu\cleaders
  \hbox{$#1\mkern-2mu\copy\z@\mkern-2mu$}\hfill
  \mkern-6mu\box\z@$}
\def\leftrightarrowfill@#1{\setboxz@h{$#1-\m@th$}\ht\z@\z@
  $#1\m@th\mathord\leftarrow\mkern-6mu\cleaders
  \hbox{$#1\mkern-2mu\box\z@\mkern-2mu$}\hfill
  \mkern-6mu\mathord\rightarrow$}
\def\overrightarrow{\mathpalette\overrightarrow@}
\def\overrightarrow@#1#2{\vbox{\ialign{##\crcr\rightarrowfill@#1\crcr
 \noalign{\kern-\ex@\nointerlineskip}$\m@th\hfil#1#2\hfil$\crcr}}}

\def\overleftarrow{\mathpalette\overleftarrow@}
\def\overleftarrow@#1#2{\vbox{\ialign{##\crcr\leftarrowfill@#1\crcr
 \noalign{\kern-\ex@\nointerlineskip}$\m@th\hfil#1#2\hfil$\crcr}}}
\def\overleftrightarrow{\mathpalette\overleftrightarrow@}
\def\overleftrightarrow@#1#2{\vbox{\ialign{##\crcr\leftrightarrowfill@#1\crcr
 \noalign{\kern-\ex@\nointerlineskip}$\m@th\hfil#1#2\hfil$\crcr}}}
\def\underrightarrow{\mathpalette\underrightarrow@}
\def\underrightarrow@#1#2{\vtop{\ialign{##\crcr$\m@th\hfil#1#2\hfil$\crcr
 \noalign{\nointerlineskip}\rightarrowfill@#1\crcr}}}

\def\underleftarrow{\mathpalette\underleftarrow@}
\def\underleftarrow@#1#2{\vtop{\ialign{##\crcr$\m@th\hfil#1#2\hfil$\crcr
 \noalign{\nointerlineskip}\leftarrowfill@#1\crcr}}}
\def\underleftrightarrow{\mathpalette\underleftrightarrow@}
\def\underleftrightarrow@#1#2{\vtop{\ialign{##\crcr$\m@th\hfil#1#2\hfil$\crcr
 \noalign{\nointerlineskip}\leftrightarrowfill@#1\crcr}}}
\message{various kinds of dots,}
\let\DOTSI\relax
\let\DOTSB\relax

\newif\ifmath@
{\uccode`7=`\\ \uccode`8=`m \uccode`9=`a \uccode`0=`t \uccode`!=`h
 \uppercase{\gdef\math@#1#2#3#4#5#6\math@{\global\math@false\ifx 7#1\ifx 8#2%
 \ifx 9#3\ifx 0#4\ifx !#5\xdef\meaning@{#6}\global\math@true\fi\fi\fi\fi\fi}}}
\newif\ifmathch@
{\uccode`7=`c \uccode`8=`h \uccode`9=`\"
 \uppercase{\gdef\mathch@#1#2#3#4#5#6\mathch@{\global\mathch@false
  \ifx 7#1\ifx 8#2\ifx 9#5\global\mathch@true\xdef\meaning@{9#6}\fi\fi\fi}}}
\newcount\classnum@
\def\getmathch@#1.#2\getmathch@{\classnum@#1 \divide\classnum@4096
 \ifcase\number\classnum@\or\or\gdef\thedots@{\dotsb@}\or
 \gdef\thedots@{\dotsb@}\fi}
\newif\ifmathbin@
{\uccode`4=`b \uccode`5=`i \uccode`6=`n
 \uppercase{\gdef\mathbin@#1#2#3{\relaxnext@
  \DNii@##1\mathbin@{\ifx\space@\next\global\mathbin@true\fi}%
 \global\mathbin@false\DN@##1\mathbin@{}%
 \ifx 4#1\ifx 5#2\ifx 6#3\DN@{\FN@\nextii@}\fi\fi\fi\next@}}}
\newif\ifmathrel@
{\uccode`4=`r \uccode`5=`e \uccode`6=`l
 \uppercase{\gdef\mathrel@#1#2#3{\relaxnext@
  \DNii@##1\mathrel@{\ifx\space@\next\global\mathrel@true\fi}%
 \global\mathrel@false\DN@##1\mathrel@{}%
 \ifx 4#1\ifx 5#2\ifx 6#3\DN@{\FN@\nextii@}\fi\fi\fi\next@}}}
\newif\ifmacro@
{\uccode`5=`m \uccode`6=`a \uccode`7=`c
 \uppercase{\gdef\macro@#1#2#3#4\macro@{\global\macro@false
  \ifx 5#1\ifx 6#2\ifx 7#3\global\macro@true
  \xdef\meaning@{\macro@@#4\macro@@}\fi\fi\fi}}}
\def\macro@@#1->#2\macro@@{#2}
\newif\ifDOTS@
\newcount\DOTSCASE@
{\uccode`6=`\\ \uccode`7=`D \uccode`8=`O \uccode`9=`T \uccode`0=`S
 \uppercase{\gdef\DOTS@#1#2#3#4#5{\global\DOTS@false\DN@##1\DOTS@{}%
  \ifx 6#1\ifx 7#2\ifx 8#3\ifx 9#4\ifx 0#5\let\next@\DOTS@@\fi\fi\fi\fi\fi
  \next@}}}
{\uccode`3=`B \uccode`4=`I \uccode`5=`X
 \uppercase{\gdef\DOTS@@#1{\relaxnext@
  \DNii@##1\DOTS@{\ifx\space@\next\global\DOTS@true\fi}%
  \DN@{\FN@\nextii@}%
  \ifx 3#1\global\DOTSCASE@\z@\else
  \ifx 4#1\global\DOTSCASE@\@ne\else
  \ifx 5#1\global\DOTSCASE@\tw@\else\DN@##1\DOTS@{}%
  \fi\fi\fi\next@}}}
\newif\ifnot@
{\uccode`5=`\\ \uccode`6=`n \uccode`7=`o \uccode`8=`t
 \uppercase{\gdef\not@#1#2#3#4{\relaxnext@
  \DNii@##1\not@{\ifx\space@\next\global\not@true\fi}%
 \global\not@false\DN@##1\not@{}%
 \ifx 5#1\ifx 6#2\ifx 7#3\ifx 8#4\DN@{\FN@\nextii@}\fi\fi\fi
 \fi\next@}}}
\newif\ifkeybin@
\def\keybin@{\keybin@true
 \ifx\next+\else\ifx\next=\else\ifx\next<\else\ifx\next>\else\ifx\next-\else
 \ifx\next*\else\ifx\next:\else\keybin@false\fi\fi\fi\fi\fi\fi\fi}
\def\dots{\RIfM@\expandafter\mdots@\else\expandafter\tdots@\fi}
\def\tdots@{\unskip\relaxnext@
 \DN@{$\m@th\mathinner{\ldotp\ldotp\ldotp}\,
   \ifx\next,\,$\else\ifx\next.\,$\else\ifx\next;\,$\else\ifx\next:\,$\else
   \ifx\next?\,$\else\ifx\next!\,$\else$ \fi\fi\fi\fi\fi\fi}%
 \ \FN@\next@}
\def\mdots@{\FN@\mdots@@}
\def\mdots@@{\gdef\thedots@{\dotso@}
 \ifx\next\boldkey\gdef\thedots@\boldkey{\boldkeydots@}\else                
 \ifx\next\boldsymbol\gdef\thedots@\boldsymbol{\boldsymboldots@}\else       
 \ifx,\next\gdef\thedots@{\dotsc}
 \else\ifx\not\next\gdef\thedots@{\dotsb@}
 \else\keybin@
 \ifkeybin@\gdef\thedots@{\dotsb@}
 \else\xdef\meaning@{\meaning\next..........}\xdef\meaning@@{\meaning@}
  \expandafter\math@\meaning@\math@
  \ifmath@
   \expandafter\mathch@\meaning@\mathch@
   \ifmathch@\expandafter\getmathch@\meaning@\getmathch@\fi                 
  \else\expandafter\macro@\meaning@@\macro@                                 
  \ifmacro@                                                                
   \expandafter\not@\meaning@\not@\ifnot@\gdef\thedots@{\dotsb@}
  \else\expandafter\DOTS@\meaning@\DOTS@
  \ifDOTS@
   \ifcase\number\DOTSCASE@\gdef\thedots@{\dotsb@}%
    \or\gdef\thedots@{\dotsi}\else\fi                                      
  \else\expandafter\math@\meaning@\math@                                   
  \ifmath@\expandafter\mathbin@\meaning@\mathbin@
  \ifmathbin@\gdef\thedots@{\dotsb@}
  \else\expandafter\mathrel@\meaning@\mathrel@
  \ifmathrel@\gdef\thedots@{\dotsb@}
  \fi\fi\fi\fi\fi\fi\fi\fi\fi\fi\fi\fi
 \thedots@}
\def\plainldots@{\mathinner{\ldotp\ldotp\ldotp}}
\def\plaincdots@{\mathinner{\cdotp\cdotp\cdotp}}
\def\dotsi{\!\plaincdots@}
\let\dotsb@\plaincdots@
\newif\ifextra@
\newif\ifrightdelim@
\def\rightdelim@{\global\rightdelim@true                                    
 \ifx\next)\else                                                            
 \ifx\next]\else
 \ifx\next\rbrack\else
 \ifx\next\}\else
 \ifx\next\rbrace\else
 \ifx\next\rangle\else
 \ifx\next\rceil\else
 \ifx\next\rfloor\else
 \ifx\next\rgroup\else
 \ifx\next\rmoustache\else
 \ifx\next\right\else
 \ifx\next\bigr\else
 \ifx\next\biggr\else
 \ifx\next\Bigr\else                                                        
 \ifx\next\Biggr\else\global\rightdelim@false
 \fi\fi\fi\fi\fi\fi\fi\fi\fi\fi\fi\fi\fi\fi\fi}
\def\extra@{%
 \global\extra@false\rightdelim@\ifrightdelim@\global\extra@true            
 \else\ifx\next$\global\extra@true                                          
 \else\xdef\meaning@{\meaning\next..........}
 \expandafter\macro@\meaning@\macro@\ifmacro@                               
 \expandafter\DOTS@\meaning@\DOTS@
 \ifDOTS@
 \ifnum\DOTSCASE@=\tw@\global\extra@true                                    
 \fi\fi\fi\fi\fi}
\newif\ifbold@
\def\dotso@{\relaxnext@
 \ifbold@
  \let\next\delayed@
  \DNii@{\extra@\plainldots@\ifextra@\,\fi}%
 \else
  \DNii@{\DN@{\extra@\plainldots@\ifextra@\,\fi}\FN@\next@}%
 \fi
 \nextii@}
\def\extrap@#1{%
 \ifx\next,\DN@{#1\,}\else
 \ifx\next;\DN@{#1\,}\else
 \ifx\next.\DN@{#1\,}\else\extra@
 \ifextra@\DN@{#1\,}\else
 \let\next@#1\fi\fi\fi\fi\next@}
\def\ldots{\DN@{\extrap@\plainldots@}%
 \FN@\next@}
\def\cdots{\DN@{\extrap@\plaincdots@}%
 \FN@\next@}

\def\dotsc{\relaxnext@
 \DN@{\ifx\next;\plainldots@\,\else
  \ifx\next.\plainldots@\,\else\extra@\plainldots@
  \ifextra@\,\fi\fi\fi}%
 \FN@\next@}
\def\cdot{\mathchar"2201 }
\def\longrightarrow{\DOTSB\relbar\joinrel\rightarrow}

\message{special superscripts,}
\def\dddot#1{{\mathop{#1}\limits^{\vbox to-1.4\ex@{\kern-\tw@\ex@
 \hbox{\rm...}\vss}}}}
\def\ddddot#1{{\mathop{#1}\limits^{\vbox to-1.4\ex@{\kern-\tw@\ex@
 \hbox{\rm....}\vss}}}}
\def\sphat{^{\mathchoice{}{}%
 {\,\,\botsmash{\hbox{\lower4\ex@\hbox{$\m@th\widehat{\null}$}}}}%
 {\,\botsmash{\hbox{\lower3\ex@\hbox{$\m@th\hat{\null}$}}}}}}

\def\spacute{^{\!\botsmash{\hbox{\lower\@ne ex\hbox{\'{}}}}}}
\def\spgrave{^{\mathchoice{}{}{}{\!}%
 \botsmash{\hbox{\lower\@ne ex\hbox{\`{}}}}}}
\def\spdot{^{\hbox{\raise\ex@\hbox{\rm.}}}}
\def\spddot{^{\hbox{\raise\ex@\hbox{\rm..}}}}
\def\spdddot{^{\hbox{\raise\ex@\hbox{\rm...}}}}
\def\spddddot{^{\hbox{\raise\ex@\hbox{\rm....}}}}
\def\spbreve{^{\!\botsmash{\hbox{\lower4\ex@\hbox{\u{}}}}}}

\message{\string\text,}
\def\textonlyfont@#1#2{\def#1{\RIfM@
 \Err@{Use \string#1\space only in text}\else#2\fi}}
\textonlyfont@\rm\tenrm
\textonlyfont@\it\tenit
\textonlyfont@\sl\tensl
\textonlyfont@\bf\tenbf
\def\oldnos#1{\RIfM@{\mathcode`\,="013B \fam\@ne#1}\else
 \leavevmode\hbox{$\m@th\mathcode`\,="013B \fam\@ne#1$}\fi}
\def\text{\RIfM@\expandafter\text@\else\expandafter\text@@\fi}
\def\text@@#1{\leavevmode\hbox{#1}}
\def\mathhexbox@#1#2#3{\text{$\m@th\mathchar"#1#2#3$}}
\def\dag{{\mathhexbox@279}}
\def\ddag{{\mathhexbox@27A}}
\def\S{{\mathhexbox@278}}
\def\P{{\mathhexbox@27B}}
\newif\iffirstchoice@
\firstchoice@true
\def\text@#1{\mathchoice
 {\hbox{\everymath{\displaystyle}\def\textfonti{\the\textfont\@ne}%
  \def\textfontii{\the\textfont\tw@}\textdef@@ T#1}}
 {\hbox{\firstchoice@false
  \everymath{\textstyle}\def\textfonti{\the\textfont\@ne}%
  \def\textfontii{\the\textfont\tw@}\textdef@@ T#1}}
 {\hbox{\firstchoice@false
  \everymath{\scriptstyle}\def\textfonti{\the\scriptfont\@ne}%
  \def\textfontii{\the\scriptfont\tw@}\textdef@@ S\rm#1}}
 {\hbox{\firstchoice@false
  \everymath{\scriptscriptstyle}\def\textfonti
  {\the\scriptscriptfont\@ne}%
  \def\textfontii{\the\scriptscriptfont\tw@}\textdef@@ s\rm#1}}}
\def\textdef@@#1{\textdef@#1\rm\textdef@#1\bf\textdef@#1\sl\textdef@#1\it}
\def\rmfam{0}
\def\textdef@#1#2{%
 \DN@{\csname\expandafter\eat@\string#2fam\endcsname}%
 \if S#1\edef#2{\the\scriptfont\next@\relax}%
 \else\if s#1\edef#2{\the\scriptscriptfont\next@\relax}%
 \else\edef#2{\the\textfont\next@\relax}\fi\fi}
\scriptfont\itfam\tenit \scriptscriptfont\itfam\tenit
\scriptfont\slfam\tensl \scriptscriptfont\slfam\tensl
\newif\iftopfolded@
\newif\ifbotfolded@
\def\topfoldedtext{\topfolded@true\botfolded@false\foldedtext@}
\def\botfoldedtext{\botfolded@true\topfolded@false\foldedtext@}
\def\foldedtext{\topfolded@false\botfolded@false\foldedtext@}
\Invalid@\foldedwidth
\def\foldedtext@{\relaxnext@
 \DN@{\ifx\next\foldedwidth\let\next@\nextii@\else
  \DN@{\nextii@\foldedwidth{.3\hsize}}\fi\next@}%
 \DNii@\foldedwidth##1##2{\setbox\z@\vbox
  {\normalbaselines\hsize##1\relax
  \tolerance1600 \noindent\ignorespaces##2}\ifbotfolded@\boxz@\else
  \iftopfolded@\vtop{\unvbox\z@}\else\vcenter{\boxz@}\fi\fi}%
 \FN@\next@}
\message{math font commands,}
\def\bold{\RIfM@\expandafter\bold@\else
 \expandafter\nonmatherr@\expandafter\bold\fi}
\def\bold@#1{{\bold@@{#1}}}
\def\bold@@#1{\fam\bffam\relax#1}
\def\slanted{\RIfM@\expandafter\slanted@\else
 \expandafter\nonmatherr@\expandafter\slanted\fi}
\def\slanted@#1{{\slanted@@{#1}}}
\def\slanted@@#1{\fam\slfam\relax#1}
\def\roman{\RIfM@\expandafter\roman@\else
 \expandafter\nonmatherr@\expandafter\roman\fi}
\def\roman@#1{{\roman@@{#1}}}
\def\roman@@#1{\fam\rmfam\relax#1}
\def\italic{\RIfM@\expandafter\italic@\else
 \expandafter\nonmatherr@\expandafter\italic\fi}
\def\italic@#1{{\italic@@{#1}}}
\def\italic@@#1{\fam\itfam\relax#1}
\def\Cal{\RIfM@\expandafter\Cal@\else
 \expandafter\nonmatherr@\expandafter\Cal\fi}
\def\Cal@#1{{\Cal@@{#1}}}
\def\Cal@@#1{\noaccents@\fam\tw@#1}
\mathchardef\Gamma="0000
\mathchardef\Delta="0001
\mathchardef\Theta="0002
\mathchardef\Lambda="0003
\mathchardef\Xi="0004
\mathchardef\Pi="0005
\mathchardef\Sigma="0006
\mathchardef\Upsilon="0007
\mathchardef\Phi="0008
\mathchardef\Psi="0009
\mathchardef\Omega="000A
\mathchardef\varGamma="0100
\mathchardef\varDelta="0101
\mathchardef\varTheta="0102
\mathchardef\varLambda="0103
\mathchardef\varXi="0104
\mathchardef\varPi="0105
\mathchardef\varSigma="0106
\mathchardef\varUpsilon="0107
\mathchardef\varPhi="0108
\mathchardef\varPsi="0109
\mathchardef\varOmega="010A
\let\alloc@@\alloc@
\def\hexnumber@#1{\ifcase#1 0\or 1\or 2\or 3\or 4\or 5\or 6\or 7\or 8\or
 9\or A\or B\or C\or D\or E\or F\fi}
\def\loadmsam{%
 \font@\tenmsa=msam10
 \font@\sevenmsa=msam7
 \font@\fivemsa=msam5
 \alloc@@8\fam\chardef\sixt@@n\msafam
 \textfont\msafam=\tenmsa
 \scriptfont\msafam=\sevenmsa
 \scriptscriptfont\msafam=\fivemsa
 \edef\next{\hexnumber@\msafam}%
 \mathchardef\dabar@"0\next39
 \edef\dashrightarrow{\mathrel{\dabar@\dabar@\mathchar"0\next4B}}%
 \edef\dashleftarrow{\mathrel{\mathchar"0\next4C\dabar@\dabar@}}%
 \let\dasharrow\dashrightarrow
 \edef\ulcorner{\delimiter"4\next70\next70 }%
 \edef\urcorner{\delimiter"5\next71\next71 }%
 \edef\llcorner{\delimiter"4\next78\next78 }%
 \edef\lrcorner{\delimiter"5\next79\next79 }%
 \edef\yen{{\noexpand\mathhexbox@\next55}}%
 \edef\checkmark{{\noexpand\mathhexbox@\next58}}%
 \edef\circledR{{\noexpand\mathhexbox@\next72}}%
 \edef\maltese{{\noexpand\mathhexbox@\next7A}}%
 \global\let\loadmsam\empty}%
\def\loadmsbm{%
 \font@\tenmsb=msbm10 \font@\sevenmsb=msbm7 \font@\fivemsb=msbm5
 \alloc@@8\fam\chardef\sixt@@n\msbfam
 \textfont\msbfam=\tenmsb
 \scriptfont\msbfam=\sevenmsb \scriptscriptfont\msbfam=\fivemsb
 \global\let\loadmsbm\empty
 }
\def\widehat#1{\ifx\undefined\msbfam \DN@{362}%
  \else \setboxz@h{$\m@th#1$}%
    \edef\next@{\ifdim\wdz@>\tw@ em%
        \hexnumber@\msbfam 5B%
      \else 362\fi}\fi
  \mathaccent"0\next@{#1}}
\def\widetilde#1{\ifx\undefined\msbfam \DN@{365}%
  \else \setboxz@h{$\m@th#1$}%
    \edef\next@{\ifdim\wdz@>\tw@ em%
        \hexnumber@\msbfam 5D%
      \else 365\fi}\fi
  \mathaccent"0\next@{#1}}
\message{\string\newsymbol,}
\def\newsymbol#1#2#3#4#5{\define#1{}%
  \count@#2\relax \advance\count@\m@ne 
 \ifcase\count@
   \ifx\undefined\msafam\loadmsam\fi \let\next@\msafam
 \or \ifx\undefined\msbfam\loadmsbm\fi \let\next@\msbfam
 \else  \Err@{\Invalid@@\string\newsymbol}\let\next@\tw@\fi
 \mathchardef#1="#3\hexnumber@\next@#4#5\space}

\def\Bbb{\RIfM@\expandafter\Bbb@\else
 \expandafter\nonmatherr@\expandafter\Bbb\fi}
\def\Bbb@#1{{\Bbb@@{#1}}}
\def\Bbb@@#1{\noaccents@\fam\msbfam\relax#1}
\message{bold Greek and bold symbols,}
\def\loadbold{%
 \font@\tencmmib=cmmib10 \font@\sevencmmib=cmmib7 \font@\fivecmmib=cmmib5
 \skewchar\tencmmib'177 \skewchar\sevencmmib'177 \skewchar\fivecmmib'177
 \alloc@@8\fam\chardef\sixt@@n\cmmibfam
 \textfont\cmmibfam\tencmmib
 \scriptfont\cmmibfam\sevencmmib \scriptscriptfont\cmmibfam\fivecmmib
 \font@\tencmbsy=cmbsy10 \font@\sevencmbsy=cmbsy7 \font@\fivecmbsy=cmbsy5
 \skewchar\tencmbsy'60 \skewchar\sevencmbsy'60 \skewchar\fivecmbsy'60
 \alloc@@8\fam\chardef\sixt@@n\cmbsyfam
 \textfont\cmbsyfam\tencmbsy
 \scriptfont\cmbsyfam\sevencmbsy \scriptscriptfont\cmbsyfam\fivecmbsy
 \let\loadbold\empty
}
\def\boldnotloaded#1{\Err@{\ifcase#1\or First\else Second\fi
       bold symbol font not loaded}}
\def\mathchari@#1#2#3{\ifx\undefined\cmmibfam
    \boldnotloaded@\@ne
  \else\mathchar"#1\hexnumber@\cmmibfam#2#3\space \fi}
\def\mathcharii@#1#2#3{\ifx\undefined\cmbsyfam
    \boldnotloaded\tw@
  \else \mathchar"#1\hexnumber@\cmbsyfam#2#3\space\fi}
\edef\bffam@{\hexnumber@\bffam}
\def\boldkey#1{\ifcat\noexpand#1A%
  \ifx\undefined\cmmibfam \boldnotloaded\@ne
  \else {\fam\cmmibfam#1}\fi
 \else
 \ifx#1!\mathchar"5\bffam@21 \else
 \ifx#1(\mathchar"4\bffam@28 \else\ifx#1)\mathchar"5\bffam@29 \else
 \ifx#1+\mathchar"2\bffam@2B \else\ifx#1:\mathchar"3\bffam@3A \else
 \ifx#1;\mathchar"6\bffam@3B \else\ifx#1=\mathchar"3\bffam@3D \else
 \ifx#1?\mathchar"5\bffam@3F \else\ifx#1[\mathchar"4\bffam@5B \else
 \ifx#1]\mathchar"5\bffam@5D \else
 \ifx#1,\mathchari@63B \else
 \ifx#1-\mathcharii@200 \else
 \ifx#1.\mathchari@03A \else
 \ifx#1/\mathchari@03D \else
 \ifx#1<\mathchari@33C \else
 \ifx#1>\mathchari@33E \else
 \ifx#1*\mathcharii@203 \else
 \ifx#1|\mathcharii@06A \else
 \ifx#10\bold0\else\ifx#11\bold1\else\ifx#12\bold2\else\ifx#13\bold3\else
 \ifx#14\bold4\else\ifx#15\bold5\else\ifx#16\bold6\else\ifx#17\bold7\else
 \ifx#18\bold8\else\ifx#19\bold9\else
  \Err@{\string\boldkey\space can't be used with #1}%
 \fi\fi\fi\fi\fi\fi\fi\fi\fi\fi\fi\fi\fi\fi\fi
 \fi\fi\fi\fi\fi\fi\fi\fi\fi\fi\fi\fi\fi\fi}
\def\boldsymbol#1{%
 \DN@{\Err@{You can't use \string\boldsymbol\space with \string#1}#1}%
 \ifcat\noexpand#1A%
   \let\next@\relax
   \ifx\undefined\cmmibfam \boldnotloaded\@ne
   \else {\fam\cmmibfam#1}\fi
 \else
  \xdef\meaning@{\meaning#1.........}%
  \expandafter\math@\meaning@\math@
  \ifmath@
   \expandafter\mathch@\meaning@\mathch@
   \ifmathch@
    \expandafter\boldsymbol@@\meaning@\boldsymbol@@
   \fi
  \else
   \expandafter\macro@\meaning@\macro@
   \expandafter\delim@\meaning@\delim@
   \ifdelim@
    \expandafter\delim@@\meaning@\delim@@
   \else
    \boldsymbol@{#1}%
   \fi
  \fi
 \fi
 \next@}
\def\mathhexboxii@#1#2{\ifx\undefined\cmbsyfam
    \boldnotloaded\tw@
  \else \mathhexbox@{\hexnumber@\cmbsyfam}{#1}{#2}\fi}
\def\boldsymbol@#1{\let\next@\relax\let\next#1%
 \ifx\next\cdot\mathcharii@201 \else
 \ifx\next\prime{{\null\mathcharii@030 \null}}\else
 \ifx\next\lbrack\mathchar"4\bffam@5B \else
 \ifx\next\rbrack\mathchar"5\bffam@5D \else
 \ifx\next\{\mathcharii@466 \else
 \ifx\next\lbrace\mathcharii@466 \else
 \ifx\next\}\mathcharii@567 \else
 \ifx\next\rbrace\mathcharii@567 \else
 \ifx\next\surd{{\mathcharii@170}}\else
 \ifx\next\S{{\mathhexboxii@78}}\else
 \ifx\next\P{{\mathhexboxii@7B}}\else
 \ifx\next\dag{{\mathhexboxii@79}}\else
 \ifx\next\ddag{{\mathhexboxii@7A}}\else
 \DN@{\Err@{You can't use \string\boldsymbol\space with \string#1}#1}%
 \fi\fi\fi\fi\fi\fi\fi\fi\fi\fi\fi\fi\fi}
\def\boldsymbol@@#1.#2\boldsymbol@@{\classnum@#1 \count@@@\classnum@        
 \divide\classnum@4096 \count@\classnum@                                    
 \multiply\count@4096 \advance\count@@@-\count@ \count@@\count@@@           
 \divide\count@@@\@cclvi \count@\count@@                                    
 \multiply\count@@@\@cclvi \advance\count@@-\count@@@                       
 \divide\count@@@\@cclvi                                                    
 \multiply\classnum@4096 \advance\classnum@\count@@                         
 \ifnum\count@@@=\z@                                                        
  \count@"\bffam@ \multiply\count@\@cclvi
  \advance\classnum@\count@
  \DN@{\mathchar\number\classnum@}%
 \else
  \ifnum\count@@@=\@ne                                                      
   \ifx\undefined\cmmibfam \DN@{\boldnotloaded\@ne}%
   \else \count@\cmmibfam \multiply\count@\@cclvi
     \advance\classnum@\count@
     \DN@{\mathchar\number\classnum@}\fi
  \else
   \ifnum\count@@@=\tw@                                                    
     \ifx\undefined\cmbsyfam
       \DN@{\boldnotloaded\tw@}%
     \else
       \count@\cmbsyfam \multiply\count@\@cclvi
       \advance\classnum@\count@
       \DN@{\mathchar\number\classnum@}%
     \fi
  \fi
 \fi
\fi}
\newif\ifdelim@
\newcount\delimcount@
{\uccode`6=`\\ \uccode`7=`d \uccode`8=`e \uccode`9=`l
 \uppercase{\gdef\delim@#1#2#3#4#5\delim@
  {\delim@false\ifx 6#1\ifx 7#2\ifx 8#3\ifx 9#4\delim@true
   \xdef\meaning@{#5}\fi\fi\fi\fi}}}
\def\delim@@#1"#2#3#4#5#6\delim@@{\if#32%
\let\next@\relax
 \ifx\undefined\cmbsyfam \boldnotloaded\@ne
 \else \mathcharii@#2#4#5\space \fi\fi}
\def\vert{\delimiter"026A30C }
\def\Vert{\delimiter"026B30D }
\let\|\Vert
\def\backslash{\delimiter"026E30F }
\def\boldkeydots@#1{\bold@true\let\next=#1\let\delayed@=#1\mdots@@
 \boldkey#1\bold@false}  
\def\boldsymboldots@#1{\bold@true\let\next#1\let\delayed@#1\mdots@@
 \boldsymbol#1\bold@false}
\message{Euler fonts,}

\def\frak{\mathfont@\frak}

\def\loadmathfont#1{%
   \expandafter\font@\csname ten#1\endcsname=#110
   \expandafter\font@\csname seven#1\endcsname=#17
   \expandafter\font@\csname five#1\endcsname=#15
   \edef\next{\noexpand\alloc@@8\fam\chardef\sixt@@n
     \expandafter\noexpand\csname#1fam\endcsname}%
   \next
   \textfont\csname#1fam\endcsname \csname ten#1\endcsname
   \scriptfont\csname#1fam\endcsname \csname seven#1\endcsname
   \scriptscriptfont\csname#1fam\endcsname \csname five#1\endcsname
   \expandafter\def\csname #1\expandafter\endcsname\expandafter{%
      \expandafter\mathfont@\csname#1\endcsname}%
 \expandafter\gdef\csname load#1\endcsname{}%
}
\def\mathfont@#1{\RIfM@\expandafter\mathfont@@\expandafter#1\else
  \expandafter\nonmatherr@\expandafter#1\fi}
\def\mathfont@@#1#2{{\mathfont@@@#1{#2}}}
\def\mathfont@@@#1#2{\noaccents@
   \fam\csname\expandafter\eat@\string#1fam\endcsname
   \relax#2}
\message{math accents,}
\def\accentclass@{7}
\def\noaccents@{\def\accentclass@{0}}
\def\makeacc@#1#2{\def#1{\mathaccent"\accentclass@#2 }}
\makeacc@\hat{05E}
\makeacc@\check{014}
\makeacc@\tilde{07E}
\makeacc@\acute{013}
\makeacc@\grave{012}
\makeacc@\dot{05F}
\makeacc@\ddot{07F}
\makeacc@\breve{015}
\makeacc@\bar{016}

\newcount\skewcharcount@
\newcount\familycount@
\def\theskewchar@{\familycount@\@ne
 \global\skewcharcount@\the\skewchar\textfont\@ne                           
 \ifnum\fam>\m@ne\ifnum\fam<16
  \global\familycount@\the\fam\relax
  \global\skewcharcount@\the\skewchar\textfont\the\fam\relax\fi\fi          
 \ifnum\skewcharcount@>\m@ne
  \ifnum\skewcharcount@<128
  \multiply\familycount@256
  \global\advance\skewcharcount@\familycount@
  \global\advance\skewcharcount@28672
  \mathchar\skewcharcount@\else
  \global\skewcharcount@\m@ne\fi\else
 \global\skewcharcount@\m@ne\fi}                                            
\newcount\pointcount@
\def\getpoints@#1.#2\getpoints@{\pointcount@#1 }
\newdimen\accentdimen@
\newcount\accentmu@
\def\dimentomu@{\multiply\accentdimen@ 100
 \expandafter\getpoints@\the\accentdimen@\getpoints@
 \multiply\pointcount@18
 \divide\pointcount@\@m
 \global\accentmu@\pointcount@}
\def\Makeacc@#1#2{\def#1{\RIfM@\DN@{\mathaccent@
 {"\accentclass@#2 }}\else\DN@{\nonmatherr@{#1}}\fi\next@}}
\def\unbracefonts@{\let\Cal@\Cal@@\let\roman@\roman@@\let\bold@\bold@@
 \let\slanted@\slanted@@}
\def\mathaccent@#1#2{\ifnum\fam=\m@ne\xdef\thefam@{1}\else
 \xdef\thefam@{\the\fam}\fi                                                 
 \accentdimen@\z@                                                           
 \setboxz@h{\unbracefonts@$\m@th\fam\thefam@\relax#2$}
 \ifdim\accentdimen@=\z@\DN@{\mathaccent#1{#2}}
  \setbox@ne\hbox{\unbracefonts@$\m@th\fam\thefam@\relax#2\theskewchar@$}
  \setbox\tw@\hbox{$\m@th\ifnum\skewcharcount@=\m@ne\else
   \mathchar\skewcharcount@\fi$}
  \global\accentdimen@\wd@ne\global\advance\accentdimen@-\wdz@
  \global\advance\accentdimen@-\wd\tw@                                     
  \global\multiply\accentdimen@\tw@
  \dimentomu@\global\advance\accentmu@\@ne                                 
 \else\DN@{{\mathaccent#1{#2\mkern\accentmu@ mu}%
    \mkern-\accentmu@ mu}{}}\fi                                             
 \next@}\Makeacc@\Hat{05E}
\Makeacc@\Check{014}
\Makeacc@\Tilde{07E}
\Makeacc@\Acute{013}
\Makeacc@\Grave{012}
\Makeacc@\Dot{05F}
\Makeacc@\Ddot{07F}
\Makeacc@\Breve{015}
\Makeacc@\Bar{016}
\def\Vec{\RIfM@\DN@{\mathaccent@{"017E }}\else
 \DN@{\nonmatherr@\Vec}\fi\next@}
\def\accentedsymbol#1#2{\csname newbox\expandafter\endcsname
  \csname\expandafter\eat@\string#1@box\endcsname
 \expandafter\setbox\csname\expandafter\eat@
  \string#1@box\endcsname\hbox{$\m@th#2$}\define
  #1{\copy\csname\expandafter\eat@\string#1@box\endcsname{}}}
\message{roots,}
\def\sqrt#1{\radical"270370 {#1}}
\let\underline@\underline
\let\overline@\overline
\def\underline#1{\underline@{#1}}
\def\overline#1{\overline@{#1}}
\Invalid@\leftroot
\Invalid@\uproot
\newcount\uproot@
\newcount\leftroot@
\def\root{\relaxnext@
  \DN@{\ifx\next\uproot\let\next@\nextii@\else
   \ifx\next\leftroot\let\next@\nextiii@\else
   \let\next@\plainroot@\fi\fi\next@}%
  \DNii@\uproot##1{\uproot@##1\relax\FN@\nextiv@}%
  \def\nextiv@{\ifx\next\space@\DN@. {\FN@\nextv@}\else
   \DN@.{\FN@\nextv@}\fi\next@.}%
  \def\nextv@{\ifx\next\leftroot\let\next@\nextvi@\else
   \let\next@\plainroot@\fi\next@}%
  \def\nextvi@\leftroot##1{\leftroot@##1\relax\plainroot@}%
   \def\nextiii@\leftroot##1{\leftroot@##1\relax\FN@\nextvii@}%
  \def\nextvii@{\ifx\next\space@
   \DN@. {\FN@\nextviii@}\else
   \DN@.{\FN@\nextviii@}\fi\next@.}%
  \def\nextviii@{\ifx\next\uproot\let\next@\nextix@\else
   \let\next@\plainroot@\fi\next@}%
  \def\nextix@\uproot##1{\uproot@##1\relax\plainroot@}%
  \bgroup\uproot@\z@\leftroot@\z@\FN@\next@}
\def\plainroot@#1\of#2{\setbox\rootbox\hbox{$\m@th\scriptscriptstyle{#1}$}%
 \mathchoice{\r@@t\displaystyle{#2}}{\r@@t\textstyle{#2}}
 {\r@@t\scriptstyle{#2}}{\r@@t\scriptscriptstyle{#2}}\egroup}
\def\r@@t#1#2{\setboxz@h{$\m@th#1\sqrt{#2}$}%
 \dimen@\ht\z@\advance\dimen@-\dp\z@
 \setbox@ne\hbox{$\m@th#1\mskip\uproot@ mu$}\advance\dimen@ 1.667\wd@ne
 \mkern-\leftroot@ mu\mkern5mu\raise.6\dimen@\copy\rootbox
 \mkern-10mu\mkern\leftroot@ mu\boxz@}
\def\boxed#1{\setboxz@h{$\m@th\displaystyle{#1}$}\dimen@.4\ex@
 \advance\dimen@3\ex@\advance\dimen@\dp\z@
 \hbox{\lower\dimen@\hbox{%
 \vbox{\hrule height.4\ex@
 \hbox{\vrule width.4\ex@\hskip3\ex@\vbox{\vskip3\ex@\boxz@\vskip3\ex@}%
 \hskip3\ex@\vrule width.4\ex@}\hrule height.4\ex@}%
 }}}
\message{commutative diagrams,}
\let\ampersand@\relax
\newdimen\minaw@
\minaw@11.11128\ex@
\newdimen\minCDaw@
\minCDaw@2.5pc
\def\minCDarrowwidth#1{\RIfMIfI@\onlydmatherr@\minCDarrowwidth
 \else\minCDaw@#1\relax\fi\else\onlydmatherr@\minCDarrowwidth\fi}
\newif\ifCD@
\def\CD{\bgroup\vspace@\relax\let\ampersand@&\iffalse}\fi
 \CD@true\vcenter\bgroup\Let@\tabskip\z@skip\baselineskip20\ex@
 \lineskip3\ex@\lineskiplimit3\ex@\halign\bgroup
 &\hfill$\m@th##$\hfill\crcr}
\def\endCD{\crcr\egroup\egroup\egroup}
\newdimen\bigaw@
\atdef@>#1>#2>{\ampersand@                                                  
 \setboxz@h{$\m@th\ssize\;{#1}\;\;$}
 \setbox@ne\hbox{$\m@th\ssize\;{#2}\;\;$}
 \setbox\tw@\hbox{$\m@th#2$}
 \ifCD@\global\bigaw@\minCDaw@\else\global\bigaw@\minaw@\fi                 
 \ifdim\wdz@>\bigaw@\global\bigaw@\wdz@\fi
 \ifdim\wd@ne>\bigaw@\global\bigaw@\wd@ne\fi                                
 \ifCD@\enskip\fi                                                           
 \ifdim\wd\tw@>\z@
  \mathrel{\mathop{\hbox to\bigaw@{\rightarrowfill@\displaystyle}}%
    \limits^{#1}_{#2}}
 \else\mathrel{\mathop{\hbox to\bigaw@{\rightarrowfill@\displaystyle}}%
    \limits^{#1}}\fi                                                        
 \ifCD@\enskip\fi                                                          
 \ampersand@}                                                              
\atdef@<#1<#2<{\ampersand@\setboxz@h{$\m@th\ssize\;\;{#1}\;$}%
 \setbox@ne\hbox{$\m@th\ssize\;\;{#2}\;$}\setbox\tw@\hbox{$\m@th#2$}%
 \ifCD@\global\bigaw@\minCDaw@\else\global\bigaw@\minaw@\fi
 \ifdim\wdz@>\bigaw@\global\bigaw@\wdz@\fi
 \ifdim\wd@ne>\bigaw@\global\bigaw@\wd@ne\fi
 \ifCD@\enskip\fi
 \ifdim\wd\tw@>\z@
  \mathrel{\mathop{\hbox to\bigaw@{\leftarrowfill@\displaystyle}}%
       \limits^{#1}_{#2}}\else
  \mathrel{\mathop{\hbox to\bigaw@{\leftarrowfill@\displaystyle}}%
       \limits^{#1}}\fi
 \ifCD@\enskip\fi\ampersand@}
\begingroup
 \catcode`\~=\active \lccode`\~=`\@
 \lowercase{%
  \global\atdef@)#1)#2){~>#1>#2>}
  \global\atdef@(#1(#2({~<#1<#2<}}
\endgroup
\atdef@ A#1A#2A{\llap{$\m@th\vcenter{\hbox
 {$\ssize#1$}}$}\Big\uparrow\rlap{$\m@th\vcenter{\hbox{$\ssize#2$}}$}&&}
\atdef@ V#1V#2V{\llap{$\m@th\vcenter{\hbox
 {$\ssize#1$}}$}\Big\downarrow\rlap{$\m@th\vcenter{\hbox{$\ssize#2$}}$}&&}
\atdef@={&\enskip\mathrel
 {\vbox{\hrule width\minCDaw@\vskip3\ex@\hrule width
 \minCDaw@}}\enskip&}
\atdef@|{\Big\Vert&&}
\atdef@\vert{\Big\Vert&&}
\def\pretend#1\haswidth#2{\setboxz@h{$\m@th\scriptstyle{#2}$}\hbox
 to\wdz@{\hfill$\m@th\scriptstyle{#1}$\hfill}}
\message{poor man's bold,}
\def\pmb{\RIfM@\expandafter\mathpalette\expandafter\pmb@\else
 \expandafter\pmb@@\fi}
\def\pmb@@#1{\leavevmode\setboxz@h{#1}%
   \dimen@-\wdz@
   \kern-.5\ex@\copy\z@
   \kern\dimen@\kern.25\ex@\raise.4\ex@\copy\z@
   \kern\dimen@\kern.25\ex@\box\z@
}
\def\binrel@@#1{\ifdim\wd2<\z@\mathbin{#1}\else\ifdim\wd\tw@>\z@
 \mathrel{#1}\else{#1}\fi\fi}
\newdimen\pmbraise@
\def\pmb@#1#2{\setbox\thr@@\hbox{$\m@th#1{#2}$}%
 \setbox4\hbox{$\m@th#1\mkern.5mu$}\pmbraise@\wd4\relax
 \binrel@{#2}%
 \dimen@-\wd\thr@@
   \binrel@@{%
   \mkern-.8mu\copy\thr@@
   \kern\dimen@\mkern.4mu\raise\pmbraise@\copy\thr@@
   \kern\dimen@\mkern.4mu\box\thr@@
}}
\def\documentstyle#1{\W@{}\input #1.sty\relax}
\message{syntax check,}
\font\dummyft@=dummy
\fontdimen1 \dummyft@=\z@
\fontdimen2 \dummyft@=\z@
\fontdimen3 \dummyft@=\z@
\fontdimen4 \dummyft@=\z@
\fontdimen5 \dummyft@=\z@
\fontdimen6 \dummyft@=\z@
\fontdimen7 \dummyft@=\z@
\fontdimen8 \dummyft@=\z@
\fontdimen9 \dummyft@=\z@
\fontdimen10 \dummyft@=\z@
\fontdimen11 \dummyft@=\z@
\fontdimen12 \dummyft@=\z@
\fontdimen13 \dummyft@=\z@
\fontdimen14 \dummyft@=\z@
\fontdimen15 \dummyft@=\z@
\fontdimen16 \dummyft@=\z@
\fontdimen17 \dummyft@=\z@
\fontdimen18 \dummyft@=\z@
\fontdimen19 \dummyft@=\z@
\fontdimen20 \dummyft@=\z@
\fontdimen21 \dummyft@=\z@
\fontdimen22 \dummyft@=\z@
\def\fontlist@{\\{\tenrm}\\{\sevenrm}\\{\fiverm}\\{\teni}\\{\seveni}%
 \\{\fivei}\\{\tensy}\\{\sevensy}\\{\fivesy}\\{\tenex}\\{\tenbf}\\{\sevenbf}%
 \\{\fivebf}\\{\tensl}\\{\tenit}}
\def\font@#1=#2 {\rightappend@#1\to\fontlist@\font#1=#2 }
\def\dodummy@{{\def\\##1{\global\let##1\dummyft@}\fontlist@}}
\def\nopages@{\output{\setbox\z@\box\@cclv \deadcycles\z@}%
 \alloc@5\toks\toksdef\@cclvi\output}
\let\galleys\nopages@
\newif\ifsyntax@
\newcount\countxviii@
\def\syntax{\syntax@true\dodummy@\countxviii@\count18
 \loop\ifnum\countxviii@>\m@ne\textfont\countxviii@=\dummyft@
 \scriptfont\countxviii@=\dummyft@\scriptscriptfont\countxviii@=\dummyft@
 \advance\countxviii@\m@ne\repeat                                           
 \dummyft@\tracinglostchars\z@\nopages@\frenchspacing\hbadness\@M}
\def\first@#1#2\end{#1}
\def\printoptions{\W@{Do you want S(yntax check),
  G(alleys) or P(ages)?}%
 \message{Type S, G or P, followed by <return>: }%
 \begingroup 
 \endlinechar\m@ne 
 \read\m@ne to\ans@
 \edef\ans@{\uppercase{\def\noexpand\ans@{%
   \expandafter\first@\ans@ P\end}}}%
 \expandafter\endgroup\ans@
 \if\ans@ P
 \else \if\ans@ S\syntax
 \else \if\ans@ G\galleys
 \else\message{? Unknown option: \ans@; using the `pages' option.}%
 \fi\fi\fi}
\def\alloc@#1#2#3#4#5{\global\advance\count1#1by\@ne
 \ch@ck#1#4#2\allocationnumber=\count1#1
 \global#3#5=\allocationnumber
 \ifalloc@\wlog{\string#5=\string#2\the\allocationnumber}\fi}
\def\document{\def\alloclist@{}\def\fontlist@{}}
\let\enddocument\bye

\let\proclaim\undefined
\let\footnote\undefined
\let\=\undefined
\let\>\undefined

\catcode`\@=\active
\message{... finished}

\expandafter\ifx\csname mathdefs.tex\endcsname\relax
  \expandafter\gdef\csname mathdefs.tex\endcsname{}
\else \message{Hey!  Apparently you were trying to
  \string\input{mathdefs.tex} twice.   This does not make sense.} 
\errmessage{Please edit your file (probably \jobname.tex) and remove
any duplicate ``\string\input'' lines}\endinput\fi




\catcode`\X=12\catcode`\@=11

\def\n@wcount{\alloc@0\count\countdef\insc@unt}
\def\n@wwrite{\alloc@7\write\chardef\sixt@@n}
\def\n@wread{\alloc@6\read\chardef\sixt@@n}
\def\r@s@t{\relax}\def\v@idline{\par}\def\@mputate#1/{#1}
\def\l@c@l#1X{\firstpart.#1}\def\gl@b@l#1X{#1}\def\t@d@l#1X{{}}

\def\crossrefs#1{\ifx\all#1\let\tr@ce=\all\else\def\tr@ce{#1,}\fi
   \n@wwrite\cit@tionsout\openout\cit@tionsout=\jobname.cit 
   \write\cit@tionsout{\tr@ce}\expandafter\setfl@gs\tr@ce,}
\def\setfl@gs#1,{\def\@{#1}\ifx\@\empty\let\next=\relax
   \else\let\next=\setfl@gs\expandafter\xdef
   \csname#1tr@cetrue\endcsname{}\fi\next}
\def\m@ketag#1#2{\expandafter\n@wcount\csname#2tagno\endcsname
     \csname#2tagno\endcsname=0\let\tail=\all\xdef\all{\tail#2,}
   \ifx#1\l@c@l\let\tail=\r@s@t\xdef\r@s@t{\csname#2tagno\endcsname=0\tail}\fi
   \expandafter\gdef\csname#2cite\endcsname##1{\expandafter
     \ifx\csname#2tag##1\endcsname\relax?\else\csname#2tag##1\endcsname\fi
     \expandafter\ifx\csname#2tr@cetrue\endcsname\relax\else
     \write\cit@tionsout{#2tag ##1 cited on page \folio.}\fi}
   \expandafter\gdef\csname#2page\endcsname##1{\expandafter
     \ifx\csname#2page##1\endcsname\relax?\else\csname#2page##1\endcsname\fi
     \expandafter\ifx\csname#2tr@cetrue\endcsname\relax\else
     \write\cit@tionsout{#2tag ##1 cited on page \folio.}\fi}
   \expandafter\gdef\csname#2tag\endcsname##1{\expandafter
      \ifx\csname#2check##1\endcsname\relax
      \expandafter\xdef\csname#2check##1\endcsname{}%
      \else\immediate\write16{Warning: #2tag ##1 used more than once.}\fi
      \multit@g{#1}{#2}##1/X%
      \write\t@gsout{#2tag ##1 assigned number \csname#2tag##1\endcsname\space
      on page \number\count0.}%
   \csname#2tag##1\endcsname}}

\def\multit@g#1#2#3/#4X{\def\t@mp{#4}\ifx\t@mp\empty%
      \global\advance\csname#2tagno\endcsname by 1 
      \expandafter\xdef\csname#2tag#3\endcsname
      {#1\number\csname#2tagno\endcsnameX}%
   \else\expandafter\ifx\csname#2last#3\endcsname\relax
      \expandafter\n@wcount\csname#2last#3\endcsname
      \global\advance\csname#2tagno\endcsname by 1 
      \expandafter\xdef\csname#2tag#3\endcsname
      {#1\number\csname#2tagno\endcsnameX}
      \write\t@gsout{#2tag #3 assigned number \csname#2tag#3\endcsname\space
      on page \number\count0.}\fi
   \global\advance\csname#2last#3\endcsname by 1
   \def\t@mp{\expandafter\xdef\csname#2tag#3/}%
   \expandafter\t@mp\@mputate#4\endcsname
   {\csname#2tag#3\endcsname\lastpart{\csname#2last#3\endcsname}}\fi}
\def\t@gs#1{\def\all{}\m@ketag#1e\m@ketag#1s\m@ketag\t@d@l p
\let\realscite\scite
\let\realstag\stag
   \m@ketag\gl@b@l r \n@wread\t@gsin
   \openin\t@gsin=\jobname.tgs \re@der \closein\t@gsin
   \n@wwrite\t@gsout\openout\t@gsout=\jobname.tgs }
\outer\def\localtags{\t@gs\l@c@l}
\outer\def\globaltags{\t@gs\gl@b@l}
\outer\def\newlocaltag#1{\m@ketag\l@c@l{#1}}
\outer\def\newglobaltag#1{\m@ketag\gl@b@l{#1}}

\newif\ifpr@ 
\def\m@kecs #1tag #2 assigned number #3 on page #4.%
   {\expandafter\gdef\csname#1tag#2\endcsname{#3}
   \expandafter\gdef\csname#1page#2\endcsname{#4}
   \ifpr@\expandafter\xdef\csname#1check#2\endcsname{}\fi}
\def\re@der{\ifeof\t@gsin\let\next=\relax\else
   \read\t@gsin to\t@gline\ifx\t@gline\v@idline\else
   \expandafter\m@kecs \t@gline\fi\let \next=\re@der\fi\next}
\def\pretags#1{\pr@true\pret@gs#1,,}
\def\pret@gs#1,{\def\@{#1}\ifx\@\empty\let\n@xtfile=\relax
   \else\let\n@xtfile=\pret@gs \openin\t@gsin=#1.tgs \message{#1} \re@der 
   \closein\t@gsin\fi \n@xtfile}

\newcount\sectno\sectno=0\newcount\subsectno\subsectno=0
\newif\ifultr@local \def\ultralocal{\ultr@localtrue}
\def\firstpart{\number\sectno}
\def\lastpart#1{\ifcase#1 \or a\or b\or c\or d\or e\or f\or g\or h\or 
   i\or k\or l\or m\or n\or o\or p\or q\or r\or s\or t\or u\or v\or w\or 
   x\or y\or z \fi}

\def\resetall{\global\advance\sectno by 1\subsectno=0
   \gdef\firstpart{\number\sectno}\r@s@t}
\def\resetsub{\global\advance\subsectno by 1
   \gdef\firstpart{\number\sectno.\number\subsectno}\r@s@t}
\def\newsection#1\par{\resetall\vskip0pt plus.3\vsize\penalty-250
   \vskip0pt plus-.3\vsize\bigskip\bigskip
   \message{#1}\leftline{\bf#1}\nobreak\bigskip}
\def\subsection#1\par{\ifultr@local\resetsub\fi
   \vskip0pt plus.2\vsize\penalty-250\vskip0pt plus-.2\vsize
   \bigskip\smallskip\message{#1}\leftline{\bf#1}\nobreak\medskip}


\newdimen\marginshift

\newdimen\margindelta
\newdimen\marginmax
\newdimen\marginmin

\def\margininit{       
\marginmax=3 true cm                  
				      
\margindelta=0.1 true cm              
\marginmin=0.1true cm                 
\marginshift=\marginmin
}    

\def\t@gsjj#1,{\def\@{#1}\ifx\@\empty\let\next=\relax\else\let\next=\t@gsjj
   \def\@@{p}\ifx\@\@@\else
   \expandafter\gdef\csname#1cite\endcsname##1{\citejj{##1}}
   \expandafter\gdef\csname#1page\endcsname##1{?}
   \expandafter\gdef\csname#1tag\endcsname##1{\tagjj{##1}}\fi\fi\next}
\newif\ifshowstuffinmargin
\showstuffinmarginfalse
\def\jjtags{\ifx\shlhetal\relax 
  \else
\ifx\shlhetal\undefinedcontrolseq
\else
\showstuffinmargintrue
\ifx\all\relax\else\expandafter\t@gsjj\all,\fi\fi \fi
}

\def\tagjj#1{\realstag{#1}\mginpar{\zeigen{#1}}}
\def\citejj#1{\rechnen{#1}\mginpar{\zeigen{#1}}}     

\def\rechnen#1{\expandafter\ifx\csname stag#1\endcsname\relax ??\else
                           \csname stag#1\endcsname\fi}

\newdimen\theight

\def\marginfont{\sevenrm}

\def\trymarginbox#1{\setbox0=\hbox{\marginfont\hskip\marginshift #1}%
		\global\marginshift\wd0 
		\global\advance\marginshift\margindelta}

\def \mginpar#1{%
\ifvmode\setbox0\hbox to \hsize{\hfill\rlap{\marginfont\quad#1}}%
\ht0 0cm
\dp0 0cm
\box0\vskip-\baselineskip
\else 
             \vadjust{\trymarginbox{#1}%
		\ifdim\marginshift>\marginmax \global\marginshift\marginmin
			\trymarginbox{#1}%
                \fi
             \theight=\ht0
             \advance\theight by \dp0    \advance\theight by \lineskip
             \kern -\theight \vbox to \theight{\rightline{\rlap{\box0}}%
\vss}}\fi}


\def\t@gsoff#1,{\def\@{#1}\ifx\@\empty\let\next=\relax\else\let\next=\t@gsoff
   \def\@@{p}\ifx\@\@@\else
   \expandafter\gdef\csname#1cite\endcsname##1{\zeigen{##1}}
   \expandafter\gdef\csname#1page\endcsname##1{?}
   \expandafter\gdef\csname#1tag\endcsname##1{\zeigen{##1}}\fi\fi\next}
\def\verbatimtags{\showstuffinmarginfalse
\ifx\all\relax\else\expandafter\t@gsoff\all,\fi}
\def\zeigen#1{\hbox{$\langle$}#1\hbox{$\rangle$}}
\def\margincite#1{\ifshowstuffinmargin\mginpar{\rechnen{#1}}\fi}

\def\(#1){\edef\dot@g{\ifmmode\ifinner(\hbox{\noexpand\etag{#1}})
   \else\noexpand\eqno(\hbox{\noexpand\etag{#1}})\fi
   \else(\noexpand\ecite{#1})\fi}\dot@g}

\newif\ifbr@ck
\def\eat#1{}
\def\[#1]{\br@cktrue[\br@cket#1'X]}
\def\br@cket#1'#2X{\def\temp{#2}\ifx\temp\empty\let\next\eat
   \else\let\next\br@cket\fi
   \ifbr@ck\br@ckfalse\br@ck@t#1,X\else\br@cktrue#1\fi\next#2X}
\def\br@ck@t#1,#2X{\def\temp{#2}\ifx\temp\empty\let\neext\eat
   \else\let\neext\br@ck@t\def\temp{,}\fi
   \def\teemp{#1}\ifx\teemp\empty\else\rcite{#1}\fi\temp\neext#2X}
\def\resetbr@cket{\gdef\[##1]{[\rtag{##1}]}}
\def\references{\resetbr@cket\newsection References\par}

\newtoks\symb@ls\newtoks\s@mb@ls\newtoks\p@gelist\n@wcount\ftn@mber
    \ftn@mber=1\newif\ifftn@mbers\ftn@mbersfalse\newif\ifbyp@ge\byp@gefalse
\def\defm@rk{\ifftn@mbers\n@mberm@rk\else\symb@lm@rk\fi}
\def\n@mberm@rk{\xdef\m@rk{{\the\ftn@mber}}%
    \global\advance\ftn@mber by 1 }
\def\rot@te#1{\let\temp=#1\global#1=\expandafter\r@t@te\the\temp,X}
\def\r@t@te#1,#2X{{#2#1}\xdef\m@rk{{#1}}}
\def\b@@st#1{{$^{#1}$}}\def\str@p#1{#1}
\def\symb@lm@rk{\ifbyp@ge\rot@te\p@gelist\ifnum\expandafter\str@p\m@rk=1 
    \s@mb@ls=\symb@ls\fi\write\f@nsout{\number\count0}\fi \rot@te\s@mb@ls}
\def\byp@ge{\byp@getrue\n@wwrite\f@nsin\openin\f@nsin=\jobname.fns 
    \n@wcount\currentp@ge\currentp@ge=0\p@gelist={0}
    \re@dfns\closein\f@nsin\rot@te\p@gelist
    \n@wread\f@nsout\openout\f@nsout=\jobname.fns }
\def\m@kelist#1X#2{{#1,#2}}
\def\re@dfns{\ifeof\f@nsin\let\next=\relax\else\read\f@nsin to \f@nline
    \ifx\f@nline\v@idline\else\let\t@mplist=\p@gelist
    \ifnum\currentp@ge=\f@nline
    \global\p@gelist=\expandafter\m@kelist\the\t@mplistX0
    \else\currentp@ge=\f@nline
    \global\p@gelist=\expandafter\m@kelist\the\t@mplistX1\fi\fi
    \let\next=\re@dfns\fi\next}
\def\symbols#1{\symb@ls={#1}\s@mb@ls=\symb@ls} 
\def\bigsymbol{\textstyle}
\symbols{\bigsymbol\ast,\dagger,\ddagger,\sharp,\flat,\natural,\star}
\def\ftnumbers{\ftn@mberstrue} \def\ftsymbols{\ftn@mbersfalse}
\def\paginal{\byp@ge} \def\resetftnumbers{\ftn@mber=1}
\def\ftnote#1{\defm@rk\expandafter\expandafter\expandafter\footnote
    \expandafter\b@@st\m@rk{#1}}

\long\def\jump#1\endjump{}
\def\ssum{\mathop{\lower .1em\hbox{$\textstyle\Sigma$}}\nolimits}

\def\qed{\nobreak\kern 1em \vrule height .5em width .5em depth 0em}
\def\newneq{\hbox{\rlap{\hbox to 1\wd9{\hss$=$\hss}}\raise .1em 
   \hbox to 1\wd9{\hss$\scriptscriptstyle/$\hss}}}
\def\subsetne{\setbox9 = \hbox{$\subset$}\mathrel{\hbox{\rlap
   {\lower .4em \newneq}\raise .13em \hbox{$\subset$}}}}
\def\supsetne{\setbox9 = \hbox{$\subset$}\mathrel{\hbox{\rlap
   {\lower .4em \newneq}\raise .13em \hbox{$\supset$}}}}

\def\vbar{\mathchoice{\vrule height6.3ptdepth-.5ptwidth.8pt\kern-.8pt}
   {\vrule height6.3ptdepth-.5ptwidth.8pt\kern-.8pt}
   {\vrule height4.1ptdepth-.35ptwidth.6pt\kern-.6pt}
   {\vrule height3.1ptdepth-.25ptwidth.5pt\kern-.5pt}}
\def\f@dge{\mathchoice{}{}{\mkern.5mu}{\mkern.8mu}}
\def\b@c#1#2{{\rm \mkern#2mu\vbar\mkern-#2mu#1}}
\def\b@b#1{{\rm I\mkern-3.5mu #1}}
\def\b@a#1#2{{\rm #1\mkern-#2mu\f@dge #1}}
\def\bb#1{{\count4=`#1 \advance\count4by-64 \ifcase\count4\or\b@a A{11.5}\or
   \b@b B\or\b@c C{5}\or\b@b D\or\b@b E\or\b@b F \or\b@c G{5}\or\b@b H\or
   \b@b I\or\b@c J{3}\or\b@b K\or\b@b L \or\b@b M\or\b@b N\or\b@c O{5} \or
   \b@b P\or\b@c Q{5}\or\b@b R\or\b@a S{8}\or\b@a T{10.5}\or\b@c U{5}\or
   \b@a V{12}\or\b@a W{16.5}\or\b@a X{11}\or\b@a Y{11.7}\or\b@a Z{7.5}\fi}}

\catcode`\X=11 \catcode`\@=12


\expandafter\ifx\csname citeadd.tex\endcsname\relax
\expandafter\gdef\csname citeadd.tex\endcsname{}
\else \message{Hey!  Apparently you were trying to
\string\input{citeadd.tex} twice.   This does not make sense.} 
\errmessage{Please edit your file (probably \jobname.tex) and remove
any duplicate ``\string\input'' lines}\endinput\fi

\sectno=-1   
\localtags
\jjtags
\NoBlackBoxes
\define\mr{\medskip\roster}
\define\sn{\smallskip\noindent}
\define\mn{\medskip\noindent}
\define\bn{\bigskip\noindent}
\define\ub{\underbar}
\define\wilog{\text{without loss of generality}}
\define\ermn{\endroster\medskip\noindent}

\define\dbcu{\dsize\bigcup}
\define \nl{\newline}
\magnification=\magstep 1
\documentstyle{amsppt}

{    
\catcode`@11

\ifx\alicetwothousandloaded@\relax
  \endinput\else\global\let\alicetwothousandloaded@\relax\fi

\gdef\subjclass{\let\savedef@\subjclass
 \def\subjclass##1\endsubjclass{\let\subjclass\savedef@
   \toks@{\def\usualspace{{\rm\enspace}}\eightpoint}%
   \toks@@{##1\unskip.}%
   \edef\thesubjclass@{\the\toks@
     \frills@{{\noexpand\rm2000 {\noexpand\it Mathematics Subject
       Classification}.\noexpand\enspace}}%
     \the\toks@@}}%
  \nofrillscheck\subjclass}
} 


\expandafter\ifx\csname alice2jlem.tex\endcsname\relax
  \expandafter\xdef\csname alice2jlem.tex\endcsname{\the\catcode`@}
\else \message{Hey!  Apparently you were trying to
\string\input{alice2jlem.tex}  twice.   This does not make sense.}
\errmessage{Please edit your file (probably \jobname.tex) and remove
any duplicate ``\string\input'' lines}\endinput\fi

\expandafter\ifx\csname bib4plain.tex\endcsname\relax
  \expandafter\gdef\csname bib4plain.tex\endcsname{}
\else \message{Hey!  Apparently you were trying to \string\input
  bib4plain.tex twice.   This does not make sense.}
\errmessage{Please edit your file (probably \jobname.tex) and remove
any duplicate ``\string\input'' lines}\endinput\fi

\def\renewcommand{\newcommand}	       
\edef\cite{\the\catcode`@}%
\catcode`@ = 11
\let\@oldatcatcode = \cite
\chardef\@letter = 11
\chardef\@other = 12
%
%
%
%
\def\@innerdef#1#2{\edef#1{\expandafter\noexpand\csname #2\endcsname}}%
%
%
\@innerdef\@innernewcount{newcount}%
\@innerdef\@innernewdimen{newdimen}%
\@innerdef\@innernewif{newif}%
\@innerdef\@innernewwrite{newwrite}%
%
%
%
\def\@gobble#1{}%
%
%
%
\ifx\inputlineno\@undefined
   \let\@linenumber = \empty 
\else
   \def\@linenumber{\the\inputlineno:\space}%
\fi
%
%
%
\def\@futurenonspacelet#1{\def\cs{#1}%
   \afterassignment\@stepone\let\@nexttoken=
}%
\begingroup 
\def\\{\global\let\@stoken= }%
\\ 
\endgroup
\def\@stepone{\expandafter\futurelet\cs\@steptwo}%
\def\@steptwo{\expandafter\ifx\cs\@stoken\let\@@next=\@stepthree
   \else\let\@@next=\@nexttoken\fi \@@next}%
\def\@stepthree{\afterassignment\@stepone\let\@@next= }%
%
%
%
\def\@getoptionalarg#1{%
   \let\@optionaltemp = #1%
   \let\@optionalnext = \relax
   \@futurenonspacelet\@optionalnext\@bracketcheck
}%
%
%
\def\@bracketcheck{%
   \ifx [\@optionalnext
      \expandafter\@@getoptionalarg
   \else
      \let\@optionalarg = \empty
      \expandafter\@optionaltemp
   \fi
}%
\def\@@getoptionalarg[#1]{%
   \def\@optionalarg{#1}%
   \@optionaltemp
}%
%
%
%
\def\@nnil{\@nil}%
\def\@fornoop#1\@@#2#3{}%
\def\@for#1:=#2\do#3{%
   \edef\@fortmp{#2}%
   \ifx\@fortmp\empty \else
      \expandafter\@forloop#2,\@nil,\@nil\@@#1{#3}%
   \fi
}%
\def\@forloop#1,#2,#3\@@#4#5{\def#4{#1}\ifx #4\@nnil \else
       #5\def#4{#2}\ifx #4\@nnil \else#5\@iforloop #3\@@#4{#5}\fi\fi
}%
\def\@iforloop#1,#2\@@#3#4{\def#3{#1}\ifx #3\@nnil
       \let\@nextwhile=\@fornoop \else
      #4\relax\let\@nextwhile=\@iforloop\fi\@nextwhile#2\@@#3{#4}%
}%
%
%
%
\@innernewif\if@fileexists
\def\@testfileexistence{\@getoptionalarg\@finishtestfileexistence}%
\def\@finishtestfileexistence#1{%
   \begingroup
      \def\extension{#1}%
      \immediate\openin0 =
         \ifx\@optionalarg\empty\jobname\else\@optionalarg\fi
         \ifx\extension\empty \else .#1\fi
         \space
      \ifeof 0
         \global\@fileexistsfalse
      \else
         \global\@fileexiststrue
      \fi
      \immediate\closein0
   \endgroup
}%
%
%
%
%
\def\bibliographystyle#1{%
   \@readauxfile
   \@writeaux{\string\bibstyle{#1}}%
}%
\let\bibstyle = \@gobble
%
%
\let\bblfilebasename = \jobname
\def\bibliography#1{%
   \@readauxfile
   \@writeaux{\string\bibdata{#1}}%
   \@testfileexistence[\bblfilebasename]{bbl}%
   \if@fileexists
      \nobreak
      \@readbblfile
   \fi
}%
\let\bibdata = \@gobble
%
%
\def\nocite#1{%
   \@readauxfile
   \@writeaux{\string\citation{#1}}%
}%
\@innernewif\if@notfirstcitation
%
%
\def\cite{\@getoptionalarg\@cite}%
%
%
\def\@cite#1{%
   \let\@citenotetext = \@optionalarg
   \printcitestart
   \nocite{#1}%
   \@notfirstcitationfalse
   \@for \@citation :=#1\do
   {%
      \expandafter\@onecitation\@citation\@@
   }%
   \ifx\empty\@citenotetext\else
      \printcitenote{\@citenotetext}%
   \fi
   \printcitefinish
}%
\newif\ifweareinprivate
\weareinprivatetrue
\ifx\shlhetal\undefinedcontrolseq\weareinprivatefalse\fi
\ifx\shlhetal\relax\weareinprivatefalse\fi
\def\@onecitation#1\@@{%
   \if@notfirstcitation
      \printbetweencitations
   \fi
   \expandafter \ifx \csname\@citelabel{#1}\endcsname \relax
      \if@citewarning
         \message{\@linenumber Undefined citation `#1'.}%
      \fi
     \ifweareinprivate
      \expandafter\gdef\csname\@citelabel{#1}\endcsname{%
\strut 
\vadjust{\vskip-\dp\strutbox
\vbox to 0pt{\vss\parindent0cm \leftskip=\hsize 
\advance\leftskip3mm
\advance\hsize 4cm\strut\openup-4pt 
\rightskip 0cm plus 1cm minus 0.5cm ?  #1 ?\strut}}
         {\tt
            \escapechar = -1
            \nobreak\hskip0pt\pfeilsw
            \expandafter\string\csname#1\endcsname
\pfeilso
            \nobreak\hskip0pt
         }%
      }%
     \else  
      \expandafter\gdef\csname\@citelabel{#1}\endcsname{%
            {\tt\expandafter\string\csname#1\endcsname}
      }%
     \fi  
   \fi
   \csname\@citelabel{#1}\endcsname
   \@notfirstcitationtrue
}%
%
%
\def\@citelabel#1{b@#1}%
%
%
\def\@citedef#1#2{\expandafter\gdef\csname\@citelabel{#1}\endcsname{#2}}%
%
%
%
\def\@readbblfile{%
   \ifx\@itemnum\@undefined
      \@innernewcount\@itemnum
   \fi
   \begingroup
      \def\begin##1##2{%
         \setbox0 = \hbox{\biblabelcontents{##2}}%
         \biblabelwidth = \wd0
      }%
      \def\end##1{}
      %
      %
      \@itemnum = 0
      \def\bibitem{\@getoptionalarg\@bibitem}%
      \def\@bibitem{%
         \ifx\@optionalarg\empty
            \expandafter\@numberedbibitem
         \else
            \expandafter\@alphabibitem
         \fi
      }%
      \def\@alphabibitem##1{%
         \expandafter \xdef\csname\@citelabel{##1}\endcsname {\@optionalarg}%
         \ifx\biblabelprecontents\@undefined
            \let\biblabelprecontents = \relax
         \fi
         \ifx\biblabelpostcontents\@undefined
            \let\biblabelpostcontents = \hss
         \fi
         \@finishbibitem{##1}%
      }%
      \def\@numberedbibitem##1{%
         \advance\@itemnum by 1
         \expandafter \xdef\csname\@citelabel{##1}\endcsname{\number\@itemnum}%
         \ifx\biblabelprecontents\@undefined
            \let\biblabelprecontents = \hss
         \fi
         \ifx\biblabelpostcontents\@undefined
            \let\biblabelpostcontents = \relax
         \fi
         \@finishbibitem{##1}%
      }%
      \def\@finishbibitem##1{%
         \biblabelprint{\csname\@citelabel{##1}\endcsname}%
         \@writeaux{\string\@citedef{##1}{\csname\@citelabel{##1}\endcsname}}%
         \ignorespaces
      }%
      %
      %
      \let\em = \bblem
      \let\newblock = \bblnewblock
      \let\sc = \bblsc
      \frenchspacing
      \clubpenalty = 4000 \widowpenalty = 4000
      \tolerance = 10000 \hfuzz = .5pt
      \everypar = {\hangindent = \biblabelwidth
                      \advance\hangindent by \biblabelextraspace}%
      \bblrm
      \parskip = 1.5ex plus .5ex minus .5ex
      \biblabelextraspace = .5em
      \bblhook
      \input \bblfilebasename.bbl
   \endgroup
}%
%
%
\@innernewdimen\biblabelwidth
\@innernewdimen\biblabelextraspace
%
%
%
\def\biblabelprint#1{%
   \noindent
   \hbox to \biblabelwidth{%
      \biblabelprecontents
      \biblabelcontents{#1}%
      \biblabelpostcontents
   }%
   \kern\biblabelextraspace
}%
%
%
%
\def\biblabelcontents#1{{\bblrm [#1]}}%
%
%
\def\bblrm{\rm}%
%
%
\def\bblem{\it}%
%
%
\def\bblsc{\ifx\@scfont\@undefined
              \font\@scfont = cmcsc10
           \fi
           \@scfont
}%
%
%
\def\bblnewblock{\hskip .11em plus .33em minus .07em }%
%
%
\let\bblhook = \empty
%
%
%
\def\printcitestart{[}
\def\printcitefinish{]}
\def\printbetweencitations{, }
\def\printcitenote#1{, #1}
%
%
%
\let\citation = \@gobble
%
%
%
\@innernewcount\@numparams
%
%
\def\newcommand#1{%
   \def\@commandname{#1}%
   \@getoptionalarg\@continuenewcommand
}%
%
%
\def\@continuenewcommand{%
   \@numparams = \ifx\@optionalarg\empty 0\else\@optionalarg \fi \relax
   \@newcommand
}%
%
%
\def\@newcommand#1{%
   \def\@startdef{\expandafter\edef\@commandname}%
   \ifnum\@numparams=0
      \let\@paramdef = \empty
   \else
      \ifnum\@numparams>9
         \errmessage{\the\@numparams\space is too many parameters}%
      \else
         \ifnum\@numparams<0
            \errmessage{\the\@numparams\space is too few parameters}%
         \else
            \edef\@paramdef{%
               \ifcase\@numparams
                  \empty  No arguments.
               \or ####1%
               \or ####1####2%
               \or ####1####2####3%
               \or ####1####2####3####4%
               \or ####1####2####3####4####5%
               \or ####1####2####3####4####5####6%
               \or ####1####2####3####4####5####6####7%
               \or ####1####2####3####4####5####6####7####8%
               \or ####1####2####3####4####5####6####7####8####9%
               \fi
            }%
         \fi
      \fi
   \fi
   \expandafter\@startdef\@paramdef{#1}%
}%
%
%
%
%
\def\@readauxfile{%
   \if@auxfiledone \else 
      \global\@auxfiledonetrue
      \@testfileexistence{aux}%
      \if@fileexists
         \begingroup
            \endlinechar = -1
            \catcode`@ = 11
            \input \jobname.aux
         \endgroup
      \else
         \message{\@undefinedmessage}%
         \global\@citewarningfalse
      \fi
      \immediate\openout\@auxfile = \jobname.aux
   \fi
}%
%
%
\newif\if@auxfiledone
\ifx\noauxfile\@undefined \else \@auxfiledonetrue\fi
%
%
%
%
\@innernewwrite\@auxfile
\def\@writeaux#1{\ifx\noauxfile\@undefined \write\@auxfile{#1}\fi}%
%
%
%
\ifx\@undefinedmessage\@undefined
   \def\@undefinedmessage{No .aux file; I won't give you warnings about
                          undefined citations.}%
\fi
%
%
\@innernewif\if@citewarning
\ifx\noauxfile\@undefined \@citewarningtrue\fi
%
%
%
\catcode`@ = \@oldatcatcode

\def\pfeilso{\leavevmode
            \vrule width 1pt height9pt depth 0pt\relax
           \vrule width 1pt height8.7pt depth 0pt\relax
           \vrule width 1pt height8.3pt depth 0pt\relax
           \vrule width 1pt height8.0pt depth 0pt\relax
           \vrule width 1pt height7.7pt depth 0pt\relax
            \vrule width 1pt height7.3pt depth 0pt\relax
            \vrule width 1pt height7.0pt depth 0pt\relax
            \vrule width 1pt height6.7pt depth 0pt\relax
            \vrule width 1pt height6.3pt depth 0pt\relax
            \vrule width 1pt height6.0pt depth 0pt\relax
            \vrule width 1pt height5.7pt depth 0pt\relax
            \vrule width 1pt height5.3pt depth 0pt\relax
            \vrule width 1pt height5.0pt depth 0pt\relax
            \vrule width 1pt height4.7pt depth 0pt\relax
            \vrule width 1pt height4.3pt depth 0pt\relax
            \vrule width 1pt height4.0pt depth 0pt\relax
            \vrule width 1pt height3.7pt depth 0pt\relax
            \vrule width 1pt height3.3pt depth 0pt\relax
            \vrule width 1pt height3.0pt depth 0pt\relax
            \vrule width 1pt height2.7pt depth 0pt\relax
            \vrule width 1pt height2.3pt depth 0pt\relax
            \vrule width 1pt height2.0pt depth 0pt\relax
            \vrule width 1pt height1.7pt depth 0pt\relax
            \vrule width 1pt height1.3pt depth 0pt\relax
            \vrule width 1pt height1.0pt depth 0pt\relax
            \vrule width 1pt height0.7pt depth 0pt\relax
            \vrule width 1pt height0.3pt depth 0pt\relax}

\def\pfeilsw{ \leavevmode 
            \vrule width 1pt height0.3pt depth 0pt\relax
            \vrule width 1pt height0.7pt depth 0pt\relax
            \vrule width 1pt height1.0pt depth 0pt\relax
            \vrule width 1pt height1.3pt depth 0pt\relax
            \vrule width 1pt height1.7pt depth 0pt\relax
            \vrule width 1pt height2.0pt depth 0pt\relax
            \vrule width 1pt height2.3pt depth 0pt\relax
            \vrule width 1pt height2.7pt depth 0pt\relax
            \vrule width 1pt height3.0pt depth 0pt\relax
            \vrule width 1pt height3.3pt depth 0pt\relax
            \vrule width 1pt height3.7pt depth 0pt\relax
            \vrule width 1pt height4.0pt depth 0pt\relax
            \vrule width 1pt height4.3pt depth 0pt\relax
            \vrule width 1pt height4.7pt depth 0pt\relax
            \vrule width 1pt height5.0pt depth 0pt\relax
            \vrule width 1pt height5.3pt depth 0pt\relax
            \vrule width 1pt height5.7pt depth 0pt\relax
            \vrule width 1pt height6.0pt depth 0pt\relax
            \vrule width 1pt height6.3pt depth 0pt\relax
            \vrule width 1pt height6.7pt depth 0pt\relax
            \vrule width 1pt height7.0pt depth 0pt\relax
            \vrule width 1pt height7.3pt depth 0pt\relax
            \vrule width 1pt height7.7pt depth 0pt\relax
            \vrule width 1pt height8.0pt depth 0pt\relax
            \vrule width 1pt height8.3pt depth 0pt\relax
            \vrule width 1pt height8.7pt depth 0pt\relax
            \vrule width 1pt height9pt depth 0pt\relax
      }


\def\widestnumber#1#2{}

\def\citewarning#1{\ifx\shlhetal\relax 
    \else
    \par{#1}\par
    \fi
}

\def\rm{\fam0 \tenrm}

\def\fakesubhead#1\endsubhead{\bigskip\noindent{\bf#1}\par}



%
%
%

%

\font\textrsfs=rsfs10
\font\scriptrsfs=rsfs7
\font\scriptscriptrsfs=rsfs5

\newfam\rsfsfam
\textfont\rsfsfam=\textrsfs
\scriptfont\rsfsfam=\scriptrsfs
\scriptscriptfont\rsfsfam=\scriptscriptrsfs

\edef\oldcatcodeofat{\the\catcode`\@}
\catcode`\@11

\def\Cal@@#1{\noaccents@ \fam \rsfsfam #1}

\catcode`\@\oldcatcodeofat


\expandafter\ifx \csname margininit\endcsname \relax\else\margininit\fi

\pageheight{8.5truein}
\topmatter
\title{Note on $\omega$--nw--nep forcing notions\\
 Sh711} \endtitle
\author {Saharon Shelah \thanks {\null\newline I would like to thank 
Alice Leonhardt for the beautiful typing. \null\newline
This research was partially supported by the Israel Science Foundation
founded by the Israel Academy of Sciences and Humanities. 
\null\newline
Publication 711 
\null\newline
First Typed at RU - 01/10/24 \null\newline 
Latest Revision - 03/Feb/6} \endthanks} \endauthor  
\affil{Institute of Mathematics\\
 The Hebrew University\\
 Jerusalem, Israel
 \medskip
 Rutgers University\\
 Mathematics Department\\
 New Brunswick, NJ  USA} \endaffil

\abstract  We prove that if $\Bbb Q$ is a nw-nep forcing then it
cannot add a dominating real.  We also prove that $A$-moeba forcing
cannot be ${\Cal P}(X)/I$ if $I$ is an $\aleph_1$-complete ideal.
 \endabstract
\endtopmatter
\document  
 
\newpage

\head {\S1 Nwnep forcing notions} \endhead  \resetall \sectno=1
\bigskip

The results were announced in \cite{Sh:666}; the original motivation
was a question of Kambureti's solved in '97.  The definitions below
are ad-hoc, to simplify the presentation for the claims here.  See
more in \cite{Sh:666}.  An older connected work is \cite{Bn95}. We
thank Jakob Kellner for reading the paper extremely carefully
resulting in considerable improvement of the presentation and readability.
\bigskip

\definition{\stag{xa.1} Definition}  1) We say $N$ is 
$x-1$-nw-candidate \ub{if}, fixing some strong limit $\chi$, (a) or (b):
\mr
\item "{$(a)$}"  $N \prec ({\Cal H}(\chi),\in)$ is countable, $x \in
N$
\sn
\item "{$(b)$}"  for some $N_1$ as in (a), $N$ is an elementary
extension of $N_1$ not increasing $\omega^N$; i.e., if $N_1 \models g <
\omega$ then $g \in N_1$ and for notion $N \models$ ``$y$ is a real,
i.e., a set of integers" $\Rightarrow y = \{n:N \models n \in y\}$ 
(so really it should have one two-place relation
$R,R^N$ is the membership relation in $N$; but we shall write $E^N$ or
$N \models x \in y$). 
\ermn
2) We say $N$ is a standard $x-2$-nw-candidate \ub{if} (for $\chi$ as
above) (a) holds or
\mr
\item "{$(b)'$}"  for some $N_1$ as in (a), $N$ is a forcing
extension of $N$.
\ermn
3) We say that $\Bbb Q$ is 1-nw-nep forcing \ub{if} $\Bbb Q$ is a pair of
formulas $\bar \varphi = (\varphi_0(x),\varphi_1(x,y))$, in the
language of set theory such that we write $\Bbb Q$-candidate instead
of $\bar \varphi-1$-nw-candidate
\mr
\item "{$(a)$}"  $\varphi_0(x)$ defines a set of reals (= set of
members of $\Bbb Q$)
\sn
\item "{$(b)$}"  $\varphi_1(x,y)$ defines a set of pairs of reals, a
quasi order on $\{x:\varphi_0(x)\}$, this is $\le_{\Bbb Q}$
\sn
\item "{$(c)$}"  $\varphi_0,\varphi_1$ are $\Sigma^1_1$-formulas;
equivalently, are upward absolute from $\Bbb Q$-candidates, i.e. for a
$\Bbb Q$-candidate $N_1 \subseteq N_2,x \in \Bbb Q^{N_1} \Rightarrow x
\in \Bbb Q^{N_2} \Rightarrow x \in \Bbb Q$ and $x \le^{N_1}_{\Bbb Q} y
\Rightarrow x \le^{N_2}_{\Bbb Q} y \Rightarrow x \le_{\Bbb Q} y$
\sn
\item "{$(d)$}"  if $N$ is a $\Bbb Q$-candidate and $p \in \Bbb Q^N$
(i.e. $N \models p \in \Bbb Q$) then there is $q$ such that: $q \in
\Bbb Q,p \le_{\Bbb Q} q$ and $q$ is $\langle N,\Bbb Q \rangle$-generic, i.e. 
$$
\align
q \Vdash ``&\Bbb Q^N \cap G \text{ is a subset of } \Bbb Q^N, \text{
directed by} \\
  &\le^N_{\Bbb Q} \text{ and } G \cap {\Cal I}^N \ne 
\emptyset" \text{ whenever} \\
  &N \models ``{\Cal I} \subseteq \Bbb Q \text{ is predense}".
\endalign
$$
\ermn
4) We say $\Bbb Q$ is 2-nw-nep-forcing if above we replace $\bar
\varphi-1$-nw-candidate by $\bar \varphi-2$-candidates. \nl
5) Let nw-nep mean 1-nw-nep.
\enddefinition
\bigskip

\proclaim{\stag{x.3} Claim}  Assume $\Bbb Q$ is nw-nep. \ub{Then}
forcing with $\Bbb Q$ does not add a dominating real.
\endproclaim
\bigskip

\demo{Proof}  Toward contradiction assume $p^* \Vdash_{\Bbb Q}
``{\underset\tilde {}\to \eta^*}$ is a dominating real". \nl
Without loss of generality, $p^* \Vdash ``{\underset\tilde {}\to
\eta^*}$  is strictly increasing, 
${\underset\tilde {}\to \eta^*}(n)>n$". Let $\Gamma_0 = \{\eta \in
{}^{\omega >} \omega:\eta$ strictly increasing and $\eta(\ell) > \ell$
for $\ell < \ell g(\eta)\}$
so $p^* \Vdash_{\Bbb Q} ``{\underset\tilde {}\to \eta^*} \in \lim(\Gamma_0)"$.
As $\Bbb Q$ is nw-nep there is
$p^{**}$ such that
\mr
\item "{$\otimes_1$}"  $p^* \le_{\Bbb Q} p^{**}$ and for each $n$ there
is a countable ${\Cal J}^*_n \subseteq \Bbb Q$ which is an antichain
predense above $p^{**}$, such that each $p \in {\Cal J}^*_n$ forces a value to
$\underset\tilde {}\to \eta(n)$ and is above $p^*$ and above
some $p' \in {\Cal J}^*_m$ for each
$m < n$. \nl
[Why?  Let $N$ be a $\Bbb Q$-candidate to which $p^*$ belongs.  Inside
$N$ by induction on $n < \omega$ we choose ${\Cal J}_n \in N$ as above
except countability.  Let $p^{**}$ be above $p^*$ and be
$\langle N,G \rangle$-generic.]
\ermn
Clearly above any $p \ge p^{**}$ there are two incompatible elements
of $\Bbb Q$ so without loss of generality
\mr
\item "{$\otimes_2$}"    if $m<n$ and $p \in {\Cal J}_m$ then there
are infinitely many members of ${\Cal J}_n$ which are above $p$.
\ermn
Let $\Gamma$ denote a subset of $\Gamma_0$ closed under initial segments
such that $\langle\rangle\in\Gamma$ and $\eta\in\Gamma\ \Rightarrow\
(\exists^\infty n)(\eta \frown \langle n \rangle \in \Gamma)$. 
We can find $\Gamma,\bar k$ and choose $\bar p^*$ such that: 
\mr
\item "{$\otimes_3(\alpha)$}"  $\bar p^* = \langle p^*_\eta:\eta \in \Gamma
\rangle$ and $p^*_\eta \in \Bbb Q$, (in fact $p^*_\eta \in \{p^*\}
\cup \dbcu_n {\Cal J}^*_n$)
\sn
\item "{$(\beta)$}"   $\nu \vartriangleleft \eta \Rightarrow \Bbb Q
\models p^*_\nu \le p^*_\eta$, 
\sn
\item "{$(\gamma)$}"   for $n \in [1,\omega)$ we have 
$\langle p^*_\eta:\eta \in \Gamma \cap {}^n \omega \rangle$ 
is an antichain of $\Bbb Q$ predense above $p^{**}$, 
\sn
\item "{$(\delta)$}"   $\bar k = \langle k_\eta:\eta \in \Gamma
\rangle$ where $k_\eta < \omega$
\sn
\item "{$(\varepsilon)$}"  if 
$\eta \in {}^{\omega >} \omega$, $\eta \frown \langle m,n \rangle\in
\Gamma$, \ub{then} $p_{\eta \frown \langle m,n\rangle}$ forces a value to
${\underset\tilde {}\to \eta^*}(m)$, which we call $k_{\eta \char
94 <n,m>}$ and $m > k_{\eta \char 94 \langle n,m\rangle} > n$, 
\sn
\item "{$(\zeta)$}"   $p^*_{\langle\rangle} = p^*$ ($=$ no
information) and $\bar k_{\langle\rangle} = 1$.
\ermn
So we choose $\Gamma \cap {}^n\omega$ and $k_\eta,p_\eta$ 
(for $\eta\in\Gamma \cap {}^n \omega$) by induction on $n$ with $\eta
\ne <> \Rightarrow p^*_\eta \in \dbcu_n {\Cal J}^*_n$. 

For $n=0$ let $p_{<>} = p^*,k_{<>}=1$.

For $n+1$, for each $\eta \in {}^n \omega$ let $m = \sup \text{
Rang}(\eta)$ and let $\langle p_{\eta,j}:j < \omega \rangle$ list the
member of ${\Cal J}_m$ which are $p_\eta$, so for some $k_{\eta,j}$ we
have $p_{\eta,j} \Vdash_{\Bbb Q} {\underset\tilde {}\to \eta^*}(m) =
k_{\eta,j}$.  Let $f_\eta:\omega \rightarrow \omega$ be strictly
increasing such that $k_{\eta,j} < f_\eta(j)$ and sup Rang$(\eta) <
f_\eta(j)$ for $j < \omega$, and lastly, let $\Gamma \cap {}^{n+1}
\omega = \{\eta \char 94 \langle f(j) \rangle:\eta \in \Gamma \cap
{}^n \omega$ and $j <k\}$ and $k_{\eta \char 94 <f(j)>} = k_{\eta,j}$.
\nl
Let ${\underset\tilde {}\to \eta'}$ be the $\Bbb Q$-name of the
$\omega$-branch of $\Gamma$ such that $p^* \Vdash_{\Bbb Q}
``p_{{\underset\tilde {}\to \eta'} \restriction n} \in
{\underset\tilde {}\to G_{\Bbb Q}}"$ for each $n<\omega$.

We claim
\mr
\item "{$\boxtimes$}"   if $h:\Gamma\longrightarrow\omega$, then
$$
p^* \Vdash_{\Bbb Q} ``\text{for every large enough } n<\omega \text{ we
have } {\underset\tilde {}\to \eta'}(n) > h
({\underset\tilde {}\to \eta'} \restriction n)".
$$
\ermn
[Why?  Let $f_h:\omega\longrightarrow\omega$ be 

$$
f_h(n)=\sup\{h(\eta):\eta\in\Gamma \text{ and sup Range}(\eta) \le n\}+1,
$$ 
\mn
Note that the supremum is over a finite set as every $\eta \in \Gamma$
is strictly increasing.  So 
assume $p^{**} \in G \subseteq \Bbb Q$, $G$ is generic over $\bold V$,
$\eta' = {\underset\tilde {}\to \eta'}[G]$, $\eta^* =
{\underset\tilde {}\to \eta}[G]$, and we shall find $n$ as required.
Clearly for some $n^* > 2$ we
have $m\in [n^*,\omega) \Rightarrow f_h(m) < {\underset\tilde {}\to
\eta^*}(m)$.  So assume $m \in [n^*,\omega)$ and we shall prove that
$h(\eta \restriction m) < \eta'(m)$, this suffices.  

So assume $m \in [n^*,\omega)$ then \nl
$h(\eta' \restriction (m+1)) \le
f_h(\eta'(m))$ by the definition of $f_h,\eta'$ being increasing,
\nl
$f_h(\eta'(m)) < \eta^*(\eta'(n))$ as $\eta'(m) \ge m \ge n^*$ and the
choice of $n^*$, \nl
$\eta^*(\eta'(m)) = k_{\eta' \restriction (m+2)}$ by clause
$(\varepsilon)$ and \nl
$k_{\eta' \restriction (m+2)} < \eta'(m+1)$ by $(\varepsilon)$.

For limit $\alpha \le \omega_1$, let $\Xi_\alpha =
\{\bar \rho:\bar \rho = \langle\rho_\delta:
\delta \le \alpha$ where $\delta$ is limit$\rangle$ and each
$\rho_\delta$ is a (strictly) increasing $\omega$-sequence 
converging to $\delta\}$.  For $\bar \rho \in \Xi_\alpha$, we 
define a function $g_{\bar \rho}$ from $\Gamma$ to $\alpha+1$, defining
$g_{\bar \rho}(\eta)$ by induction on ${\text{\rm lg\/}}(\eta)$
as follows:
\mr
\item "{$(B)(a)$}"   $g_{\bar{\rho}}(\langle\rangle)=\alpha$, 
\sn
\item "{${{}}(b)$}"   if $g_{\bar \rho}(\eta) = \beta+1$ 
and $\eta \frown \langle\ell \rangle \in \Gamma$, \ub{then} 
$g_{\bar \rho}(\eta \frown \langle\ell\rangle) = \beta$,
\sn
\item "{${{}}(c)$}"   if 
$g_{\bar \rho}(\eta)=\delta$, $\delta$ a limit ordinal and 
$\eta \frown \langle\ell\rangle\in\Gamma$, \ub{then} $g_{\bar \rho}(\eta
\frown \langle\ell\rangle)=\rho_\delta(\ell)$, 
\sn
\item "{${{}}(d)$}"   if 
$g_{\bar \rho}(\eta) = 0$ and $\eta \frown \langle\ell\rangle) \in
\Gamma$ then $g_{\bar \rho}(\eta \char 94 \langle \ell \rangle) =
\alpha$.
\ermn
Let $A_{n,\bar \rho} = \{\eta\in \Gamma:g_{\bar \rho}(\eta) = \alpha$ and
$|\{\ell< \,{\text{\rm lg\/}}(\eta):
g_{\bar \rho}(\eta\restriction\ell)=\alpha\}|=n\}$,
so $A_{n,\bar \rho}$ is a front of $\Gamma$, above $A_{m,\bar \rho}$ for
$m < \omega$. Hence for each $m$ we have $p^* 
\Vdash ``{\underset\tilde {}\to \eta'}$ 
has an initial segment in $A_{m,\bar \rho}$". 
Define a $\Bbb Q$-name ${\underset\tilde {}\to h_{\bar \rho}}$ 
of a function from $\omega$ to $\omega$ by
${\underset\tilde {}\to h_{\bar \rho}}[G](n) =
\sup \text{ Rang}({\underset\tilde {}\to \eta^*} \restriction \ell)$
where $\ell$ is the unique natural number such that 
${\underset\tilde {}\to \eta'} \restriction \ell \in A_{n,\bar
\rho}$. 
Clearly ${\Cal I}_{n,\bar \rho} = \{p^*_\eta:\eta \in A_{n,\bar
\rho}\}$  is predense above $p^*$, so 
\footnote{ the ${\underset\tilde {}\to h_{\bar \rho}}$-s just clarify
here but are necessary in the proof of \scite{x.4} below.}
$\Vdash_{\Bbb Q} 
``{\underset\tilde {}\to h_{\bar \rho}} \in{}^\omega\omega$".

Now there is $N$ is a $\Bbb Q$-candidate, with 
$(\omega_1)^N$ not well ordered, and $p^*,\bar p^* \in N$ 
so $N \models$ ``all the above statement on $p^*,\rho \in
\Xi_\alpha$". \nl
[Why?  Let $N_1 \prec ({\Cal H}(\chi),\in)$ be such that $\Bbb
Q,p^*,\bar p^* \in N_1$ and $N_1 \prec N,N$ as above; $\omega^N =
\omega^{N_1}$ as in the Definition (see Keisler \cite{Ke71}).]

Choose $\alpha^*_n$ for $n < \omega$ such that $N \models 
``\alpha^*_{n+1}<\alpha^*_n$ are countable ordinals".  Without loss of
generality $N \models ``\alpha^*_n$ is a limit ordinal", and so for
some $\bar \rho$ we have $N \models ``\bar \rho \in \Xi_{\alpha^*_0}$".  So
clearly $N \models ``{\Cal I}_{n,\bar \rho}$ is predense above
$p^*$" for each $n$.  
Assume toward contradiction that the conclusion fails, so there is
$r^* \in {\Bbb Q}$ above $p^{**}$ which is $\langle N,
{\Bbb Q} \rangle$-generic. So
\mr
\item "{$\boxtimes_1$}"  ${\Cal I}_{n,\bar \rho}$ is predense 
above $r^*$ for each $n < \omega$. 
\ermn
There is in $\bold V$ (not in $N$!) a function 
$f_0:\Gamma \longrightarrow \omega$ such that 
$$
\align
\dsize \bigvee_{n<\omega}
\alpha^*_n \le^N g_{\bar \rho}(\eta) \Rightarrow (\forall k)
(k &\ge f(\eta) \and \eta^\frown <k> \in \Gamma\\
  &\rightarrow \dsize \bigvee_{n<\omega} \alpha^*_n \le^N 
g_{\bar \rho}(\eta^\frown \langle k \rangle)).
\endalign
$$
\mn
So by $\boxtimes$ (as $p^* \le_{\Bbb Q} r^*$) 
there is $q$ such that $r^* \le_{\Bbb Q} q$ 
and $q \Vdash_{\Bbb Q}$ ``for every $\ell$ large enough, 
${\underset\tilde {}\to \eta'}(\ell) > 
f_0({\underset\tilde {}\to \eta'} \restriction \ell)$".  
So $q$ forces that for some $\ell_0$ we have
$$
\ell_0 < \ell<\omega\ \Rightarrow \dsize \bigvee_n \alpha^*_n <^N
g_{\bar \rho}({\underset\tilde {}\to \eta'} \restriction \ell)^N 
\Rightarrow ({\underset\tilde {}\to \eta'} \restriction \ell) 
\notin \dbcu_n A_{n,\bar \rho}.
$$ 
\mn
Hence
\mr
\item "{$\boxtimes_2$}"  $q \Vdash_{\Bbb Q}$ ``the number of $n$ such
that $(\exists \ell)({\underset\tilde {}\to \eta'} \restriction \ell
\in A_{n,\bar \rho})$ is finite". 
\ermn
But $\boxtimes_1 + \boxtimes_2$ gives a contradiction.
\hfill$\square_{\scite{x.3}}$\margincite{x.3}
\enddemo
\bigskip

\demo{\stag{x.4} Proposition}  Assume that
\mr
\item "{$(a)$}"  $\Bbb Q$ is normal nep (need a weak version)
\sn
\item "{$(b)$}"   ${\underset\tilde {}\to \eta^*} \in {}^\omega \omega$ 
is a $\Bbb Q$--name, 
\sn
\item "{$(c)$}"  $p^* \in \Bbb Q$ forces ${\underset\tilde {}\to \eta^*}$ 
is a dominating real, 
\sn
\item "{$(d)$}"  ``$\{p_1,\ldots,p_n\}$ is predense over $\{q,q_1\}$, 
$n<\omega$'' as well as ``$p \in \Bbb Q$", ``$p \le_{\Bbb Q} q^*$" are 
upward absolute from $\Bbb Q$-candidates in the sense of Definition
\scite{xa.1}(1), so not necessarily well founded.
\ermn
\ub{Then} $p^*$ forces that ${\frak b}=\aleph_1$ in $\bold V^{\Bbb Q}$.  
\enddemo
\bigskip

\remark{\stag{x.4a} Remark}  Concerning clause $(d)$, of 
course, ``${\Cal I}$ is predense above
${\Cal J}$" and ${\Cal I},{\Cal J} \subseteq \Bbb Q$ means that: if $p
\in \Bbb Q$ is above every $q \in {\Cal J}$ \ub{then} $p$ is
compatible with some $r \in {\Cal I}$.
\endremark
\bigskip

\demo{Proof}  It is enough to prove that some condition above $p^*$
forces ${\frak b} = \aleph_1$.

By the nep \wilog \, there are $p^{**}$ and $\langle {\Cal J}^*_n:n <
\omega \rangle$ as in the beginning of \scite{x.3}.
Let ${\underset\tilde {}\to \eta^*},\Gamma,
\bar p^* =\langle p^*_\eta:\eta\in\Gamma\rangle,
{\underset\tilde {}\to \eta'},\Xi_\alpha$ for $\alpha \le \omega_1$
and $g_{\bar \rho},
{\underset\tilde {}\to h_{\bar \rho}}$ for $\bar \rho \in
\Xi_\alpha,\alpha < \omega_1$ be as in the
proof of \scite{x.3}.  Choose $\bar \rho \in \Xi_{\omega_1}$ and let 
$\bar \rho_\delta = \bar \rho \restriction (\delta +1)$ 
for limit $\delta<\omega_1$. If 

$$
p^{**} \Vdash ``\{{\underset\tilde {}\to h_{{\bar \rho}_\delta}}:
\delta < \omega_1 \text{ limit}\} \text{ is not dominated}",
$$
\mn
we are done. So assume not, hence for some ${\underset\tilde {}\to
h^*}$ and $q^*$, we have $p^{**} \le_{\Bbb Q} q^*$ and 

$$
q^* \Vdash_{\Bbb Q} ``{\underset\tilde {}\to h^*} \in 
{}^\omega\omega \text{ dominates } 
\{{\underset\tilde {}\to h_{{\bar \rho}_\delta}}:
\delta < \omega_1 \text{ limit }\}".
$$
Without loss of generality ${\underset\tilde {}\to h^*}$ is hc
(hereditarily countable) $\Bbb Q$-name above $q^*$, more 
specifically as before for each $n < \omega$
we have $\langle r^*_{n,\ell}:\ell<\omega\rangle$
an antichain of $\Bbb Q$ predense over $q^*$ such that 
$r^*_{n,\ell} \Vdash_{\Bbb Q} ``{\underset\tilde {}\to h^*}(n) =
k_{n,\ell}"$. So ${\underset\tilde {}\to h^*}$
is 

$$
\langle (n,\ell,{\underset\tilde {}\to r^*_{n,\ell}},k_{n,\ell}):n,\ell<
\omega \rangle.
$$ 
\mn
Choose a $\Bbb Q$-candidate $M$ to which 
${\underset\tilde {}\to h^*},\bar \rho,p^*,p^{**},\bar p^*$ belong 
(e.g., $M \prec ({\Cal H}(\chi),\in)$ is countable in
normal cases), and choose a countable elementary extension $N$ of $M$ such
that in $N$ there are $\alpha^*_n$ for $n<\omega$ as in the proof of
\scite{x.3}.  

So in $N,\bar \rho = \bar \rho^N_{\alpha^*_0}$ is well defined, 
so as $M \prec N$ there are $n^*,r^*$ such that 
\mr
\item "{$(*)_0$}"  $N \models ``r^* \in {\Bbb Q}$ and
$q^* \le_{\Bbb Q} r^*$ and $n^*<\omega$ and $r^*$ forces (for
$\Vdash_{\Bbb Q}$) that ${\underset\tilde {}\to h_{\bar \rho}}
\restriction [n^*,\omega) < {\underset\tilde {}\to h^*} 
\restriction [n^*,\omega)$". 
\ermn
Let $n \in [n^*,\omega)$ and $\ell < \omega$.  Recalling that (in
$\bold V$ hence in $M$ hence in $N$) we have $r^*_{n,\ell}
\Vdash_{\Bbb Q}
``{\underset\tilde {}\to h^*}(n) = k_{n,\ell}"$ and recalling that every
$\eta \in \Gamma$ is strictly increasing and the definition of
$h_{\bar \rho}$ in $N$, clearly
\mr
\item "{$(*)^1_{n,\ell}$}"  the set 
$A_{n,\ell} = \{\nu \in A_{n,{\bar \rho}}$: max Rang$(\nu) < k_{n,\ell}\}$ 
is finite (we get the same set
in $N$ and in $\bold V$).
\ermn
Hence (using $(*)_0$)
\mr
\item "{$(*)^2_{n,\ell}$}"  in $N$ the set $\{p_\nu:\nu \in
A_{n,\ell}\}$ is predense (in $\Bbb Q^N$) above
$\{r^*,r^*_{n,\ell}\}$.
\ermn
But by the clause (d) of the assumption this amount of predensity is
upward absolute hence
\mr
\item "{$(*)^3_{n,\ell}$}"  in $\bold V$ the set $\{p_\nu:\nu \in
A_{n,\ell}\}$ is predense (in $\Bbb Q$) above
$\{r^*,r^*_{n,\ell}\}$.
\ermn
But $q^* \le_{\Bbb Q} r^*$ and $\{r^*_{n,\ell}:\ell < \omega\}$ is
predense (in $\Bbb Q$) above $q^*$, hence
\mr
\item "{$(*)^4_n$}"  for each $n < \omega$ in $\bold V$ the set $\dbcu_{\ell < \omega}
\{p_\nu:\nu \in A_{n,\ell}\}$ is predense in $\Bbb Q$ above $r^*$.
\ermn
Now $\dbcu_{\ell < \omega} \{p_\nu:\nu \in A_{n,\ell}\} \subseteq \{p_\nu:\nu
\in A_{n,\bar \rho}\} = {\Cal I}_{n,\bar \rho}$; note that $N \models
``{\Cal I}_{n,\bar \rho}$ is a countable subset of $\Bbb Q$" hence
${\Cal I}_{n,\bar \rho}$ is a countable subset of $\Bbb Q$ so
\mr
\item "{$(*)_5$}"  for every $n,{\Cal I}_{n,\bar \rho} = 
{\Cal I}^N_{n,\bar \rho}$ is predense in $\Bbb Q$ above $r^*$. \nl
[Why?  If $G \subseteq \Bbb Q$ is generic over $\bold V$ and $r^* \in
G$ then $q^* \in
G$ (as $q^* \le_{\Bbb Q} r^*$) and for some $\ell,r^*_{n,\ell} \in G$
(as $\{r^*_{n,\ell}:\ell < \omega\}$ is predense over $r^*$).  But
$\{p_\nu:\nu \in A_{n,\ell}\}$ is predense over $\{r^*,r^*_{n,\ell}\}$
in $\bold V$ then for some $\nu \in A_{n,\ell}$ we have $p_\nu \in G$,
but $p_\nu \in A_{n,\ell} \Rightarrow p_\nu \in {\Cal I}_{n,\rho}$ as
said earlier, so we are done proving $(*)_5$.]
\ermn
This means that in the proof of \scite{xa.1} the statement
$\boxtimes_1$ holds and continues as in the proof of \scite{x.3}.
\hfill$\square_{\scite{x.4}}$\margincite{x.4}
\enddemo
\bigskip

\proclaim{\stag{1.8} Claim}   Amoeba forcing forces ${\frak b}=\aleph_1$, 
and similarly dominating real 
forcing (= Hechler forcing) and universal meagre forcing.
\endproclaim
\bigskip

\demo{Proof}  We like to apply \scite{x.4}.  Let $\Bbb Q$ be an amoeba
forcing which is $\{T:T \subseteq {}^{\omega >} 2$ is non empty closed
under initial segments and Leb(lim$(T)) > 1/2\}$, ordered by inverse
inclusion; note that for notational simplicity we allow trees with
maximal nodes.

Clearly $p \in \Bbb Q,p \le_{\Bbb Q} q$ 
are Borel relations and any $p,q \in \Bbb
Q$ has a l.u.b.: $p \cap q$, ``$p,q$ compatible" is Borel.  
The main point is to show that
``$\{p_\ell:\ell < n\}$ is predense above $\{q_1,q_2\}$" is upward
absolute for nw-candidates; we can replace $\{q_1,q_2\}$ by $\{q\}$
where $q = q_1 \cap q_2$.  Define for $m,k < \omega$ but $> 0$ and $s
\subseteq q_1 \cap q_2 \cap {}^m 2$:

$a_{m,k}(s) = \text{ Min}\{|s \cap p_\ell|/2^m:\ell < n\}$, this is a 
real number $\in [0,1]$ and

$$
a_{m,k} = \text{ Min}\{a_{n,k}(s):s \subseteq q_1 \cap q_2 \cap {}^m 2
\text{ and } |s|/2^m \ge \frac 12 + \frac 1k\}.
$$
\mn
Note that if $s$ exemplifies $a_{n,k} \le b$ then $s^+ = \{\sigma
\char 94 <\ell>:\ell = 0,1 \text{ and } \nu \in s\}$ exemplifies
$a_{n,k} \le b$ hence
\mr
\item "{$(*)$}"  $a_{m,k} \ge a_{m+1,k}$.
\ermn
We shall show that the following statements are equivalent:
\mr
\item "{$(\alpha)$}"  there is $r \in \Bbb Q$ above $q$ incompatible
with $p_0,\dotsc,p_{n-1}$
\sn
\item "{$(\beta)$}"  for some $r \in \Bbb Q$ we have Leb(lim$(p_\ell
\cap r)) \le \frac 12 - \frac 1k$ for $\ell < n$ and Leb(lim$(r)) >
\frac 12 + \frac 1k$ for some $k \in (0,u)$
\sn
\item "{$(\gamma)$}"  lim$\langle a_{m,k}:m < \omega \rangle \le 
\frac 12 - \frac 1k$ for some $k \in (0,\omega)$.
\ermn
If $(\alpha)$ holds, let $r$ exemplify it, so for some $\varepsilon_1
> 0$, Leb(lim$(r)) > \frac 12 + \varepsilon_1$, and Leb(lim$(p_\ell
\cap r)) \le \frac 12$ for $\ell < n$.  We cn find, for $\ell < n$, a
clopen subset $B_\ell$ of ${}^\omega 2$ such that Leb(lim$(p_\ell \cap
r) \cap B_\ell) > 0$, Leb$(B_\ell) < \frac{\varepsilon_1}{n+3}$.  
Let $r' = \{\eta
\in r$: there is $\rho \in {}^\omega 2 \backslash \dbcu_\ell B_\ell$
above $\eta\}$, and $k$ be large enough, they exemplify $(\beta)$.

If $(\beta)$ holds, exemplified by $r,k$, \ub{then}
$\ell < n \Rightarrow \text{ Leb}(\text{Lim}(p_\ell \cap r) \le \frac
12$ -  hence by the definition of Leb measure
\mr
\item "{$(*)$}"  $\frac 12 - \frac 1k \ge \text{ Lim}\langle
|p_\ell \cap r \cap {}^m 2|/2^m:m < k \rangle$ \nl
hence
\nl
$\frac 12 - \frac 1k \ge \text{ Lim}\langle 
\underset {\ell < n} {}\to {\text{Max}}|p_\ell 
\cap r \cap {}^m2|/2^m:m < \omega \rangle$.
\ermn
But $a_{m,k} \le a_{m,k}(r \cap {}^m 2)$ because $|r \cap {}^m 2|/2^m
\ge \text{ Leb(Lim } r) \ge \frac 12 + \frac 1k$ and $a_{m,k}(r \cap
{}^m 2) = \underset {\ell < n} {}\to {\text{Max}}(|p_\ell \cap r \cap
{}^m 2|/2^m)$. \nl
Putting together those inequalities and $(*)_1$ we have $\frac 12 >
\text{ lim} \langle a_{m,k}:m < \omega \rangle$ as required, so
$(\gamma)$ holds, i.e. we have proved $(\beta) \Rightarrow (\gamma)$.

Lastly, assume $(\gamma)$ and we shall prove $(\alpha)$.  For each $m$
let $s_m \subseteq q \cap {}^m2$ be such that $a_{m,k}(s_m) =
a_{m,k}$.  Let $m$ be large enough such that $a_{m,k} < \frac 12 -
\frac 1{4k}$ and $|r \cap {}^m 2|/2^m - \text{ Leb(lim }r) < 1/4k$.
Let $r = \{\rho \in q$: if $\ell g(\rho) \ge m$ then $\rho \in s_m\}$.

Clearly $r \subseteq q$ is a subtree, and

$$
\align
\text{Leb}(\text{lim } r) &\ge \text{ Leb}(\text{lim}(q)) -
\text{ Leb}\{\eta \in {}^\omega 2:\eta \in \text{ lim}(q) \text{ but }
\eta \restriction m \in a_{m,k})\} \\
  &\ge \text{ Leb}(\text{lim}(q)) - (|q \cap {}^m 2|/2^m -
|a_{m,k}|(2^m) \\
  &\ge \text{ Leb}(\text{lim}(q)) - ((\text{Leb}(\text{lim}(q) +
1/4^k) - (a_{m,k})/2^m) \\
  &= |a_{m,k}|/2^m + -1/4k \ge \frac 12 + \frac 1k - 1/4k > \frac 12 - .
\endalign
$$
\mn
So $r \in \Bbb Q$, also for $\ell < n$, the conditions $r,p_\ell$ are
incompatible as Leb(lim$(p_\ell \cap r)) = \text{
Leb}(\text{lim}(\{\eta \in p_\ell$: if $\ell g(\eta) \ge n$ then $\eta
\restriction m \in a_{m,k})\}) \le \text{ Leb}(\{\eta \in {}^\omega
2:\eta \restriction m \in p_\ell \cap {}^m 2 \cap a_{m,r}\}) = p_\ell
\cap a_{m,k}|/2^{< m} \le \frac 12 - \frac 1k < \frac 12$. \nl
So we have finished proving $(\gamma) \Rightarrow (\alpha)$ hence
proving $(\alpha) \Leftrightarrow (\beta) \Rightarrow \delta)$.
\hfill$\square_{\scite{1.8}}$\margincite{1.8}
\enddemo
\bigskip

\demo{\stag{x.6} Conclusion}  1) For no $\aleph_1$-complete ideal $J$ on a
set $X$ is the Boolean algebra ${\Cal P}(X)/J$ isomorphic to the
Boolean algebra of the Amoeba forcing. \nl
2) The following is impossible
\mr
\item "{$(a)$}"  $J$ is a $(< \kappa)$-complete ideal on a set $X$ not
$\kappa^+$-complete 
\sn
\item "{$(b)$}"  ${\Cal P}(X)/J$ is isomorhpic to the Boolean algebra
of the forcing notion $\Bbb Q$ which satisfies the $\kappa^+$-c.c. say
$g:\Bbb Q \rightarrow {\Cal P}(X)/J$
\sn
\item "{$(c)$}"  forcing with $\Bbb Q$ adds a dominating real
\sn
\item "{$(d)$}"  forcing with $\Bbb Q$ makes ${\frak b} \le \kappa$.
\endroster
\enddemo
\bigskip

\demo{Proof}  1) Follows by part (2) for $\kappa = \aleph_1$
below and proposition
\scite{x.4}. \nl
2) Let $\kappa_1$ be maximal such that $J$ is $(< \kappa_1)$-complete,
so $J$ is not $(< \kappa^+_1)$-complete, now replacing $\kappa$ by
$\kappa_1$, clauses (a)-(d) are still satisfied, so without loss of
generality  $J$ is not $\kappa^+$-complete.   The proof is close to 
\cite[3.1]{GiSh:357} but see \cite{GiSh:412} but we give a
self contained proof.

Let $G \subseteq \Bbb Q$ be generic over $\bold V$, so we can define
$\underset\tilde {}\to D[G] = \{Y \subseteq X:g(p) \subseteq Y$ mod
$J$ for some $p \in G\}$ is an ultrafilter on $X$, i.e., on the
Boolean algebra ${\Cal P}(X)^{\bold V}$ disjoint to $J$.  As
${\Cal P}(X)/J$ satisfies the $\kappa^+$-c.c. and $J$ 
is $\kappa$-complete (in $\bold V$) 
clearly the ultrapower
$\bold V^X/\underset\tilde {}\to D[G] = \{f/\underset\tilde {}\to
D[G]:f \in {}^X \bold V$ is from $\bold V\}$ is well founded so we
identify it with its Mostowski collapse $M$.  
Let $\bold j$ be the
natural elementary embedding of $\bold V$ into $M$.
Clearly in $\bold V[G]$ the model $M$ is closed
under taking sequences of length $\le \kappa$.    
In particular $M$ contains all $\omega$-sequences of natural numbers
from $\bold V[G]$ hence $({}^\omega \omega)^M = 
({}^\omega \omega)^{\bold V[G]}$.
As $M$ contains all $\le \kappa$-sequences of reals from $\bold
V[G]$ and $\bold V[G] \models {\frak b} \le \kappa$, clearly in $\bold V[G]$
there is a sequence $\bar f = \langle f_\alpha:\alpha < \theta
\rangle$ exemplifying ${\frak b} = \theta \le \kappa$, 
hence $\bar f \in M$.  So
necessarily $M \models ``{\frak b} \le \theta"$ but $M \models \theta
\le \kappa < \bold j(\kappa)$ hence by Los theorem
also $\bold V \models {\frak b} \le \theta$.  So let $\bar f$ be such that
$\bold V \models ``\bar f
= \langle f_\alpha:\alpha < \theta \rangle$ exemplify ${\frak b} \le
\theta"$, hence as $\theta < \kappa$, clearly $\bold j(\bar f) = \bar f$, so
$\{f_\alpha:\alpha < \theta\} \subseteq ({}^\omega \omega)^M$ is
unbounded in $({}^\omega \omega)^M$ but the latter is $({}^\omega
\omega)^{\bold V[G]}$ in which there is a $\underset\tilde {}\to \eta
\in {}^\omega \omega$ dominating $({}^\omega \omega)^{\bold V}$ hence
$\bar f$, contradiction. \hfill$\square_{\scite{x.6}}$\margincite{x.6}
\enddemo
\newpage

\head {\S2} \endhead  \resetall \sectno=2
\bn
In \cite{Sh:480} it is proved that if $\Bbb Q$ is a
Souslin-c.c.c. forcing adding a non-dominated real, \ub{then} it adds
a Cohen real.  We here try to extend the result to nep forcing.

More fully we use the following: let $N$ be a countable elementary
submodel of $({\Cal H}(\chi),\in)$ to which $\Bbb P,\Bbb Q$ belongs,
$G$ a subset of $\Bbb P^N$ generic over $N$ \ub{then} $N[G]$ is a
$\Bbb Q$-candidate. 
It may be clearer to let $M$ be the ordinal collapse of $N,\bold j:N
\rightarrow M$ the isomorphism and demand $M[\bold j'' G]$ is a $\Bbb
Q$-candidate. 
\bigskip

\proclaim{\stag{y.1} Claim}  1) Assume
\mr
\item "{$(a)$}"  $\Bbb Q$ is nep-forcing and is c.c.c. 
\sn
\item "{$(b)$}"  $\Vdash_{\Bbb Q} ``\underset\tilde {}\to f \in
{}^\omega \omega$ is not dominated by any old $f \in {}^\omega \omega"$
\sn
\item "{$(c)$}"  $\Bbb P$ is a forcing notion adding a dominating
$\underset\tilde {}\to g \in {}^\omega \omega$
\sn
\item "{$(d)$}"  forcing with $\Bbb P$ preserves \footnote{i.e., in
$\bold V$ if $N$ is a $\Bbb Q$-candidate, $G \subseteq \Bbb Q^N$ is
generic over $N$ so $G \in \bold V$, then \nl
$N<G>$ is a $\Bbb Q$-candidate; we can use transitive candidates; stated to
ease keeping track of the assumptions needed} being a $\Bbb Q$-candidate (it
follows also that ``$q \in \Bbb Q$", ``$p \le_{\Bbb Q} q$" are preserved).
\ermn
\ub{Then} in $\bold V^{\bold P}$
\mr
\item "{$\circledast_1$}"  forcing with $\Bbb Q$ add a Cohen real.
\ermn
2) Assume we replace above clause (b) by
\mr
\item "{$(b)'$}"  $\Vdash_{\Bbb Q} ``\underset\tilde {}\to \eta
\in{}^\omega 2$ is not equal to any old member of ${}^\omega 2"$.
\ermn
\ub{Then} (in $\bold V^{\Bbb P}$)
\mr
\item "{$\circledast_2$}"  for some strictly increasing sequence
$\langle n_i:i < \omega \rangle$ we have: \nl
$2^i \le |\{\eta \in {}^{n_i} 2$: some $q'$ above $q$ forces that
$\eta = \underset\tilde {}\to \eta \restriction n_i\}|$.
\endroster
\endproclaim
\bigskip

\demo{Proof}  For part (1) let $\bold t$ be 1,
${\underset\tilde {}\to f^{\bold t}} = \underset\tilde {}\to
f$ and for part (2) let $\bold t=2$ and 
${\underset\tilde {}\to f^{\bold t}} =
\underset\tilde {}\to \eta$.
So ${\underset\tilde {}\to f^{\bold t}}$ is actually $\langle
(r^*_{n,\ell},k_{n,\ell}:n < \omega,\ell < \omega \rangle$ where
$r^*_{n,\ell} \Vdash {\underset\tilde {}\to f^{\bold t}}(n) 
= k_{n,\ell}$ and $\langle r^*_{n,\ell}:\ell < \omega \rangle$ 
a maximal antichain of $\Bbb Q$; similarly $\underset\tilde {}\to
\eta$.   Without loss of generality ${\underset\tilde {}\to f^0},
{\underset\tilde {}\to f^1}$ (forced to be) ${\underset\tilde {}\to g}$ are 
strictly increasing; note that 
for ${\underset\tilde {}\to f^{\bold t}},\bold t =1$ this just
means that $n_1 < n_2 \and k_{n_1,\ell_1} \ge k_{n_2,\ell_2} \Rightarrow
(r^*_{n_1,k_1},r^*_{n_2,k_2}$ are incompatible) so is absolute enough.
Let ($\chi$ be strong limit and)
$N\prec ({\Cal H}(\chi),\in),N$ be countable such that $\{\Bbb
P,\underset\tilde {}\to g,q^*,\Bbb Q\} \in N$ 
and ${\underset\tilde {}\to f^{\bold t}} \in N$, 
i.e., $\langle (r^*_{n,\ell},k_{n,\ell}):n < \omega,\ell <
\omega \rangle$ belongs to $N$.  Now
\mr
\item "{$(*)_1$}"  $N \models ``\Bbb P$ is a forcing notion,
$\underset\tilde {}\to g$ is a $\Bbb P$-name of an increasing 
member of ${}^\omega \omega$ dominating all old ones". 
\ermn
Observe:
\mr
\item "{$(*)_2$}"  if $M$ is a $\Bbb Q$-candidate, $\bar r = \langle
r_\ell:\ell < \omega \rangle$ is a maximal antichain of $\Bbb Q$ and $\bar
r \in M,r_\ell \in \Bbb Q^M$, \ub{then} $M \models ``\bar r$ is a
maximal antichain of $\Bbb Q$". \nl
[Why?  First, if $n < m < \omega,M \models ``r_n,r_m$ are compatible in
$\Bbb Q$" let $r \in \Bbb Q^M$ be a common upper bound by $\le^M_{\Bbb Q}$,
it is a common upper bound in $\Bbb Q$, contradiction.  Second, if $N
\models ``q \in \Bbb Q$ is incompatible with each $r_\ell"$, let $q_1 \in
\Bbb Q$ be $(M,\Bbb Q)$-generic such that $q \le q^*_1$.  But $q_1$ is
necessarily $\le_{\Bbb Q}$-compatible with $r_n$ for some $n$ so for
some $q_2$ we have $r_n \le_Q q_2 \and q_1 \le_{\Bbb Q} q_2$, so $q_2
\Vdash ``\{q_1,r_n\} \subseteq {\underset\tilde {}\to G_{\Bbb Q}}
\cap M"$ but $q_1 \le_{\Bbb Q} q_2$ and $q_1 \Vdash_{\Bbb Q} ``
{\underset\tilde {}\to G_{\Bbb Q}} \cap N$ is $\le_{\Bbb
Q}^M$-directed, contradiction.] 
\ermn
Continuing $(*)_1$, for $\bold t =1$:
\mr
\item "{$(*)_3$}"  $N \models$ ``forcing with $\Bbb P$ preserve the
property of $(\Bbb Q,{\underset\tilde {}\to f^{\bold t}}$), 
i.e., ${\underset\tilde {}\to f^{\bold t}}$ not dominated". \nl
[Why?  First being a $\Bbb Q$-name of a member of
${}^\omega \omega$ is preserved after forcing with $\Bbb P$ follows:
by $(*)_2 +$ assumption $(d)$; if not, then we can find $G \subseteq
\Bbb P^N,G \in \bold V$ generic over $\bold V$, let \footnote{can use
the ord collapse of $M$} $M = N[G]$ it is a $\Bbb Q$-candidate and we
can apply $(*)_2$.  Second, assume toward contradiction that $(*)_3$
fails.
Let $p^* \in \Bbb P^N$ force the negation (in $N$) and
choose, in $\bold V$ a set $G \subseteq \Bbb P^N$
generic over $N$ to which $p^*$ belongs so $N \subseteq N[G] \in \bold
V,\bold V \models ``N[G]$ a $\Bbb Q$-candidate".  As the conclusion of
$(*)_3$ fails we can find $q_1 \in \Bbb Q^{N[G]}$ and $h \in ({}^\omega
\omega)^{N[G]}$ such that $N[G] \models ``q_1 \in \Bbb Q,q^* \le_{\Bbb
Q} q_1$ and $q_1$ forces ($\Vdash_{\Bbb Q}$) that 
$\underset\tilde {}\to f \le h$".  Let
$q_2$ be $(N[G],\Bbb Q)$-generic condition satisfying 
$q_1 \le_{\Bbb Q} q_2$ hence $q_2
\Vdash_{\Bbb Q} ``\underset\tilde {}\to f \le g \in {}^\omega
\omega"$, contradicting the choice of $q^*,\underset\tilde {}\to f$.]
\sn
\item "{$(*)_4$}"  also in $\bold V^{\Bbb P},{\underset\tilde {}\to
f^{\bold t}}$ is not dominated if $\bold t=1$. \nl
[Why?  As $N \prec ({\Cal H}(\chi),\in)$.]
\ermn
If $\bold t = 1$, let 
the $\Bbb Q$-name ${\underset\tilde {}\to \eta^{\bold t}} 
\in {}^\omega 2$ be such that for every $\ell < \omega$ we have
$\underset\tilde {}\to \eta(\ell) = 1 \Leftrightarrow \ell 
\in \text{ Rang}(\underset\tilde {}\to  f)$ and 
${\underset\tilde {}\to \eta^{\bold t}} = \underset\tilde {}\to \eta$.

Without loss of generality
\mr
\item "{$(*)_5$}"  for every $h:{}^{\omega >} \omega \rightarrow
\omega$ from $\bold V$ we have $\Vdash_{\Bbb P} (\forall^\infty n)
\underset\tilde {}\to g(n) > h(\underset\tilde {}\to g \restriction n)$. \nl
[Why?  As, e.g., we can replace ${\underset\tilde {}\to \eta^*}$ by 
$\underset\tilde {}\to \nu$, where
$\underset\tilde {}\to \nu(0)=1$ and
$\underset\tilde {}\to \nu(n+1) = {\underset\tilde {}\to \eta^*}
(\underset\tilde {}\to \nu(n))$, note that $\underset\tilde {}\to \nu$
is strictly increasing as ${\underset\tilde {}\to \eta^*}$ is and
$\underset\tilde {}\to \nu(0) > 0$.]
\sn
\item "{$(*)_6$}"  $N \models \Vdash_{\Bbb P} ``{\underset\tilde {}\to
\eta^{\bold t}} \in {}^\omega 2$ is new". \nl
[Why?  For $\bold t = 1$ by $(*)_3$, for $\bold t=2$ even easier.]
\ermn
For $q \in \Bbb Q$ let

$$
T_q[{\underset\tilde {}\to \eta^{\bold t}}] = \{\nu \in {}^{\omega >} 2:q
\nVdash_{\Bbb Q} \nu \ntriangleleft \underset\tilde {}\to \eta\}
$$
so $T_q[{\underset\tilde {}\to \eta^{\bold t}}]$ is a non empty subtree of
${}^{\omega >}2$ with no maixmal nodes. \nl
We say that $u \in [\omega]^{\aleph_0}$ is large 
for $(q,{\underset\tilde {}\to \eta^{\bold t}})$ \ub{if} for 
every $r \in \Bbb Q$ above $q$, 
we have $\otimes_{r,u} = \otimes^{\underset\tilde {}\to \eta}_{r,u}$
the following holds:
\sn
\ub{Case $\bold t = 1$}:  For some $n^* \in u$ for every $n^* < m \in
u,n,m \in u$ such that $n^* < n < m$ for some $\nu \in {}^m 2 \cap T_r
[{\underset\tilde {}\to \eta^{\bold t}}]$ we have $\ell \in [n^*,m)
\Rightarrow \nu(\ell) = 0$ but $(\exists \ell)(n \le \ell < m \and
\eta(\ell)=1)$. 
\bn
\ub{Case $\bold t =2$}:  For some $n^* \in u$ if $n,m \in u$ and $n^*
\le n < m$ then for some $\nu_1,\nu_2 \in {}^m 2 \cap T_r
[{\underset\tilde {}\to \eta^{\bold t}}]$, and for $\ell \in (n,m)$ we
have $\nu_1 \restriction \ell = \nu_2 \restriction \ell,\nu_1(\ell)
\ne \nu_2(\ell)$.
\bn
Let ${\underset\tilde {}\to g^*} \in {}^\omega \omega$, a $\Bbb
P$-name, be defined by ${\underset\tilde {}\to g^*}(0) = 0,
{\underset\tilde {}\to g^*}(n+1) = \underset\tilde {}\to g
(n+1+{\underset\tilde {}\to g^*}(n))$.
\mn
\ub{Subclaim}:  Let $q^* \in \Bbb Q^N$ then $N \models \Vdash_{\Bbb P}$ ``some
$u \in [\omega]^{\aleph_0}$ is large for
$(q^*,{\underset\tilde {}\to \eta^{\bold t}})$". \nl
It will take us awhile.  Let $\underset\tilde {}\to u =
\text{ Rang}({\underset\tilde {}\to g^*})$, a $\Bbb P$-name in $N$,
and assume that $\underset\tilde {}\to u$ is not as required, so
for some $p^* \in \Bbb P^N$ and $\Bbb P$-names 
$\underset\tilde {}\to q,n$ we have \nl
$N \models ``p^* \Vdash_{\Bbb P} [\underset\tilde {}\to q \in \Bbb Q$ is above
$q^*$ and $\neg \otimes_{\underset\tilde {}\to q,\underset\tilde {}\to
u}$.]
\nl
Let $G(\in \bold V)$ be a subset of $\Bbb P^N$ generic over $N$ such
that $p^* \in G$
\mr
\item "{$(*)_7$}"  if $\bold t = 1$ then for some $n_* < \omega$: \nl
for every $m \in (n_*,\omega)$ we have 
$N[G] \models \underset\tilde {}\to q[G]
\nVdash_Q (\exists \ell)(n_* \le \ell < m 
\and {\underset\tilde {}\to \eta^{\bold t}}(\ell) = 1)$; note that
therefore for every $m > n_*$ there is an $s \in T_{\underset\tilde
{}\to q[G]}({\underset\tilde {}\to \eta^{\bold t}})$ such that for
every $\ell \in [n_*,m)$, we have $s(\ell)=0$. \nl
[Why?  By $(*)_3$.] 
\ermn
So in $N[G]$ we define:
\sn
\ub{Case $\bold t = 1$}:  Letting $n_*$ be as in $(*)_7$, we let

$$
\align
u' = \{n \in \underset\tilde {}\to u[G]:&n > n_* \text{ and letting }
m = \text{ Min}(\underset\tilde {}\to u[G] \backslash (n+1)) 
\text{ \ub{there is}} \\
  &\eta \in {}^m 2 \cap T_r[{\underset\tilde {}\to \eta^{\bold t}}
\text{ such that} \\
  &\ell \in [n_*,n) \Rightarrow \eta(\ell) = v \text{ but }
(\exists \ell)(n \le \ell < m \and \eta(\ell) = 1).
\endalign
$$
\mn
Note that by $(*)_7$, $u'$ is infinite.
\bn
\ub{Case $\bold t=2$}:
$$
\align
u' = \{n \in \underset\tilde {}\to u[G]:&\text{letting } m = 
\text{ Min}(\underset\tilde {}\to u[G] \backslash (n+1)) 
\text{ \ub{there are} no} \\
  &\nu_1 \ne \nu_2 \in {}^m 2 \cap T_r[{\underset\tilde {}\to
\eta^{\bold t}}] \text{ such that } \nu_1 \restriction n = \nu_2
\restriction n\}.
\endalign
$$
\mn
  So by our present
assumptions, $N[G] \models ``u'$ is infinite". \nl
Let $n_i$ be the $i$-th member of $u'$ and $m_i = \text{ Min}(u
\backslash (n_i+1))$.  Let $u' = {\underset\tilde {}\to u'}[G],n_i =
{\underset\tilde {}\to n_i}[G],m_i = {\underset\tilde {}\to m_i}[G]$
for some $\Bbb P$-name 
${\underset\tilde {}\to u'} \in N$ and a sequence $\langle
{\underset\tilde {}\to n_i},{\underset\tilde {}\to m_i}:i < \omega
\rangle \in N$ of $\Bbb P$-name.

Without loss of generality 
$p^* \in \Bbb P$ is such that it forces all the above in
particular that ${\underset\tilde {}\to u'} \subseteq \omega$ is
infinite.

Let $\langle {\Cal I}_n:n < \omega \rangle$ list the dense open
subsets of $\Bbb P$ which belongs to $N$.  Let $\langle J_k,\rho_k:k <
\omega \rangle$ be such that: $J_k$ is a finite front of ${}^{\omega
>}2,J_0 = \{<>\},\rho_k \in J_k,J_{k+1} 
= (J_n \backslash \{\rho_k\}) \cup \{\rho_k
\char 94 <0>,\rho_k \char 94 <1>\}$ and $n < \omega \and \rho \in J_n
\Rightarrow (\exists m \ge k)(\rho_m = \rho)$.  We choose $\bar p^k =
\langle p_\rho,m_\rho,n_\rho:\rho \in J_k \rangle$ by induction on
$k$ such that (so choosing $\bar p^k$ we have already chosen $\bar p^{k+1}
\restriction (J_{k+1} \backslash \{\rho_k \char 94 <0>,\rho_k \char 94 <1>\})$:
\mr
\item "{$(a)$}"  $p^* \le_{\Bbb P} p_\rho \in N$
\sn
\item "{$(b)$}"   $m_\rho < n_\rho$
\sn
\item "{$(c)$}"  if $m > n_\rho$ then for some $q$ we have $p_\rho
\le_{\Bbb P} q \in N$ \nl
$q \Vdash (\exists i)({\underset\tilde {}\to n_i} \le n_\rho \wedge
{\underset\tilde {}\to m_i} > m)$
\sn
\item "{$(d)$}"  $p_{\rho_k} \le_{\Bbb P} p_{\rho_k \char 94 <\ell>}
\in {\Cal I}_{\ell g(\rho_k)}$ for $\ell=0,1$
\sn
\item "{$(e)$}"  $p_{\rho \char 94 <\ell>} \Vdash ``(\exists i)(
{\underset\tilde {}\to n_i} \le n_\rho \and {\underset\tilde {}\to
m_i} > m_{\rho \char 94 <\ell>})"$
\sn
\item "{$(f)$}"   $m_{\rho_k \char 94 <\ell>} > \sup\{n_\nu:\nu \in J_k\}$.
\ermn
Let us carry the induction.

In step $k=0$ let $p_{<>} = p^*,n_{<>}$ is chosen as below, $m_{<>}$
is immaterial.  If we have defined for $k$, first choose $m_{\rho_k
\char 94 <\ell>}$ to satisfy clause (f), then choose $p'_{\rho_k
\char 94 <\ell>} \ge p_{\rho_k}$ to satisfy clauses (e), possible by
clause (c) and choose $p_{\rho_k \char 94 <\ell>} \ge p'_{\rho_k \char
94 <\ell>}$ to satisfy clause (d).  Lastly, choose $n_{\rho_k \char 94
<\ell>}$ to satisfy clause (c); this is possible by the observation
below (for $\bold t=0$ we use $(*)_7$). 

For each $\rho \in {}^\omega 2$ let $G_\rho = \{p \in \Bbb P^N:p
\le_{\Bbb P} p_{\rho \restriction \ell}$ for some $\ell < \omega\}$, clearly
$G_\rho(\in \bold V)$ is a subset of $\Bbb P^N$ generic over $N$ (by
clause (d)).  Now 
$q_\rho = \underset\tilde {}\to q[G_\rho]$ and $T_\rho =
T_{q_\rho}[{\underset\tilde {}\to \eta}]$, are well defined in $N[G_\rho]$
hence in $\bold V$.  It is easy to see that
\mr
\item "{$\circledast_1$}"  if $\nu_1 \ne \nu_2 \in {}^\omega 2 \and
\nu_1 \restriction k \ne \nu_2 \restriction k$ and $\eta \in T_{\nu_1}
\cap T_{\nu_2}$ and $\ell g(\eta) > n_{\nu_1 \restriction k},n_{\nu_2
\restriction k}$, \ub{then} $\eta$ has at most one successor in
$T_{\nu_1} \cap T_{\nu_2}$.
\nl
Hence
\sn
\item "{$\circledast_2$}"  if $\nu_1 \ne \nu_2 \in {}^\omega 2$ then
lim$(T_{\nu_1}) \cap \text{ lim}(T_{\nu_2})$ is finite. 
\ermn
Now for each $\nu \in {}^\omega 2,N[G_\nu]$ is a $\Bbb Q$-candidate in
$\bold V$ and $N[G_\nu] \models ``q_\nu \in \Bbb Q"$, hence there is
$q^+_\nu \in \Bbb Q$ which is $\langle N[G_\nu],\Bbb Q
\rangle$-generic and $\Bbb Q \models ``q_\nu \le q^+_\nu"$.  Hence really
$q^+_\nu \Vdash_{\Bbb Q} ``\underset\tilde {}\to \eta \in \text{
lim}(T_\nu)"$, so
\mr
\item "{$\circledast_3$}"  if $\nu_1 \ne \nu_2 \in {}^\omega 2$
\ub{then} $q^+_{\nu_1},q^+_{\nu_2}$ are incompatible in $\Bbb Q$.  
\ermn
But this contradicts assumption (a) of \scite{y.1}, i.e., the c.c.c.
${{}}$  \hfill$\square_{\text{Subclaim}}$
\enddemo
\bigskip

\demo{Observation}  Assume
\mr
\item "{$(*)$}"  $p^* \Vdash_{\Bbb P} ``{\underset\tilde {}\to n_i} <
n_{i+1} < \omega,{\underset\tilde {}\to n_i} <
{\underset\tilde {}\to m_i} < \omega$ for every $i < \omega$ and for
every $h \in {}^\omega \omega$ for infinitely many $i$ we have
$h({\underset\tilde {}\to n_i}) < {\underset\tilde {}\to m_i}$.
\ermn
Then we can find $n^* < \omega$ such that for every $m \in [n^*,\omega)$
there is $q$ satisfying $p^* \le_{\Bbb P} q$ and $q \Vdash (\exists
i)({\underset\tilde {}\to n_i} \le n^* \and {\underset\tilde {}\to
m_i} \ge m)$.
\enddemo
\bigskip

\demo{Proof}  We define a function $h:\omega \rightarrow (\omega +1)$
by

$$
\align
h(n) = \sup\{m:&n < m \text{ and for some } q \in \Bbb P \text{ we
have} \\
  &p \le_{\Bbb P} q \text{ and } q \Vdash_{\Bbb P} ``(\exists i)(
{\underset\tilde {}\to n_i} \le n \wedge m \le {\underset\tilde {}\to
m_i})\}.
\endalign
$$
\mn
If for some $n,h(n) = \omega$ we are done.  Otherwise $h \in {}^\omega
\omega$ hence $p^* \Vdash_{\Bbb P} (\exists i)(h({\underset\tilde
{}\to n_i}) < {\underset\tilde {}\to m_i})$ hence there are
$n,m,q,i$ such that $p^* \le_{\Bbb P} q,q \Vdash
``h({\underset\tilde {}\to n_i}) < {\underset\tilde {}\to m_i} \cap
{\underset\tilde {}\to n_i} = n \and {\underset\tilde {}\to m_i} =
m"$.  But by the definition of $q$, $q$ witnesses that $h(n) \ge m$,
contradiction. \nl
${{}}$  \hfill$\square_{\text{observation}}$
\enddemo
\bn
\ub{Continuation of the proof of \scite{y.1}}: 
\bn
The conclusion of the subclaim holds in $\bold V$ as $N \prec ({\Cal
H}(\chi),\in)$, and this gives the conclusion of part (2) of
\scite{y.1} when $\bold t=2$, and the conclusion of part (1) of
\scite{y.1} when $\bold t=1$ is similar to \cite[1.12,p.168]{Sh:480}.
By the subclaim, as $N \prec ({\Cal H}(\chi),\in)$, clearly in $\bold
V^{\Bbb P}$ we have: for every $q^* \in \Bbb Q$ some infinite $u
\subseteq w$ is large for $(q^*,{\underset\tilde {}\to \eta^{\bold
t}})$.  Fix such $q^*,u$.

Let $u \backslash \{0\} = \{n_i:1 \le i < \omega\}$ such that $n_0 =:
0 < n_1 < n_2 < \ldots$, let $\langle k(i,\ell):\ell < \omega \rangle$
be such that $i = \Sigma_\ell k(i,\ell)2^\ell$ where $k(i,\ell) \in \{0,1\}$,
so $k(i,\ell)=0$ when $2^\ell > i$.  Let $\rho^*_m = \langle
k(i,\ell):\ell \le [\log_2(i+1)]\rangle$ where $i = i_u(m)$ is the
unique $i$ such that $n_i \le m < n_{i+1}$.  We define a $\Bbb Q$-name
$\underset\tilde {}\to \rho$ 
(of a member of $({}^\omega 2)^{{\bold V}^{\Bbb Q}})$: let
$\{{\underset\tilde {}\to k_i}:i < \omega\}$ list in increasing order
$\{k < \omega:\eta^{\bold t}(k)=1\}$ and $\underset\tilde {}\to \rho$ 
be $\rho^*_{\underset\tilde {}\to k_0} \char 94 
\rho^*_{\underset\tilde {}\to k_1} \char 94 \rho^*_{\underset\tilde
{}\to k_2} \char 94 \ldots$.

Clearly for every $p \in \Bbb Q$ and $n < \omega$ we have $p \nVdash
``{\underset\tilde {}\to \eta^{\bold t}}(k)=0$ for every $k \ge n"$.  Hence
$\Vdash_{\Bbb Q} 
``\{k < \omega:{\underset\tilde {}\to \eta^{\bold t}}(k)=1\}$ is 
infinite", hence $\Vdash_{\Bbb Q} ``\underset\tilde
{}\to \rho \in {}^{\omega}2"$. \nl
It is enough to prove that $q^* \Vdash_{\Bbb Q} ``\underset\tilde
{}\to \rho$ is a Cohen real over $\bold V^{\Bbb P}"$.  So let $T \in
\bold V^{\Bbb P}$ be a given subtree of ${}^{\omega >}2$ which is
nowhere dense, i.e., $(\forall \eta \in T)(\exists \nu)[\eta
\triangleleft \nu \in {}^{\omega >} 2 \backslash T]$, and we should
prove $q^* \Vdash_{\Bbb Q} ``\underset\tilde {}\to \rho \notin \lim
T"$ so assume $q^* \le q \in \Bbb Q$ 
and we shall find $q',q \le q' \in \Bbb Q$ such that $q'
\Vdash_{\Bbb Q} ``{\underset\tilde {}\to \rho} \notin \lim(T)$, 
this suffices.   We apply
the choice of $u$ so for some $n_* \in u$, if $n,m \in u,n'_* < n <m$
then for some $\nu \in {}^n \omega \cap T_q[{\underset\tilde {}\to
\eta^*}]$ we have $\ell \in (n^*,m) \Rightarrow \nu(\ell)=0$ but
$(\exists \ell)(m \le \ell < n \and {\underset\tilde {}\to
\eta^*}(\ell)=1]$.  As $n_* \in u$ for some $i(*)$ we have
$n_*=n_{i(*)}$ so $\Xi =: \{\rho_m:m < n_*\}$ is finite hence $\Xi' =:
\{\rho^*_{k_0} \char 94 \ldots \char 94 \rho^*_{k_1}:k_0 < \ldots <
k_\ell < n_*\}$ is finite.  As $T$ is nowhere dense we can find a
sequence $\rho^* \in {}^{\omega >}2$ such that: $\rho \in \Xi'
\Rightarrow \rho \char 94 \rho^* \notin T$ and choose $i > i(*)$ such
that $\rho^* \triangleleft \rho_i$.  This is possible by the
definition of $\rho_i$, i.e., it is enough that $i > 2^{\ell
g(\rho^*)}$ and $i = \Sigma\{\rho^*(\ell)2^\ell:\ell < \ell
g(\rho^*)\}$ mod $2^{\ell g(\rho^*)}$.

As said above we can find $q' \ge q$ such that $q' \Vdash$ ``if $n_*
\le \ell < n_i$ then $\underset\tilde {}\to \eta(\ell)=0$ but for some
$\ell \in [n_i,n_{i+1})$ we have $\underset\tilde {}\to \eta(\ell)=1$".
So $q'$ forces that $\rho_{n_i}$ appears in the choice of
$\underset\tilde {}\to \rho$ and before it appears in a concanation of
finite sequences which belong to $\Xi'$, so we are done.
\hfill$\square_{\scite{y.1}}$\margincite{y.1}
\newpage

     \shlhetal 

\newpage
    
REFERENCES.  
\bibliographystyle{lit-plain}
\bibliography{lista,listb,listx,listf,liste}

\def\germ{\frak} \def\scr{\cal} \ifx\documentclass\undefinedcs
  \def\bf{\fam\bffam\tenbf}\def\rm{\fam0\tenrm}\fi 
  \def\defaultdefine#1#2{\expandafter\ifx\csname#1\endcsname\relax
  \expandafter\def\csname#1\endcsname{#2}\fi} \defaultdefine{Bbb}{\bf}
  \defaultdefine{frak}{\bf} \defaultdefine{mathfrak}{\frak}
  \defaultdefine{mathbb}{\bf} \defaultdefine{mathcal}{\cal}
  \defaultdefine{beth}{BETH}\defaultdefine{cal}{\bf} \def\bbfI{{\Bbb I}}
  \def\mbox{\hbox} \def\text{\hbox} \def\om{\omega} \def\Cal#1{{\bf #1}}
  \def\pcf{pcf} \defaultdefine{cf}{cf} \defaultdefine{reals}{{\Bbb R}}
  \defaultdefine{real}{{\Bbb R}} \def\restriction{{|}} \def\club{CLUB}
  \def\w{\omega} \def\exist{\exists} \def\se{{\germ se}} \def\bb{{\bf b}}
  \def\equivalence{\equiv} \let\lt< \let\gt> \def\implies{\Rightarrow}
\begin{thebibliography}{GiSh 357}
\makeatletter \renewcommand{\@biblabel}[1]{[#1]} \makeatother
\def\eprintfootnotetext{References of the form {\tt math.XX/$\cdots$}
 refer to the {\tt xxx.lanl.gov} archive}
\ifx\documentstyle\undefinedcontrolsequence
   \def\anyfootnote{\footnote{*}}
   \else\def\anyfootnote{\footnote}\fi
\def\eprintfn{\ifEprint\anyfootnote{\eprintfootnotetext}\fi\Eprintfalse }
\newif\ifEprint  \Eprinttrue

\bibitem[Bn95]{Bn95}Joerg Brendle.
\newblock {Combinatorial properties of classical forcing notions}.
\newblock {\em Annals of Pure and Applied Logic}, {\bf 73}:143--170, 1995.

\bibitem[GiSh 357]{GiSh:357}Moti Gitik and Saharon Shelah.
\newblock {Forcings with ideals and simple forcing notions}.
\newblock {\em {Israel Journal of Mathematics}}, {\bf 68}:129--160, 1989.

\bibitem[GiSh 412]{GiSh:412}Moti Gitik and Saharon Shelah.
\newblock {More on simple forcing notions and forcings with ideals}.
\newblock {\em {Annals of Pure and Applied Logic}}, {\bf 59}:219--238, 1993.

\bibitem[Ke71]{Ke71}Jerome~H. Keisler.
\newblock {\em {Model theory for infinitary logic. Logic with countable
  conjunctions and finite quantifiers}}, volume~62 of {\em {Studies in Logic
  and the Foundations of Mathematics}}.
\newblock North--Holland Publishing Co., Amsterdam--London, 1971.

\bibitem[Sh 480]{Sh:480}Saharon Shelah.
\newblock {How special are Cohen and random forcings i.e. Boolean algebras of
  the family of subsets of reals modulo meagre or null}.
\newblock {\em {Israel Journal of Mathematics}}, {\bf 88}:159--174, 1994.
\newblock math.LO/9303208.

\bibitem[Sh 666]{Sh:666}Saharon Shelah.
\newblock {On what I do not understand (and have something to say)}.
\newblock {\em {Fundamenta Mathematicae}}, {\bf 166}:1--82, 2000.
\newblock math.LO/9906113.

\end{thebibliography}

\enddocument